\documentclass[11pt]{article}
\usepackage{pdfsync}
\usepackage{hyperref}
\usepackage{array} 
\usepackage{amsfonts}
\usepackage{amsmath}
\usepackage{amssymb}
\usepackage{graphicx}
\usepackage{color}

\newtheorem{theorem}{Theorem}[section]

\newtheorem{definition}[theorem]{Definition}

\numberwithin{equation}{section}

\newcommand{\RR}{\mathbb{R}}

\begin{document}

\title{\textbf{Bistable reaction equations \\
               with doubly nonlinear diffusion} \\[3mm]}

\author{{\Large Alessandro Audrito\footnote{Also affiliated with Universit\`a degli Studi di Torino (Italy) and  Universidad Aut\'{o}noma de Madrid (Spain).} }\\ [4pt]
{\small Dipartimento di Matematica ``Giuseppe Luigi Lagrange'' (DISMA), Politecnico di Torino, Italy.}
}
\date{\vspace{-5ex} }

\maketitle

\centerline{\emph{Dedicated to Professor Juan Luis V\'azquez}}

\begin{abstract}
Reaction-diffusion equations appear in biology and chemistry, and combine linear diffusion with different kind of reaction terms. Some of them are remarkable from the mathematical point of view, since they admit \normalcolor families of travelling waves that describe the asymptotic behaviour of a larger class of solutions $0\leq u(x,t)\leq 1$ of the problem posed in the real line. We investigate here the existence of waves with constant propagation speed, when the linear diffusion is replaced by the ``slow'' doubly nonlinear diffusion. In the present setting we consider bistable reaction terms, which present interesting differences w.r.t. the Fisher-KPP framework recently studied in \cite{AA-JLV:art}. We find different families of travelling waves that are employed to describe the wave propagation of more general solutions and to study the stability/instability of the steady states, even when we extend the study to several space dimensions. A similar study is performed in the critical case that we call ``pseudo-linear'', i.e., when the operator is still nonlinear but has homogeneity one. With respect to the classical model and the ``pseudo-linear'' case, the travelling waves of the ``slow'' diffusion setting exhibit free boundaries.
\\
Finally, as a complement of \cite{AA-JLV:art}, we study the asymptotic behaviour of more general solutions in the presence of a ``heterozygote superior'' reaction function and doubly nonlinear diffusion (``slow'' and ``pseudo-linear'').
\end{abstract}


%
%
%
%
%
%
%
%
%
%
%
\section{Introduction}
In this paper we study the reaction initial-value problem with doubly nonlinear diffusion posed in the whole Euclidean space
\begin{equation}\label{eq:ALLENCAHNPME}
\begin{cases}
\partial_tu = \Delta_p u^m + f(u) \quad &\text{in } \RR^N\times(0,\infty) \\
u(x,0) = u_0(x) \quad &\text{in } \RR^N,
\end{cases}
\end{equation}
where $N \geq 1$, $m > 0$ and $p > 1$. We first discuss the problem of the existence of travelling wave solutions and, later, we use that information to establish the asymptotic behaviour for large times of the solution $u = u(x,t)$ with general initial data and for different ranges of the parameters $m > 0$ and $p > 1$. This work is the natural follow-up of \cite{AA-JLV:art}, where a similar study has been carried out for Fisher-KPP reactions type. As we will see in a moment, the nature of the reaction $f = f(\cdot)$ strongly influences the asymptotic behaviour of the solutions to \eqref{eq:ALLENCAHNPME}. The goal of this paper is to study problem \eqref{eq:ALLENCAHNPME} when the reaction term is not of the Fisher-KPP type, but comes from different biological phenomena. We anticipate that significant differences from the Fisher-KPP setting can be found in both the ODEs analysis (see Theorem \ref{THEOREMEXISTENCEOFTWSPMEREACTIONTYPEC}) and in the asymptotic behaviour of the solutions (see Theorem \ref{ASYMPTOTICBEHAVIOURTHEOREMTYPEC} and \ref{ASYMPTOTICBEHAVIOURTHEOREMTYPECCPRIME}), where ``threshold effects'' and ``non-saturation'' phenomena appear.

We recall that the $p$-Laplacian is a nonlinear operator defined for all $1 \leq p < \infty$ by the formula
\[
\Delta_p v := \nabla\cdot(|\nabla v|^{p-2}\nabla v)
\]
and we consider the more general diffusion term $\Delta_p u^m := \Delta_p(u^m) = \nabla\cdot(|\nabla (u^m)|^{p-2}\nabla (u^m))$, that we call ``doubly nonlinear'' operator since it presents a double power-like nonlinearity. Here,  $\nabla$ is the spatial gradient while $\nabla\cdot$ is the spatial divergence. The doubly nonlinear operator (which can be thought as the composition of the $m$-th power and the $p$-Laplacian) is much used in the elliptic and parabolic literature (see the interesting applications presented in \cite{C-D-D-S-V:art,EstVaz:art,Leib:art}) and allows to recover the Porous Medium operator choosing $p = 2$ or the $p$-Laplacian operator choosing $m = 1$. Of course, choosing $m=1$ and $p = 2$ we obtain the classical Laplacian.

W.r.t. the Porous Medium setting or the $p$-Laplacian one, problem \eqref{eq:ALLENCAHNPME} with doubly nonlinear diffusion is less studied. However, the basic theory of existence, uniqueness and regularity is known. Results about existence of weak solutions of the pure diffusive problem and its generalizations, can be found in the survey \cite{Kal:survey} and the large number of  references therein. The problem of uniqueness was studied later, see for instance \cite{DBen-Her1:art, DBen-Her2:art, Li:art, V2:book, Wu-Yin-Li:art}). For what concerns the regularity, we refer to \cite{V1:book, V2:book} for the Porous Medium setting, while for the $p$-Laplacian case we suggest \cite{DB:book, Lindq:art} and the references therein. Finally, in the doubly nonlinear setting, we refer to \cite{Ivan:art, PorVes:art, Ves:art} and, for the ``pseudo-linear'' case, \cite{Kuusi-Sil-Urb:art}. Finally, we mention \cite{DB:book, V2:book, Wu-Yin-Li:art, Zhao:art} for a proof of the Comparison Principle, which will be an essential technical tool in the proofs of the PDEs part.

\begin{figure}[!ht]
\centering
  \includegraphics[scale=0.5]{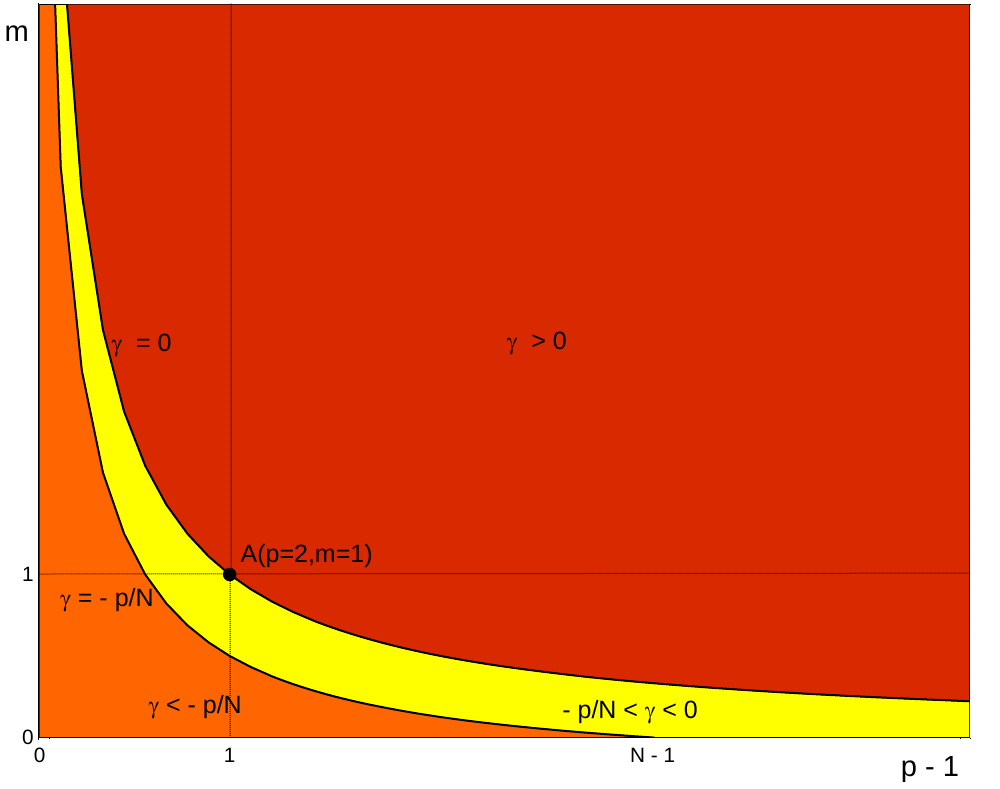}
  \caption{The ``slow diffusion'' area and the ``pseudo-linear'' line in $(m,p-1)$-plane. The yellow and orange area are called ``fast diffusion'' and ``very fast diffusion'' range, respectively, and they will not be studied in this paper.}\label{fig:PARAMETERSMP}
\end{figure}

In order to fix the notations and avoid cumbersome expressions in the rest of the paper, we introduce the constant
\[
\gamma := m(p-1)-1,
\]
which will play an important role in our study. The importance of the constant $\gamma$ is related to the properties of the fundamental solutions of the ``purely diffusive'' doubly nonlinear parabolic equation and we refer the reader to \cite{V1:book}. From the beginning, we consider parameters $m > 0$ and $p > 1$ such that
\[
\gamma \geq 0.
\]
This is an essential restriction. We refer to the assumption $\gamma > 0$ (i.e. $m(p-1)>1$) as the ``slow diffusion'' assumption, while ``pseudo-linear'' assumption when we consider $\gamma = 0$ (i.e. $m(p-1)=1$). Note that $\gamma > 0$ means $m > 1$ if $p = 2$ (Porous Medium case ``slow'' diffusion), while $p > 2$ if $m = 1$ ($p$-Laplacian setting ``slow'' diffusion), i.e., the study of the doubly nonlinear setting covers at the same time, two important models with nonlinear diffusion. Moreover, in the range $\gamma = 0$, we extend the results known in the linear case ($m = 1$ and $p = 2$). In Figure \ref{fig:PARAMETERSMP} the corresponding ranges in the $(m,p-1)$-plane are reported.

The reaction term $f(\cdot)$ is modeled on the function $f(u) = u(1-u)(u-a)$, where $0 < a < 1$ is a fixed parameter and $0 \leq u \leq 1$. More precisely, we assume
\begin{equation}\label{eq:ASSUMPTIONSONTHEREACTIONTERMTYPECPME}
\begin{cases}
f(0) = f(a) = f(1) = 0, \quad f(u) < 0 \text{ in } (0,a), \;\; f(u) > 0 \text{ in } (a,1) \\
f \in C^1([0,1]), \qquad\qquad\quad\; f'(0) < 0, \; f'(a) > 0, \; f'(1) < 0 \\
\int_0^1u^{m-1}f(u)du > 0.
\end{cases}
\end{equation}

Note that the classical reaction $f(u) = u(1-u)(u-a)$ with $0 < a < 1/2 $ satisfies \eqref{eq:ASSUMPTIONSONTHEREACTIONTERMTYPECPME} in the case $m=1$. Furthermore, to complement the work done in \cite{AA-JLV:art} by V\'azquez and the author, we consider also reaction functions satisfying 
\begin{equation}\label{eq:ASSUMPTIONSONTHEREACTIONTERMTYPEDPME}
\begin{cases}
f(0) = f(a) = f(1) = 0, \quad 0 < f(u) \leq f'(0)u \; \text{ in } (0,a), \; f(u) < 0 \text{ in } (a,1)\\
f \in C^1([0,1]), \qquad\qquad\quad\;  f'(0) > 0, \; f'(a)<0, \; f'(1) > 0.
\end{cases}
\end{equation}

We point out that in this second case, the basic model for the reaction is $f(u) = u(1-u)(a-u)$, $0 \leq u \leq 1$ and $0 < a < 1$ is again a fixed parameter. As we will see in a moment, functions which satisfy \eqref{eq:ASSUMPTIONSONTHEREACTIONTERMTYPEDPME} are essentially Fisher-KPP reactions. Since this framework has been already studied in \cite{AA-JLV:art}, we will present the full details of the proofs only in the PDEs part, in which the two models present more significative differences.

\normalcolor

Differently from the reactions of the Fisher-KPP type (or type A) (\cite{Fisher:art,K-P-P:art}), there is not a standard way to indicate them:  bistable reactions \normalcolor, Fitzhugh-Nagumo model or Nagumo's equation in \cite{BoscDamPap2015:art,FifeMcLeod1977:art, McKean1970:art, NagumoArimYoshi1962:art}, ``heterozygote inferior'' reaction in \cite{Aro-Wein1:art}, reaction of type C in \cite{BerestNiren1992:art}, or Allen-Cahn reaction \cite{MatanoPunTes2015:art}, for reaction terms like \eqref{eq:ASSUMPTIONSONTHEREACTIONTERMTYPECPME}. We will refer to them following the notation proposed in \cite{BerestNiren1992:art}, i.e., reaction of type C.
\\
According to the previous choice, we will refer to a function satisfying \eqref{eq:ASSUMPTIONSONTHEREACTIONTERMTYPEDPME} as reaction of type C', even though it was proposed as ``heterozygote superior'' in \cite{Aro-Wein1:art}. It is the least studied of the two models. This is due to the fact that reactions satisfying \eqref{eq:ASSUMPTIONSONTHEREACTIONTERMTYPEDPME} are Fisher-KPP reactions (or reaction of type A) on the interval $[0,a]$, i.e., they satisfy
\begin{equation}\label{eq:TYPEC'RESTRICTECTO0A}
\begin{cases}
f(0) = f(a) = 0, \quad &0 < f(u) \leq f'(0)u, \text{ in } (0,a) \\
f \in C^1([0,a]), \quad &f'(0) > 0, \; f'(a) < 0, \\
\end{cases}
\end{equation}
and so, part of the theory concerning reactions \eqref{eq:ASSUMPTIONSONTHEREACTIONTERMTYPEDPME} is similar to the study of models with Fisher-KPP reactions type. Let us see this fact through a scaling technique. Let us fix $0 < a < 1$ and let us suppose for a moment that $u = u(x,t)$ satisfies the equation
\[
\partial_tu = \Delta_p u^m + f(u) \quad \text{in } \RR^N\times(0,\infty),
\]
where now $f(\cdot)$ is of the Fisher-KPP type (or type A), i.e.
\[
\begin{cases}
f(0) = f(1) = 0, \quad &f(u) > 0, \text{ in } (0,1) \\
f \in C^1([0,1]), \quad &f'(0) > 0, \; f'(1) < 0.
\end{cases}
\]
Then the re-scaled $u_a = u_a(y,s)$ of $u = u(x,t)$ defined by
\[
u(x,t) = a^{-1}u_a(y,t), \quad \text{ with } \; y = a^{\gamma/p}x,
\]
satisfies the equation
\[
\partial_t u_a = \Delta_p u_a^m + f_a(u_a) \quad \text{in } \RR^N\times(0,\infty),
\]
where $f_a(u_a) := af(a^{-1}u_a)$ is of type C' in $[0,a]$, i.e., it satisfies \eqref{eq:TYPEC'RESTRICTECTO0A} with $f_a'(a) = f'(1)$. This property will be very helpful both in the ODEs and PDEs analysis, where we will highlight the connections and the significant differences between the type C' setting and the Fisher-KPP one.

\begin{figure}[!ht]
\centering
  \includegraphics[scale=0.4]{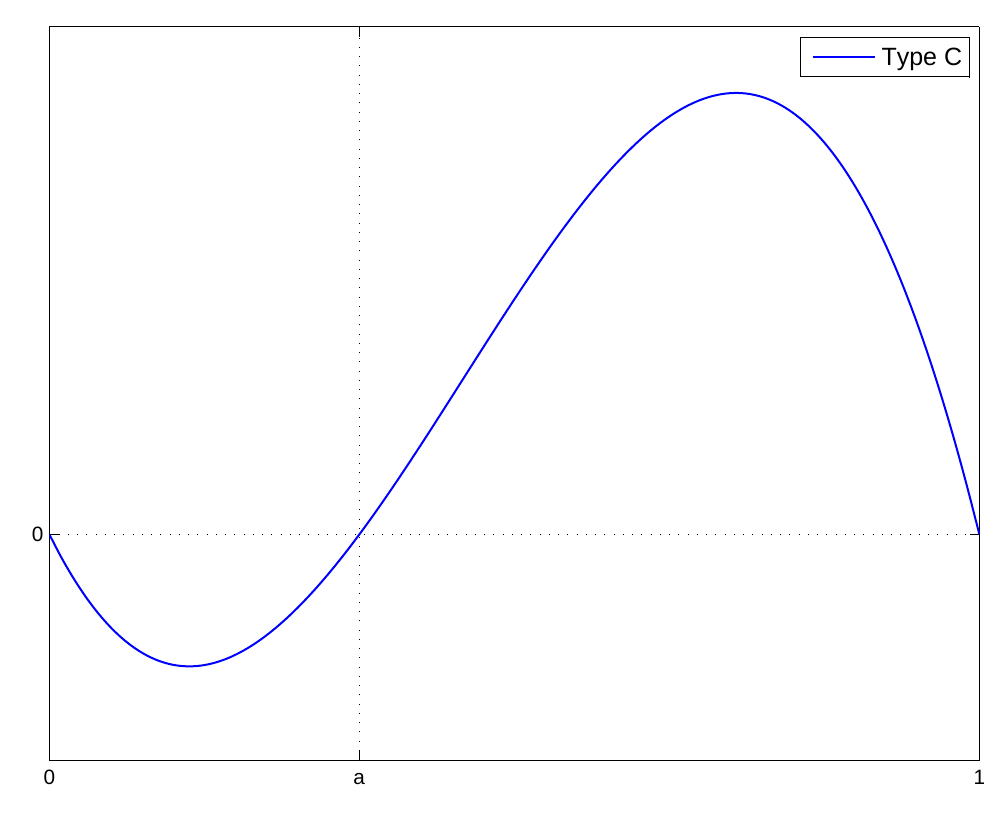}
  \includegraphics[scale=0.4]{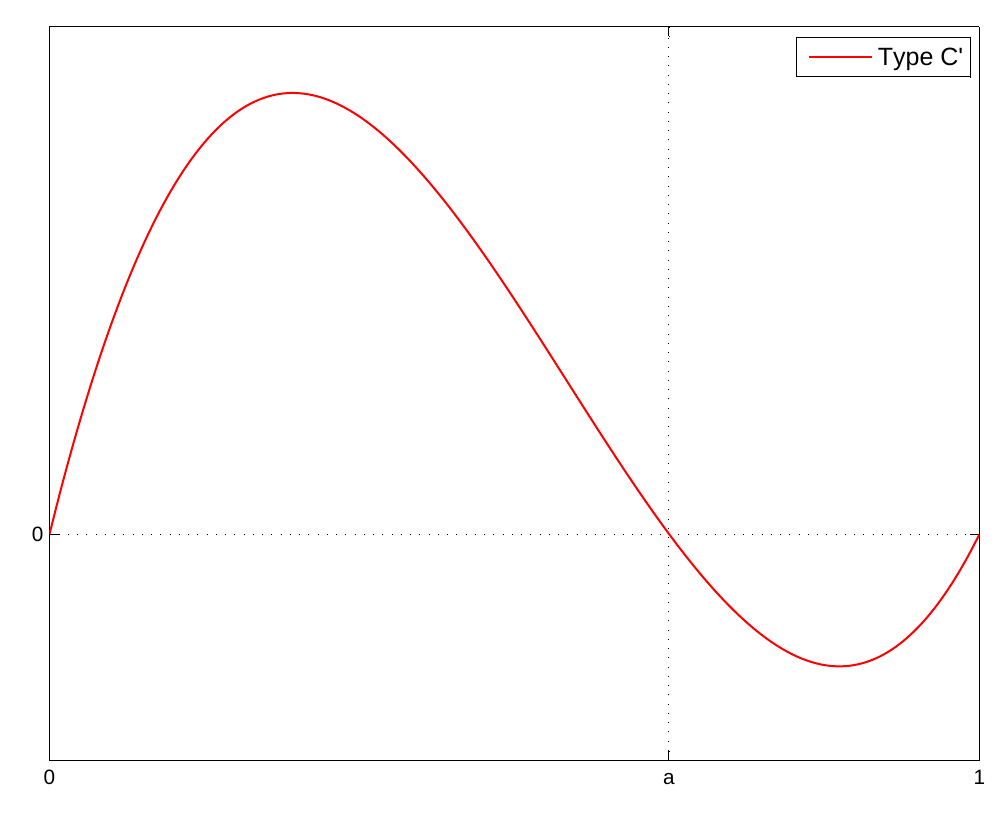}
  \caption{Qualitative representation of the reactions of type C and type C', respectively.}\label{fig:TYPECCPRIME}
\end{figure}

Finally, typical assumptions on the initial datum are
\begin{equation}\label{eq:ASSUMPTIONSONTHEINITIALDATUMBISTABLE}
\begin{cases}
u_0 : \RR^N \to \RR \text{ is continuous with compact support: } u_0 \in \mathcal{C}_c(\RR^N) \\
u_0 \not \equiv 0 \; \text{ and } \; 0 \leq u_0 \leq 1.
\end{cases}
\end{equation}
We point out that, thanks to the Comparison Principle, the assumption $0 \leq u_0 \leq 1$ implies that the solution $u = u(x,t)$ of problem \eqref{eq:ALLENCAHNPME} with reaction \eqref{eq:ASSUMPTIONSONTHEREACTIONTERMTYPECPME} or \eqref{eq:ASSUMPTIONSONTHEREACTIONTERMTYPEDPME} and initial datum \eqref{eq:ASSUMPTIONSONTHEINITIALDATUMBISTABLE} satisfies $0 \leq u \leq 1$ in $\RR^N\times(0,\infty)$. This property has remarkable consequences. First of all it introduces the main goal of this paper, which is studying the stability/instabilty of the steady state $u = 0$, $u = a$, and $u = 1$ of the equation in \eqref{eq:ALLENCAHNPME}, and the rates of convergence of general solutions $u = u(x,t)$ to these constant solutions. Secondly, the restriction $0 \leq u \leq 1$ makes sense from physical viewpoint, since $u = u(x,t)$ stands for the density of a substance evolving in time through the space, according to the nature of the reaction, see once more \cite{Aro-Wein1:art}.

%
%
%
%
%
%
%
%
%
\subsection{Travelling Waves}\label{SUBSECTIONDEFINITIONOFTRAVELLINGWAVES}
They are special solutions with remarkable applications, and there is a huge mathematical literature devoted to them. Let us review the main concepts and definitions.

\noindent Fix $m > 0$ and $p > 1$ such that $\gamma \geq 0$, and assume that we are in space dimension 1 (note that when $N = 1$, the DNL operator has the simpler expression $\Delta_p u^m = \partial_x\left(|\partial_xu^m|^{p-2}\partial_xu^m\right)$. A TW solution to
\begin{equation}\label{eq:REACDIFFTYPCDIM1}
\partial_t u = \partial_x\left(|\partial_xu^m|^{p-2}\partial_xu^m\right) + f(u) \quad \text{in } \RR\times[0,\infty),
\end{equation}
is a solution of the form $u(x,t) = \varphi(\xi)$, where $\xi = x - ct$, $c > 0$ and the \emph{profile} $\varphi(\cdot)$ is a real function. In our reaction-diffusion setting, we will need the profile to satisfy
\begin{equation}\label{eq:CONDITIONONPHIADMISSIBLETWINTRO}
0 \leq \varphi \leq a, \quad \varphi(-\infty) = a, \; \varphi(\infty) = 0 \quad \text{and} \quad \varphi' \leq 0,
\end{equation}
for some $0 < a \leq 1$. In the case in which $a = 1$ we say that $u(x,t) = \varphi(\xi)$ is an \emph{admissible} TW solution, whilst if $0 < a < 1$, we will talk about \emph{a-admissible} TW solution. Depending on the reaction term, these to classes of wave solutions play a role in the PDEs analysis or not. If $f(\cdot)$ is of bistable type we will look for \emph{admissible} TW solutions, while if it is of type C' we will look for \emph{a-admissible} TWs. \normalcolor
\\
Similarly, one can consider TWs of the form $u(x,t) = \varphi(\xi)$ with $\xi = x + ct$, $\varphi$ nondecreasing and such that $\varphi(-\infty) = 0$ and $\varphi(\infty) = a$. It is easy to see that these two options are equivalent, since the profile of the second one can be obtained by reflection of the first one, and it moves in the opposite direction of propagation. Even though in the rest of the paper we will prevalently focus on the first kind of \emph{admissible}/\emph{a-admissible} \eqref{eq:CONDITIONONPHIADMISSIBLETWINTRO}, the ``reflected'' TWs will play an important role in the PDEs part, too. Moreover, note that by definition of the \emph{mobile coordinate} $\xi = x - ct$, if $u(x,t) = \varphi(\xi)$, then $u(x,t+\tau) = u(x-c\tau,t)$ for any $\tau > 0$, which simply means that a TW solution is determined up to ``horizontal displacement'' or, in other words, ``temporal shift''. This property will be essential in both the ODEs and PDEs part (cfr. with Theorem \ref{THEOREMEXISTENCEOFTWSPMEREACTIONTYPEC}, and the proofs of Theorem \ref{ASYMPTOTICBEHAVIOURTHEOREMTYPEC} and Theorem \ref{ASYMPTOTICBEHAVIOURTHEOREMTYPECCPRIME}).

\begin{figure}[!ht]
 \centering
  \includegraphics[scale =0.4]{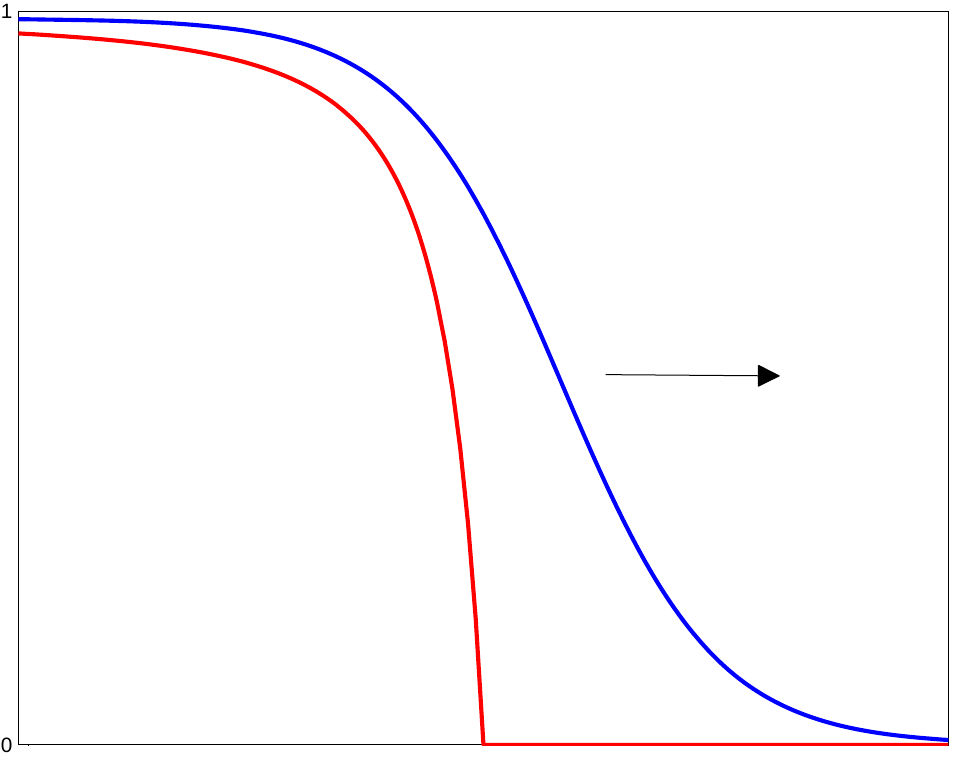}
  \caption{Examples of admissible TWs: Finite and Positive types}\label{fig:ADMISSIBLEANDREFLECTEDFINITETW}
\end{figure}

\noindent Finally, an \emph{admissible}/\emph{a-admissible} TW is said \emph{finite} (or \emph{sharp}) if $\varphi(\xi) = 0$ for $\xi \geq \xi_0$ and/or $\varphi(\xi) = 1$ for $\xi \leq \xi_1$, or \emph{positive} if $\varphi(\xi) > 0$, for all $\xi \in \RR$. The line $x = \xi_0 + ct$ that separates the regions of positivity and vanishing of $u(x,t)$ is then called the \emph{free boundary}. Same name would be given to the line $x = \xi_1 + ct$ and $\varphi(\xi) = 1$ for $\xi \leq \xi_1$ with  $\xi_1$ finite, but this last situation will not happen. Before moving forward, let us mention that travelling fronts with free boundaries were already found in literature in the Fisher-KPP setting and Porous Medium diffusion (see for instance the work of Aronson \cite{Aronson1:art}, De Pablo and V\'azquez \cite{DePablo-Vazquez1:art}, Sanchez-Gardu\~{n}o and Maini \cite{SanchezGardunoMaini1994:art} and the more recent \cite{AA-JLV:art}). See also the paper of Jin et al. \cite{Jin-Yin-Zheng2013:art} concerning bistable reactions with time delay and Porous Medium diffusion.

%
%
%
%
%
%
%
%
%
\subsection{Main results}
The paper is divided in sections as follows:

 Section \ref{EXISTENCEOFTWSTYPECDPME} is devoted to the study of the existence/non-existence of \emph{admissible}/\emph{a-admissible} TW solutions for \eqref{eq:REACDIFFTYPCDIM1}:
\[
\partial_t u = \partial_x\left(|\partial_xu^m|^{p-2}\partial_xu^m\right) + f(u) \quad \text{in } \RR\times(0,\infty),
\]
with reaction term $f(\cdot)$ satisfying \eqref{eq:ASSUMPTIONSONTHEREACTIONTERMTYPECPME} and/or \eqref{eq:ASSUMPTIONSONTHEREACTIONTERMTYPEDPME}. This will be done through a fine analysis of the ODE
\[
c\varphi' + \left(|(\varphi^m)'|^{p-2}(\varphi^m)'\right)' + f(\varphi) = 0 \quad \text{in } \RR,
\]
which is the equation of the profile of a wave solution $u(x,t) = \varphi(x-ct)$ to \eqref{eq:REACDIFFTYPCDIM1} (here $\varphi'$ denotes the derivative of $\varphi$ w.r.t. $\xi = x-ct$). \normalcolor
The following theorem precisely states for which speed/speeds of propagation and reactions terms, equation \eqref{eq:REACDIFFTYPCDIM1} possesses \emph{admissible/a-admissible} TWs, and gives meaningful information on the qualitative shape of these special solutions.
\begin{theorem}\label{THEOREMEXISTENCEOFTWSPMEREACTIONTYPEC}
Fix $N = 1$, $m > 0$, and $p > 1$.

(i) If $\gamma > 0$ and the reaction $f(\cdot)$ is of type C, i.e., it satisfies \eqref{eq:ASSUMPTIONSONTHEREACTIONTERMTYPECPME}, then there exists a unique $c_{\ast} = c_{\ast}(m,p,f) > 0$ such that equation \eqref{eq:REACDIFFTYPCDIM1} possesses a unique admissible TW for $c = c_{\ast}$ and does not have admissible TWs for $0 \leq c \not= c_{\ast}$. Moreover, the TW corresponding to the value $c = c_{\ast}$ is finite (it vanishes in an infinite half-line).

\noindent If $\gamma = 0$, the same conclusions hold except for the fact that the TW corresponding to $c_{\ast}$ is positive everywhere.

(ii) If $\gamma > 0$ and the reaction $f(\cdot)$ is of type C', i.e., it satisfies \eqref{eq:ASSUMPTIONSONTHEREACTIONTERMTYPEDPME}, then there exists a unique $c_{\ast} = c_{\ast}(m,p,f) > 0$ such that equation \eqref{eq:REACDIFFTYPCDIM1} possesses a unique a-admissible TW for all $c \geq c_{\ast}$ and does not have a-admissible TWs for $0 < c < c_{\ast}$. The TWs corresponding to values $c > c_{\ast}$ are positive everywhere while, the TW corresponding to the value $c = c_{\ast}$ is finite.

\noindent Again, if $\gamma = 0$, the same conclusions hold except for the fact that the TW corresponding to $c_{\ast}$ is positive everywhere.

Finally, in both part (i) and (ii), the uniqueness of the TW is understood up to reflection and horizontal displacement (cfr. with the definition of TW solutions, Subsection \ref{SUBSECTIONDEFINITIONOFTRAVELLINGWAVES}).
\end{theorem}
The existence/non-existence of travelling wave solutions for reaction-diffusion equations has been widely studied and still nowadays it is an important field of research. Due to this fact, a bibliographical survey is now in order. In the linear setting ($m = 1$ and $p = 2$), a version of Theorem \ref{THEOREMEXISTENCEOFTWSPMEREACTIONTYPEC} was proved by Aronson and Weinberger in \cite{Aro-Wein1:art,Aro-Wein2:art}, and by Fife and McLeod in \cite{FifeMcLeod1977:art}. Before these works, wave fronts had been studied by McKean in \cite{McKean1970:art}. We have generalized it to the all range $\gamma = 0$ and extended it to the range $\gamma > 0$, where it is proved the existence of \emph{finite} TWs and \emph{free boundaries}, which are the fundamental novelties respect to the classical case.

\noindent Passing to the nonlinear diffusion setting, the existence of \emph{free boundaries} was already observed in Porous Medium setting ($m > 1$ and $p = 2$) in \cite{DePablo-Vazquez1:art} for Fisher-KPP reactions and only more recently in \cite{Laister-Peaplow-Beardmore2004:art} and, later, in \cite{Jin-Yin-Zheng2013:art} for reactions of type C with time delay. See also the recent preprint \cite{Garriz2018:art} where the author studies the long time behaviour of solutions to a Porous Medium reaction-diffusion equation with a larger class of reaction functions. Part (i) of Theorem \ref{THEOREMEXISTENCEOFTWSPMEREACTIONTYPEC} extends the results of \cite{Laister-Peaplow-Beardmore2004:art,Jin-Yin-Zheng2013:art} to the doubly nonlinear setting with reaction satisfying \eqref{eq:ASSUMPTIONSONTHEREACTIONTERMTYPECPME} (here we do not consider reactions with time delay). For reactions of the Fisher-KPP type and nonlinear diffusion we quote \cite{DePablo-Vazquez1:art} for the Porous Medium setting, \cite{Eng-Gav-San:art} and the more recent \cite{GarrioneStrani:art} for the $p$-Laplacian framework and possible generalization and, finally, \cite{AA-JLV:art} for the ``slow'' diffusion range, while \cite{AA-JLV:art1} for the ``fast'' diffusion one.  Actually, part (ii) has been essentially proved in \cite{AA-JLV:art} (in the Fisher-KPP setting one looks for \emph{admissible} TWs instead of \emph{a-admissible} ones). We present a very short sketch of the proof for completeness. \normalcolor

\noindent As mentioned above, TW solutions appear in other kind of reaction-diffusion equations. We mention the fundamental works of \cite{BerestHamel2002:art,BerestHamelNadir2004:art,BerestHamelNadir2009:art} for reactions equations in non homogeneous media, \cite{Alfaro-Coville-Raoul2013:art,BerestRodrig2017:art,Gourley2000:art} for equations with linear diffusion and ``non-local reactions'', whilst \cite{AchleitnerKuehn2015:art,C-R2:art,GuiHuan2015:art,MelletRoqueSire2014:art} for reaction equations with ``non-local'' diffusion of Fractional Laplacian type and \cite{S-V:art} with ``non-local and nonlinear'' diffusion.

In Section \ref{THRESHOLDREACTIONSOFTYPEC} the PDEs part begins. We study the so called ``threshold properties'' and the asymptotic behaviour of radial solutions to problem \eqref{eq:ALLENCAHNPME}-\eqref{eq:ASSUMPTIONSONTHEREACTIONTERMTYPECPME},  depending on the initial datum \eqref{eq:ASSUMPTIONSONTHEINITIALDATUMBISTABLE}. We prove the following result.
\begin{theorem}\label{ASYMPTOTICBEHAVIOURTHEOREMTYPEC}
Let $m > 0$ and $p > 1$ such that $\gamma \geq 0$, and let $N \geq 1$.

\noindent Let $u = u(x,t)$ a radial solution to problem \eqref{eq:ALLENCAHNPME} with reaction of type C (satisfying \eqref{eq:ASSUMPTIONSONTHEREACTIONTERMTYPECPME}). Then:

\noindent (i) There are initial data satisfying \eqref{eq:ASSUMPTIONSONTHEINITIALDATUMBISTABLE} such that
\[
u(x,t) \to 0 \text{ pointwise in } \RR^N, \quad \text{ as } t \to +\infty.
\]
\noindent (ii) There are initial data satisfying \eqref{eq:ASSUMPTIONSONTHEINITIALDATUMBISTABLE} such that
\[
u(x,t) \to 1 \text{ pointwise in } \RR^N, \quad \text{ as } t \to +\infty.
\]
\noindent (iii) Asymptotic behaviour:

\noindent $\bullet$ For all radially decreasing initial data satisfying \eqref{eq:ASSUMPTIONSONTHEINITIALDATUMBISTABLE} and for all $c > c_{\ast}$ it holds
\[
u(x,t) \to 0 \text{ uniformly in } \{|x| \geq ct\}, \quad \text{ as } t \to +\infty.
\]
 Moreover, in the ``slow'' diffusion range $\gamma > 0$, for all $c > c_{\ast}$, there exists a time $\overline{t} > 0$ such that $u(x,t) = 0$ in $\{|x| \geq ct\}$ for all $t \geq \overline{t}$. \normalcolor

\noindent $\bullet$ For the same class of initial data of part (ii) and for all $0 < c < c_{\ast}(m,p,f)$, it holds
\[
u(x,t) \to 1 \text{ uniformly in } \{|x| \leq ct\}, \quad \text{ as } t \to +\infty.
\]
Here $c_{\ast} = c_{\ast}(m,p,f)$ is the critical speed found in Theorem \ref{THEOREMEXISTENCEOFTWSPMEREACTIONTYPEC}, part (i).
\end{theorem}
The previous statement is very significant in terms of stability/instability of the steady states $u = 0$, $u = a$, and $u = 1$, since it explains that both $u = 0$ and $u = 1$ are ``attractors'' (part (i) and (ii)) for the space of nontrivial initial data $u_0 \in \mathcal{C}_c(\RR^N)$, $0 \leq u_0 \leq 1$. This is an important difference respect to the Fisher-KPP setting, where the steady state $u = 1$ is globally stable, whilst $u = 0$ is unstable (cfr. with Theorem 2.6 of \cite{AA-JLV:art}). We ask the reader to note the part (ii) not only asserts that $u = 1$ is an ``attractor'' for a suitable class of initial data, but also gives the rate of convergence $c_{\ast} = c_{\ast}(m,p,f)$ of the solutions to this steady state, for large times. The precise classes of initial data in part (i) and (ii) will be given later (cfr. with Definition \ref{NONREACTIONGINITIALDATA} and Definition \ref{REACTIONGINITIALDATA}).
\\
Even threshold properties of reaction diffusion equations have been largely investigated since the first results proved in \cite{Aro-Wein2:art}. We quote the quite recent works \cite{DuMatano2010:art,MuratovZhong2017:art,Polacik2011:art} for the proof of sharp threshold theorems in the case of linear diffusion. As the reader can see, Theorem \ref{ASYMPTOTICBEHAVIOURTHEOREMTYPEC} is not sharp, but we will see how some special kind of TW solutions found in the fine ODEs analysis carried out in Section \ref{EXISTENCEOFTWSTYPECDPME} can be employed as barriers to show the existence of a threshold effect, which is known in the linear setting but not in the nonlinear one. We stress that, at least to our knowledge, in the case of nonlinear or non-local diffusion, sharp threshold results are not known.
\\

We finally mention that statement (i) (of Theorem \ref{ASYMPTOTICBEHAVIOURTHEOREMTYPEC} of course) is almost immediate if we take $0 \leq u_0 \leq a$ (it easily follows since $0 \leq u(x,t) \leq a$ for any $t > 0$ by comparison, and so $f(u) \leq 0$). Using TWs, we will also prove that there are initial data satisfying \eqref{eq:ASSUMPTIONSONTHEINITIALDATUMBISTABLE}, but not $u_0 \leq a$, such that statement (i) holds true. \normalcolor

In Section \ref{REACTIONSOFTYPECPRIME} we prove the second PDEs result, stated in the following theorem.
\begin{theorem}\label{ASYMPTOTICBEHAVIOURTHEOREMTYPECCPRIME}
Let $m > 0$ and $p > 1$ such that $\gamma \geq 0$, and let $N \geq 1$.

\noindent Let $u = u(x,t)$ a radial solution to problem \eqref{eq:ALLENCAHNPME} with radially decreasing initial datum \eqref{eq:ASSUMPTIONSONTHEINITIALDATUMBISTABLE} and reaction of type C' (satisfying \eqref{eq:ASSUMPTIONSONTHEREACTIONTERMTYPEDPME}). Then:

\noindent For all $0 < c < c_{\ast}$,
\[
u(x,t) \to a \text{ uniformly in } \{|x| \leq ct\}, \quad \text{ as } t \to +\infty.
\]
\noindent For all $c > c_{\ast}$,
\[
u(x,t) \to 0 \text{ uniformly in } \{|x| \geq ct\}, \quad \text{ as } t \to +\infty.
\]
where $c_{\ast} = c_{\ast}(m,p,f)$ is the critical speed found in Theorem \ref{THEOREMEXISTENCEOFTWSPMEREACTIONTYPEC}, Part (ii).  Again, in the ``slow'' diffusion range $\gamma > 0$, for all $c > c_{\ast}$, there exists a time $\overline{t} > 0$ such that $u(x,t) = 0$ in $\{|x| \geq ct\}$ for all $t \geq \overline{t}$. \normalcolor
\end{theorem}
Even in this setting, the previous theorem gives relevant information on the stability/instability of the steady states $u = 0$, $u = a$ and $u = 1$. Possibly, the most important one is that the state $u = a$ is globally stable w.r.t. the class of initial data $u_0 \in \mathcal{C}_c(\RR^N)$, $0 \leq u_0 \leq 1$, whilst both $u = 0$ and $u = 1$ are unstable. This is a strong departure from the previous case of reaction of Type C and of the Fisher-KPP type. Furthermore, as in Theorem \ref{THEOREMEXISTENCEOFTWSPMEREACTIONTYPEC} (Part (ii) and (iii)), it is shown that the rate of (uniform) convergence to the stable steady state is approximately constant for large times and it coincides with the critical speed of propagation $c_{\ast} = c_{\ast}(m,p,f)$ found in the ODEs analysis.

\noindent Theorem \ref{ASYMPTOTICBEHAVIOURTHEOREMTYPECCPRIME} was known for Fisher-KPP reactions and $a=1$ (see \cite{AA-JLV:art} and the references therein) and was proved for the linear case in \cite{Aro-Wein2:art}, together with a so called ``hair-trigger effect'' result that we do not study in this paper.  Finally, let us stress that, as in the ODEs part, our methods relies on the proof of Theorem 2.6 of \cite{AA-JLV:art}. The main difference w.r.t. to that framework is to prove that initial data $0 \leq u_0 \leq 1$ (not necessarily $0 \leq u_0 \leq a$) generate solutions that converge to the steady state $u = a$ for large times. \normalcolor

\paragraph{Remarks.} First of all, we note that in the statements of both Theorem \ref{ASYMPTOTICBEHAVIOURTHEOREMTYPEC} and \ref{ASYMPTOTICBEHAVIOURTHEOREMTYPECCPRIME}, when the spatial dimension is $N=1$, the initial data are not needed to be radially decreasing (this fact will be clarified later, in the proofs of the above theorems).
\\
Secondly, in order to simplify the reading, we have decided to state Theorem \ref{ASYMPTOTICBEHAVIOURTHEOREMTYPEC} and \ref{ASYMPTOTICBEHAVIOURTHEOREMTYPECCPRIME} for radial solutions to problem \ref{eq:ALLENCAHNPME} (generated by radially decreasing initial data). A simple comparison with ``sub'' and ``super'' initial data shows that these theorems hold true for initial data satisfying \eqref{eq:ASSUMPTIONSONTHEINITIALDATUMBISTABLE}. Indeed, if $u_0 = u_0(x)$ satisfies \eqref{eq:ASSUMPTIONSONTHEINITIALDATUMBISTABLE}, there are $\underline{u}_0 = \underline{u}_0(|x|)$ and $\overline{u}_0 = \overline{u}_0(|x|)$ radially decreasing satisfying \eqref{eq:ASSUMPTIONSONTHEINITIALDATUMBISTABLE} such that $\underline{u}_0 \leq u_0 \leq \overline{u}_0$ in $\RR^N$. Consequently, if $\underline{u} = \underline{u}(x,t)$ and $\overline{u} = \overline{u}(x,t)$ are radial solutions to problem \eqref{eq:ALLENCAHNPME} with initial data $\underline{u}_0$ and $\overline{u}_0$, respectively, it follows $\underline{u}(x,t) \leq u(x,t) \leq \overline{u}(x,t)$ for all $x \in \RR^N$ and $t > 0$, thanks to the comparison principle. So, since Theorem \ref{ASYMPTOTICBEHAVIOURTHEOREMTYPEC} and \eqref{ASYMPTOTICBEHAVIOURTHEOREMTYPECCPRIME} hold for $\underline{u} = \underline{u}(x,t)$ and $\overline{u} = \overline{u}(x,t)$, they will hold for $u = u(x,t)$, too.

%
%
%
%
%
%
%
%
%
\subsection{Preliminaries on doubly nonlinear diffusion}\label{SUBSECTIONPRELIMINARIESTYPECCPRIME}
In this brief subsection we recall some important features about doubly nonlinear diffusion, needed in the PDEs part. In particular, we focus on the so called Barenblatt solutions.
\paragraph{Barenblatt solutions.} Fix $m > 0$ and $p > 1$ such that $\gamma \geq 0$ and consider the ``purely diffusive'' doubly nonlinear problem:
\begin{equation}\label{eq:PARABOLICPLAPLACIANEQUATIONINTRO}
\begin{cases}
\begin{aligned}
\partial_tu = \Delta_p u^m \;\quad &\text{in } \RR^N\times(0,\infty) \\
u(t) \to M\delta_0 \quad\;\, &\text{in } \RR^N \text{ as } t \to 0,
\end{aligned}
\end{cases}
\end{equation}
where $M\delta_0(\cdot)$ is the Dirac's function with mass $M>0$ in the origin of $\RR^N$ and the convergence has to be intended in the sense of measures.
\paragraph{Case $\boldsymbol{\gamma > 0}$.} It has been proved (see \cite{V1:book}) that problem \eqref{eq:PARABOLICPLAPLACIANEQUATIONINTRO} admits continuous weak solutions in self-similar form $B_M(x,t) = t^{-\alpha}F_M(xt^{-\alpha/N})$, called Barenblatt solutions, where the \emph{profile} $F_M(\cdot)$ is defined by the formula:
\[
F_M(\xi) = \Big(C_M - k|\xi|^{\frac{p}{p-1}} \Big)_{+}^{\frac{p-1}{\gamma}}
\]
where
\[
\alpha = \frac{1}{\gamma + p/N}, \quad k = \frac{\gamma}{p}\Big(\frac{\alpha}{N}\Big)^{\frac{1}{p-1}}
\]
and $C_M>0$ is determined in terms of the mass choosing $M = \int_{\RR^N}B_M(x,t)dx$ (see \cite{V1:book} for a complete treatise). We remind the reader that the solution has a \emph{free boundary} which separates the set in which the solution is positive from the set in which it is identically zero (``slow'' diffusion case).
\paragraph{Case $\boldsymbol{\gamma = 0}$.} Again we have Barenblatt solutions in self-similar form. The new profile can be obtained passing to the limit as $\gamma \to 0$:
\[
F_M(\xi) = C_M \exp \big(-k|\xi|^{\frac{p}{p-1}} \big),
\]
where $C_M > 0$ is a free parameter and it is determined fixing the mass, while now $k = (p-1)p^{-p/(p-1)}$. Note that, in this case the constant $\alpha = N/p$ and for the values $m=1$ and $p=2$, we have $\alpha = N/2$ and $F_M(\cdot)$ is the Gaussian profile. The main difference with the case $\gamma > 0$ is that now the Barenblatt solutions have no \emph{free boundary} but are always positive. This fact has repercussions on the shape of the TW solutions. Indeed, we will find finite TWs in the case $\gamma > 0$ whilst positive TWs in the case $\gamma = 0$.
%
%
%
%
%
%
%
%
%
%
\section{Existence of TWs}\label{EXISTENCEOFTWSTYPECDPME}
We consider equation \eqref{eq:REACDIFFTYPCDIM1} (that we rename for the reader's convenience) with reaction satisfying \eqref{eq:ASSUMPTIONSONTHEREACTIONTERMTYPECPME}:
\begin{equation}\label{eq:APPROXIMATEEQUATIONSUPERSOLHYP}
\partial_t u = \partial_x\left(|\partial_xu^m|^{p-2}\partial_xu^m\right) + f(u) \quad \text{in } \RR\times[0,\infty),
\end{equation}
and we look for \emph{admissible} TW solutions $u(x,t) = \varphi(\xi)$, where $\xi = x - ct$, $c > 0$, and $\varphi(\cdot)$ satisfying $0 \leq \varphi \leq 1$, $\varphi(-\infty) = 1$, $\varphi(\infty) = 0$ and $\varphi' \leq 0$. Note that there is a second option in which $\varphi' \geq 0$ and the wave moves in the opposite direction, but we can skip this case since it is obtained from the previous one by reflection.
\paragraph{Proof of Theorem \ref{THEOREMEXISTENCEOFTWSPMEREACTIONTYPEC}: Part (i), range $\boldsymbol{\gamma > 0}$.} Fix $m > 0$ and $p > 1$ such that $\gamma > 0$. Substituting $u(x,t) = \varphi(x - ct)$ in \eqref{eq:APPROXIMATEEQUATIONSUPERSOLHYP}, the equation of the profile reads
\[
-c\varphi' = \left[|\left(\varphi^m\right)'|^{p-2}\left(\varphi^m\right)'\right]' + f(\varphi) \quad \text{in } \RR,
\]
where $\varphi'$ stands for the derivative of $\varphi$ w.r.t. $\xi =x - ct$. Proceeding as in \cite{DePablo-Vazquez1:art} and \cite{AA-JLV:art}, we consider the variables
\begin{equation}\label{eq:NSTANDCVARTWSHYPPME}
X = \varphi \qquad \text{ and } \qquad Z = -\left(\frac{m(p-1)}{\gamma}\varphi^{\frac{\gamma}{p-1}}\right)' = -mX^{\frac{\gamma}{p-1}-1}X'.
\end{equation}
They correspond to the density and the derivative of the pressure profile (see \cite{EstVaz:art} and \cite{V2:book}, Chapter 2). Assuming $X \geq 0$, we obtain the first-order ODEs system
\begin{equation}\label{eq:SYSTEMNONSINGULARTWTYPECPME1}
-m\frac{dX}{d\xi} = X^{1- \frac{\gamma}{p-1}}Z, \quad\quad -m(p-1)X^{\frac{\gamma}{p-1}}|Z|^{p-2} \frac{dZ}{d\xi} = cZ - |Z|^p - mX^{\frac{\gamma}{p-1}-1}f(X),
\end{equation}
that we re-write as the non-singular system
\begin{equation}\label{eq:SYSTEMNONSINGULARTWSTYPECPME2}
\frac{dX}{d\tau} = (p-1)X|Z|^{p-2}Z, \quad\quad \frac{dZ}{d\tau} = cZ - |Z|^p - mX^{\frac{\gamma}{p-1}-1}f(X),
\end{equation}
where we have used the re-parametrization  $d\xi = -m(p-1)X^{\frac{\gamma}{p-1}}|Z|^{p-2}d\tau$. Systems \eqref{eq:SYSTEMNONSINGULARTWTYPECPME1} and \eqref{eq:SYSTEMNONSINGULARTWSTYPECPME2} are equivalent outside the critical points $O(0,0)$, $S(1,0)$, $A(a,0)$, $R_c(0,c^{1/(p-1)})$, and their trajectories correspond to the solutions to the equation
\begin{equation}\label{eq:EQUATIONOFTHETRAJECTORIESHYP}
\frac{dZ}{dX} = \frac{cZ - |Z|^p - f_{m,p}(X)}{(p-1)X|Z|^{p-2}Z} := H(X,Z;c),
\end{equation}
called \emph{equation of the trajectories}. To simplify the notation, in the previous formula we have introduced the function
\[
f_{m,p}(X) = mX^{\frac{\gamma}{p-1}-1}f(X), \qquad 0 \leq X \leq 1,
\]
with $f_{m,p}(0) = f_{m,p}(a) = f_{m,p}(1) = 0$ and $f_{m,p}(X) < 0$ for $0 < X < a$, while $f_{m,p}(X) > 0$ for $a < X < 1$.

According to the statement of the theorem, we prove the existence of a special speed $c_{\ast} = c_{\ast}(m,p,f)$ with corresponding trajectory linking $S(1,0)$ and $R_{c_{\ast}}(0,c_{\ast}^{1/(p-1)})$ and lying in the strip $[0,1]\times[0,+\infty)$ of the $(X,Z)$-plane. We will show that this connection is the \emph{finite} TW we are looking for.
\\
To do this, we have to understand the qualitative behaviour of the trajectories of system \eqref{eq:SYSTEMNONSINGULARTWSTYPECPME2} (or, equivalently, the solutions of equation \eqref{eq:EQUATIONOFTHETRAJECTORIESHYP}) in dependence of the parameter $c > 0$. This will be done in some steps as follows: in the first one, we consider the simpler case $c = 0$, which is fundamental to exclude the existence of \emph{admissible} TWs for small speeds of propagation. The assumption $\int_0^1 u^{m-1}f(u)du > 0$ (cfr. with \eqref{eq:ASSUMPTIONSONTHEREACTIONTERMTYPECPME}) plays an important role in what follows. Then we study the local behaviour of the trajectories near the critical points and we prove more global monotonicity properties of the trajectories w.r.t. the speed $c>0$. Finally, we employ them to show the existence or non-existence of trajectories linking the critical points $S(1,0)$ and $R_c(0,c^{1/(p-1)})$, which correspond to a \emph{finite} TW (see \emph{Step 4}).

\emph{Step 0: Case $c = 0$.} As we have explained in the previous paragraph, we begin by taking $c = 0$ and we show that for the null speed, there are not \emph{admissible} TW profiles. With this choice, system \eqref{eq:SYSTEMNONSINGULARTWSTYPECPME2} and equation \eqref{eq:EQUATIONOFTHETRAJECTORIESHYP} become
\[
\begin{cases}
\dot{X} = (p-1)X|Z|^{p-2}Z, \\
\dot{Z} = - |Z|^p - f_{m,p}(X),
\end{cases}
\qquad
\text{ and }\;\;
\qquad
\frac{dZ}{dX} = -\frac{|Z|^p + f_{m,p}(X)}{(p-1)X|Z|^{p-2}Z} = H(X,Z;0),
\]
respectively (here $\dot{X}$ means $dX/d\tau$). The critical points are $O(0,0)$, $A(a,0)$, and $S(1,0)$ (note that the point $R_c$ ``collapses'' to $O(0,0)$).

\noindent Respect to the linear case, our system does not conserve the energy along the solutions (see \cite{Aro-Wein1:art}). Consequently, excluding the existence of a trajectory, contained in the strip $(0,1)\times(0,\infty)$ in the $(X,Z)$-plane and linking $O(0,0)$ and $S(1,0)$, is done by studying more qualitative properties of the trajectories in the $(X,Z)$-plane.

\noindent So, we begin by analyzing the \emph{null isoclines} $\widetilde{Z} = \widetilde{Z}(X)$ of our system, i.e. the solutions of the equation:
\[
|\widetilde{Z}|^p + mX^{\frac{\gamma}{p-1}-1}f(X) = 0, \qquad 0 \leq X \leq 1.
\]
They are composed by two branches linking the points $O(0,0)$ and $A(a,0)$, lying in the strip $[0,a]\times(0,\infty)$ and $[0,a]\times(-\infty,0)$, respectively, and they satisfy
\[
\widetilde{Z}(X) \sim \pm\sqrt[p]{-mf'(0)}X^{\frac{\gamma}{p(p-1)}}, \quad  \text{ for } X \sim 0.
\]
 Now, there are two symmetric trajectories: one positive and one negative in a right-neighbourhood of $O(0,0)$, the first ``leaving'' $O(0,0)$ while the second ``entering'' $O(0,0)$ (this follows from study of the \emph{null isoclines} and the sign of the derivative $dZ/dX$ in the $(X,Z)$-plane). \normalcolor Moreover, since $H(X,-Z;0) = -H(X,Z;0)$, the two trajectories coincide and we obtain a unique trajectory linking $O(0,0)$ with itself. Now, let us focus on the part lying in $[0,1)\times[0,\infty)$, $T^+ = T^+(X)$ and let $T_0 = T_0(X)$ be the trajectory ``coming into'' $S(1,0)$. If $T^+ = T^+(X)$ and $T_0 = T_0(X)$ touch at a point, they coincide in $[0,1]$ and the resulting trajectory has the shape of an \emph{admissible} profile. In the next paragraphs, we show that $T^+$ and $T_0$ must be two distinct trajectories and the just described case cannot happen.  Let us stress that the uniqueness of the trajectory $T_0$ is not trivial. It is a consequence of the fact that the equation of $Z = Z(X)$ can be transformed into a first order linear equation for $S(X) := X^{2-\gamma/(p-1)}Z^p(X)$ with smooth coefficients near $X = 1$ (cfr. with equation \eqref{eq:FIRSTORDERODES}). \normalcolor

\begin{figure}[!ht]
  \centering
  \includegraphics[scale =0.41]{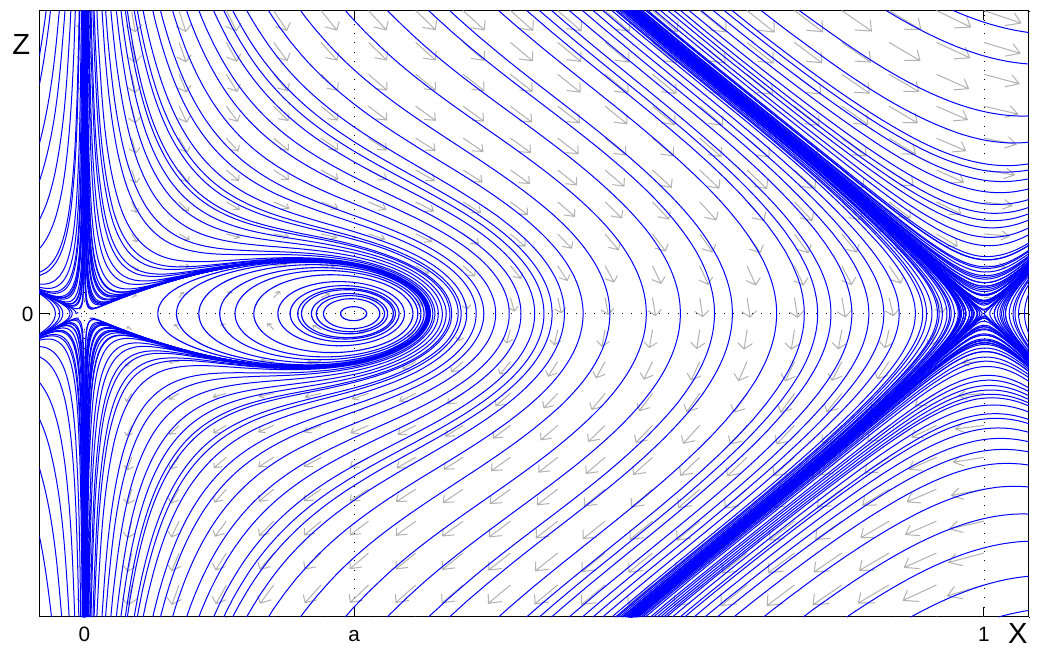}
  \includegraphics[scale =0.408]{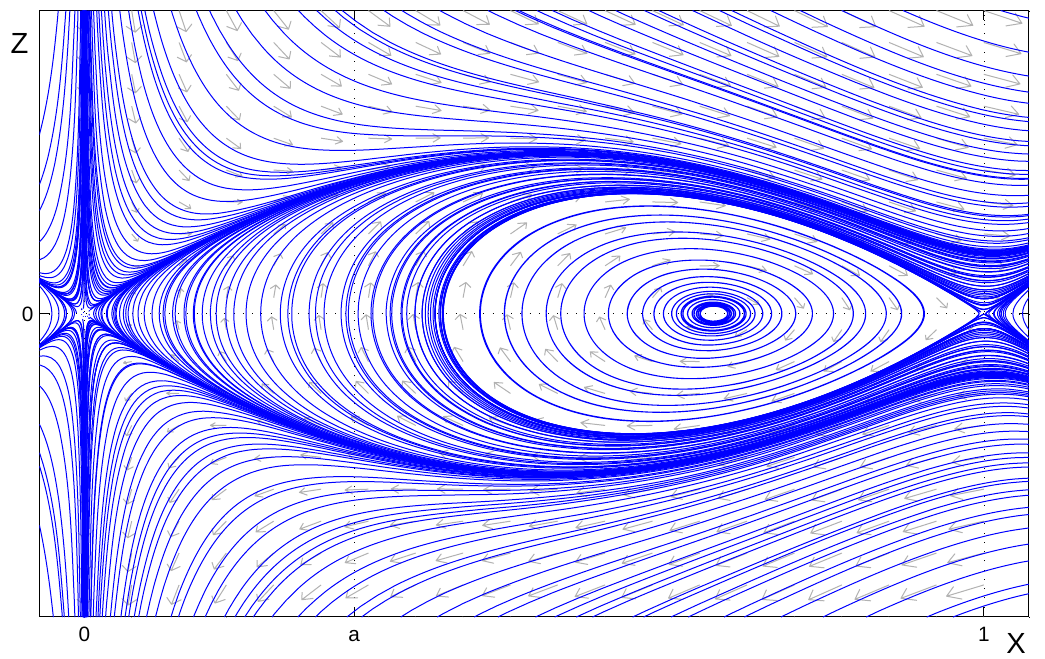}
  \caption{Reactions of type C, range $\gamma > 0$, case $c = 0$. Qualitative behaviour of the trajectories in the $(X,Z)$-plane for $f(u) = u(1-u)(u - a)$, $a = 0.3, 0.7$. The second case is excluded by the assumption $\int_0^1u^{m-1}f(u)du > 0$.}\label{fig:QUALBEHAVTRAJPMETYPEC0}
\end{figure}

\noindent As first observation, since the solution $T^+ = T^+(X)$ stays below the positive branch $\widetilde{Z} = \widetilde{Z}(X)$ for $X \sim 0$, a simple approximation argument shows that
\[
T^+(X) \sim \sqrt[p]{-\frac{mpf'(0)}{\gamma + p}} X^{\frac{\gamma}{p(p-1)}}, \quad \text{ for } X \sim 0.
\]
Hence, substituting it in the first equation of system \eqref{eq:SYSTEMNONSINGULARTWTYPECPME1}, we obtain (up to a multiplicative constant):
\[
-\frac{dX}{d\xi} \sim X^{1- \frac{\gamma}{p-1}}T^+(X) \quad \Leftrightarrow \quad X^{\frac{\gamma}{p}}(\xi) = \varphi^{\frac{\gamma}{p}}(\xi) \sim \xi_0 - \xi, \quad \text{ for } \xi \sim \xi_0^-,
\]
which contradicts then Darcy law of the \emph{free boundary}  (cfr. with \eqref{eq:ADMISSIBLEPROFILECASTDARCYLAWHYP} and \cite{V2:book}, Chapter 4 for the Porous Medium case)
\[
X^{\frac{\gamma}{p-1}}(\xi) = \varphi^{\frac{\gamma}{p-1}}(\xi) \sim \xi_0 - \xi, \qquad \text{for } \xi \sim \xi_0^-.
\]
\normalcolor
Consequently, if $T^+ = T^+(X)$ and $T_0 = T_0(X)$ coincide, we immediately conclude that the resulting trajectory linking $O(0,0)$ and $S(1,0)$ cannot be an \emph{admissible finite} TW and we conclude the non existence of \emph{admissible} TWs for $c = 0$. The qualitative behaviour of the trajectories in the $(X,Z)$-plane is shown in Figure \ref{fig:QUALBEHAVTRAJPMETYPEC0}.

\noindent However, in what follows, we will need to exclude the case in which the trajectory $T_0 = T_0(X)$ ``coming into'' $S(1,0)$ has either a closed curve or $S(1,0)$ as \emph{negative} limit set, or crosses at some point the negative half-line $X=1$ (cfr. with the right picture of Figure \ref{fig:QUALBEHAVTRAJPMETYPEC0}). To achieve this, we will show that $T_0 = T_0(X) \sim +\infty$ as $X \sim 0$, using our initial assumption on the reaction term (see \eqref{eq:ASSUMPTIONSONTHEREACTIONTERMTYPEDPME}) that we rename for convenience:
\begin{equation}\label{eq:NOHOMOCLINICASSUMPTIONONF}
\int_0^1u^{m-1}f(u)\,du > 0.
\end{equation}
For $0 < X \leq 1$ and $Z > 0$, the equation of the trajectories can be re-written as
\[
\frac{dZ}{dX} = -\frac{Z^p + mX^{\frac{\gamma}{p-1}-1}f(X)}{(p-1)XZ^{p-1}} \quad \Leftrightarrow \quad pX^{2-\frac{\gamma}{p-1}}Z^{p-1}\frac{dZ}{dX} = -\frac{pX^{1-\frac{\gamma}{p-1}}Z^p + mpf(X)}{(p-1)}
\]
Using that
\[
\frac{d}{dX}\left( X^{2-\frac{\gamma}{p-1}}Z^p\right) =
\left( 2-\frac{\gamma}{p-1} \right)X^{1-\frac{\gamma}{p-1}} Z^p +
pX^{2-\frac{\gamma}{p-1}}Z^{p-1}\frac{dZ}{dX},
\]
and the previous equation, we deduce that $S(X) := X^{2-\frac{\gamma}{p-1}}Z^p$ satisfies the equation
\begin{equation}\label{eq:FIRSTORDERODES}
\frac{dS}{dX} = \frac{1-m}{X}\,S - \frac{mp}{p-1}f(X), \quad 0 < X \leq 1,
\end{equation}
where we have used the definition of $\gamma := m(p-1)-1$. Now, assume for a moment $m\not=1$. It is simple to integrate the previous equation obtaining
\[
S(X) = X^{1-m}\left[k - \frac{mp}{p-1}\int_0^X u^{m-1}f(u)du\right], \quad 0 < X \leq 1,
\]
where $k$ is a free parameter. Now, coming back to the function $Z = Z(X)$, we get
\begin{equation}\label{eq:EQUATIONFORSXZPLANETYPEC}
Z(X) = X^{-\frac{1}{p-1}}\left[k - \frac{mp}{p-1}\int_0^X u^{m-1}f(u)du\right]^{\frac{1}{p}}, \quad 0 < X \leq 1,
\end{equation}
and, thanks to our assumption \eqref{eq:NOHOMOCLINICASSUMPTIONONF}, we can take
\[
\int_0^1 u^{m-1}f(u)du := \frac{p-1}{mp}\,h > 0.
\]
Furthermore, choosing $k = h$ in \eqref{eq:EQUATIONFORSXZPLANETYPEC}, we deduce $Z(1) = 0$. Hence, if $T_0 = T_0(X)$ is the trajectory ``coming into'' $S(1,0)$, we have by uniqueness of this solution
\[
T_0(X) = X^{-\frac{1}{p-1}}\left[h - \frac{mp}{p-1}\int_0^X u^{m-1}f(u)du\right]^{\frac{1}{p}} \sim +\infty, \quad \text{ as } X \sim 0,
\]
proving our claim (cfr. with the left diagram shown in Figure \ref{fig:QUALBEHAVTRAJPMETYPEC0}). The case $m = 1$ is very similar and formula \eqref{eq:EQUATIONFORSXZPLANETYPEC} holds with $m=1$.  As we mentioned above, the uniqueness of $T_0 = T_0(X)$ follows from the fact that the r.h.s. of \eqref{eq:FIRSTORDERODES} is of class $C^1$ in a neighborhood of $S(1,0)$. \normalcolor

\noindent We end this paragraph pointing out that, thanks to the continuous dependence of the solutions w.r.t. to the parameter $c \geq 0$, we deduce that there are not \emph{admissible} TWs for values of $c > 0$ small enough.

\begin{figure}[!ht]
  \centering
  \includegraphics[scale =0.4]{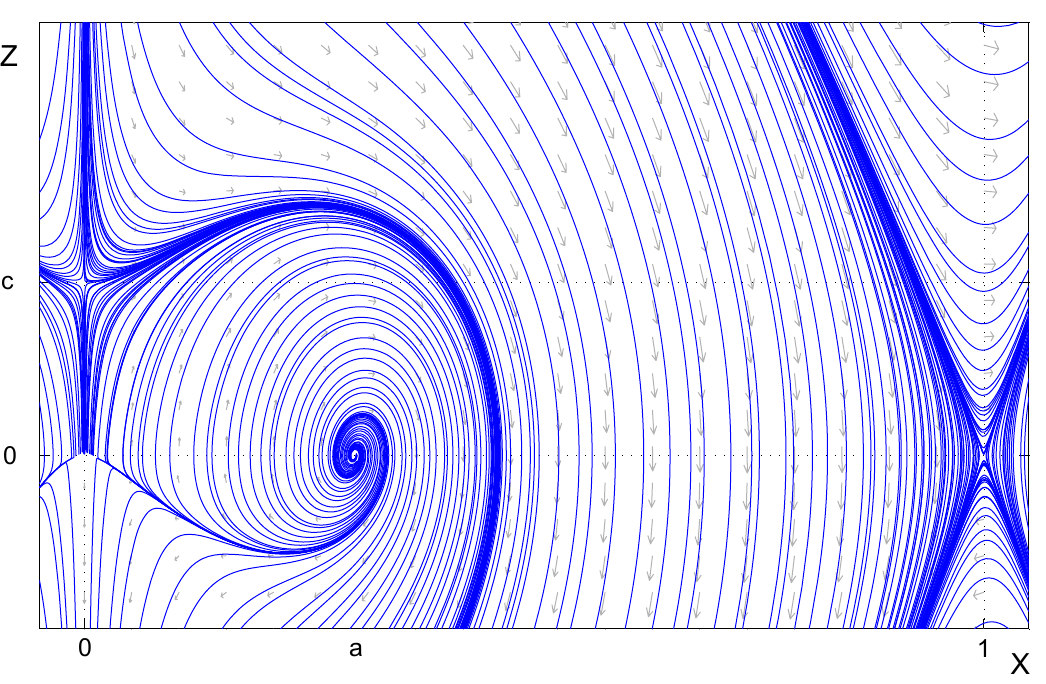}
  \includegraphics[scale =0.4]{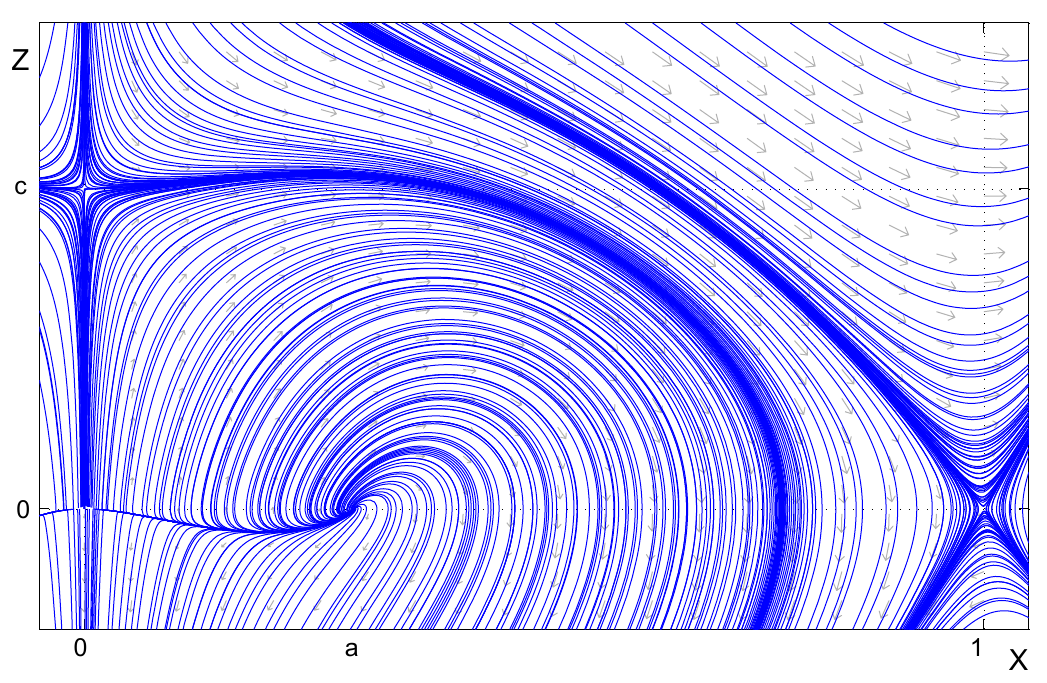}
  \includegraphics[scale =0.4]{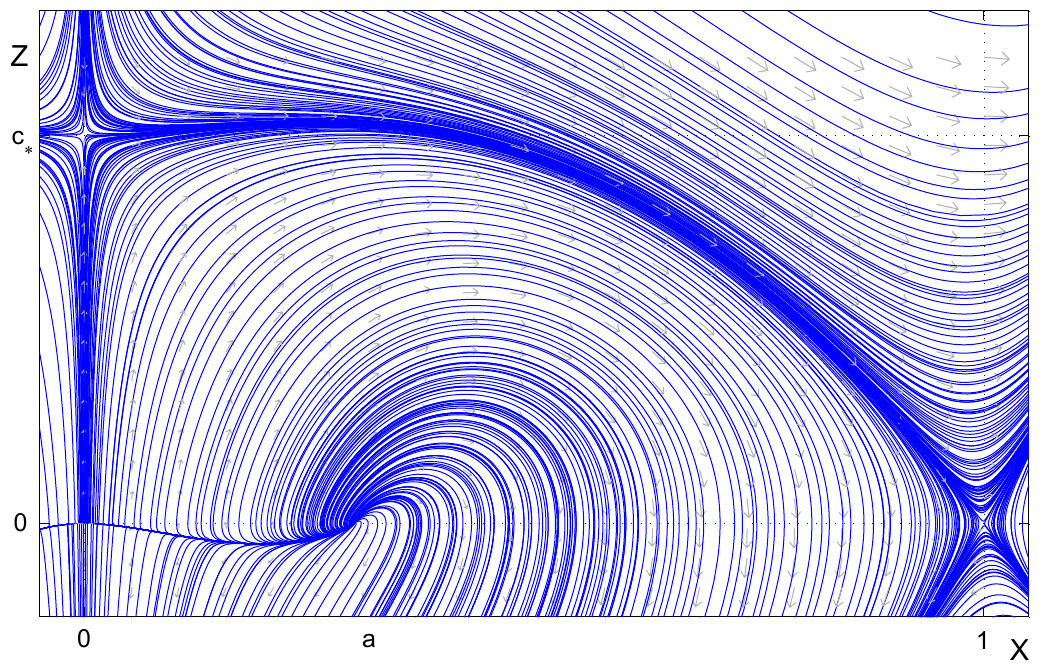}
  \includegraphics[scale =0.4]{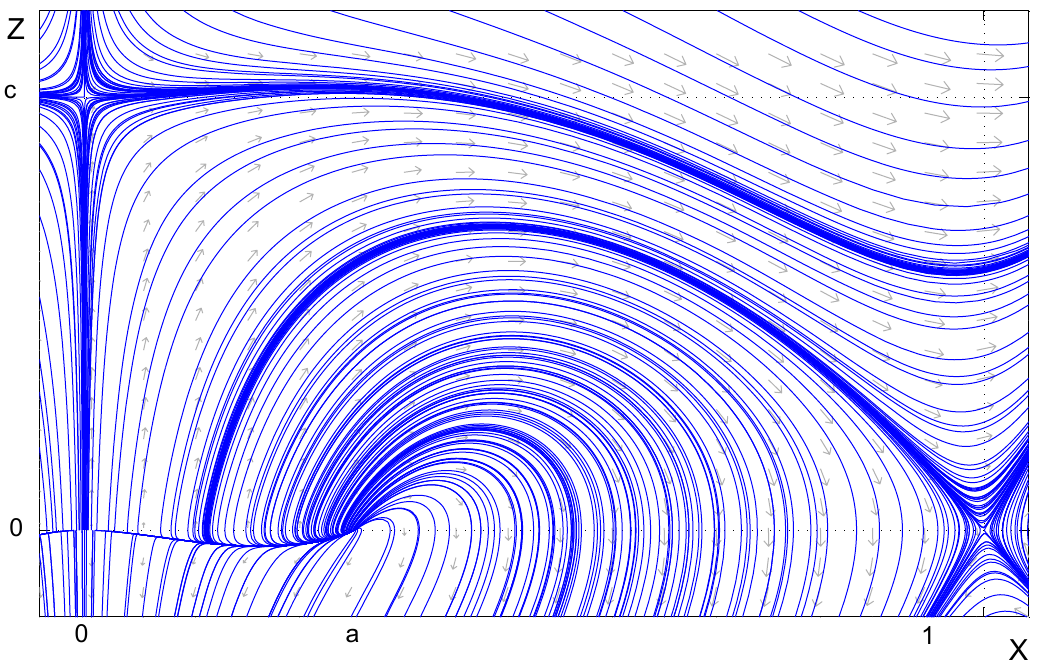}
  \caption{Reactions of type C, range $\gamma > 0$. Qualitative behaviour of the trajectories in the $(X,Z)$-plane for $f(u) = u(1-u)(u - a)$, $a = 0.3$. The first two pictures show the case $0< c < c_{\ast}$, while the others the cases $c = c_{\ast}$ and $c > c_{\ast}$, respectively.}\label{fig:QUALBEHAVTRAJPMETYPEC}
\end{figure}

\emph{Step 1: Local analysis of $S(1,0)$.} From now on, we consider $c > 0$. The local analysis near the point $S(1,0)$ has been carried out in \cite{AA-JLV:art} (see Theorem 2.1 \emph{Step 2}) where the authors proved that there exists a unique trajectory $T_c = T_c(X)$ ``coming into'' $S(1,0)$ and its asymptotic behaviour near $X = 1$ is
\begin{equation}\label{eq:ASYMPTOTICBEHAVIOURNEAR1TCTYPEC}
T_c(X) \sim
\begin{cases}
\lambda_S^-(1-X)^{2/p} \quad &\text{if } 1 < p < 2 \\
\lambda_S(1-X)         \quad &\text{if } p = 2 \\
\lambda_S^+(1-X)         \quad &\text{if } p > 2
\end{cases}
\qquad \text{ for } X \sim 1^-,
\end{equation}
for suitable positive numbers $\lambda_S^-$, $\lambda_S$, and $\lambda_S^+$. For instance, when $p > 2$, substituting the expression $Z(X) = \lambda(1-X)$ into \eqref{eq:EQUATIONOFTHETRAJECTORIESHYP}, we easily obtain
\[
-\lambda = H(X,\lambda(1-X)) \sim \frac{c\lambda(1-X) + mf'(1)(1-X)}{(p-1)\lambda^{(p-1)}(1-X)^{p-1}}, \quad \text{ for } X \sim 1,
\]
i.e.
\[
-(p-1)\lambda^{p}(1-X)^{p-2} \sim c\lambda + mf'(1), \quad \text{for } X \sim 1.
\]
Since the left side goes to zero as $X \to 1$, the previous relation is satisfied only if $\lambda = -mc^{-1}f'(1) := \lambda_S^+ > 0$. The cases $p=2$ and $1 < p < 2$ can be treated similarly, obtaining different values $\lambda_S^-$ and $\lambda_S$ which, as $\lambda_S^+$, depend on $c$, $f'(1)$, $m$ and $p$.
The local analysis of the point $A(a,0)$ is less important in this setting and we skip it. In the Porous Medium case $p=2$ and $m > 1$, it is not difficult to see that $A(a,0)$ is a \emph{focus unstable} if $c < \sqrt{4ma^{m-1}f'(a)}$, while a \emph{node unstable} if $c \geq \sqrt{4ma^{m-1}f'(a)}$ (cfr. with Figure \ref{fig:QUALBEHAVTRAJPMETYPEC}).

\begin{figure}[!ht]
  \centering
  \includegraphics[scale =0.4]{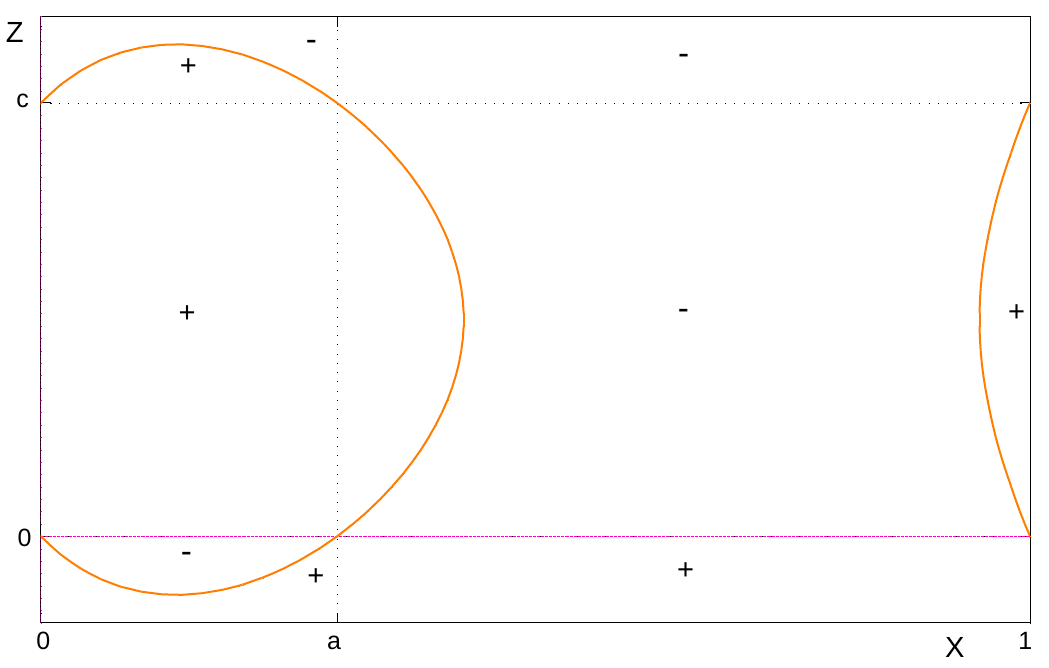}
  \includegraphics[scale =0.4]{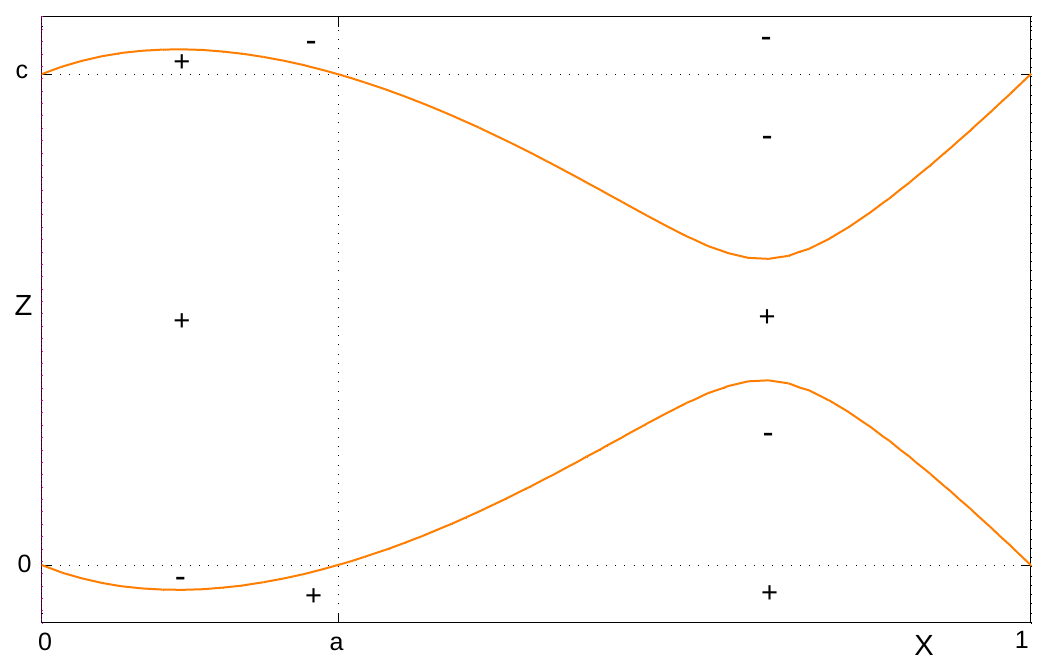}
  \caption{Reactions of type C, range $\gamma > 0$. Null isoclines in the $(X,Z)$-plane for $f(u) = u(1-u)(u - a)$, $a = 0.3$, in the cases $0< c < c_0$ and $c > c_0$, respectively.}\label{fig:ISOCLINESTYPEC}
\end{figure}

\emph{Step 2: Study of the null isoclines.} To obtain a clear picture of the trajectories of the system, we study the \emph{null isoclines} of system \eqref{eq:SYSTEMNONSINGULARTWSTYPECPME2}, i.e., the curve $\widetilde{Z} = \widetilde{Z}(X)$ satisfying
\[
H(X,\widetilde{Z};c) = 0 \qquad \text{ i.e. } \qquad c\widetilde{Z} - |\widetilde{Z}|^p = mX^{\frac{\gamma}{p-1}-1}f(X), \qquad \text{in } [0,1]\times(-\infty,\infty).
\]
First of all, even though it is not of class $C^1$, the function $f_{m,p}(X):=mX^{\frac{\gamma}{p-1}-1}f(X)$ preserves the zeros and the sign of $f(\cdot)$ in $[0,1]$.
\\
Now, let $\{\overline{X}_j\}_{j=1,\ldots,k}$ the set of local maximum points of $f_{m,p}(\cdot)$ in $(a,1)$, $M_j := f_{m,p}(\overline{X}_j)$ for $j =1,\ldots,k$, and $M := \max_{j = 1,\ldots,k} M_j$. Take $c_0 > 0$ so that
\[
\max_{\widetilde{Z} \in [0,c_0^{1/(p-1)}]} c_0\widetilde{Z} - |\widetilde{Z}|^p = M \qquad \text{i.e.} \qquad c_0 = c_0(m,p) := p \left(\frac{M}{p-1}\right)^{(p-1)/p}.
\]
Assume for a moment that there exists a unique $j \in \{1,\ldots,k\}$ such that $M = M_j$. Then it is not difficult to see that for $0 < c < c_0$, the \emph{null isocline} is composed by two disjoint branches: the left one, linking the points $O(0,0)$, $A(a,0)$, $(a,c^{1/(p-1)})$ and $R_c(0,c^{1/(p-1)})$, and the right one, connecting $S(1,0)$ and $(1,c^{1/(p-1)})$. The two branches approach as $c \to c_0$, until they touch at the point $(X_M,(c_0/p)^{1/(p-1)})$ when $c = c_0$, where $f_{m,p}(X_M) = M$. Finally, when $c > c_0$, there are again two disjoint branches: the upper one linking $R_c(0,c^{1/(p-1)})$, $(a,c^{1/(p-1)})$ and $(1,c^{1/(p-1)})$, whilst the lower one joining $O(0,0)$, $A(a,0)$ and $S(1,0)$.
\\
If $M = M_j$ for some $j \in \{j_1,\ldots,j_h\} \subseteq \{1,\ldots,k\}$, the same holds true except for the fact that, when $c = c_0$, the two branches touch at $h$ points belonging to the region $(a,1)\times(0,c^{1/(p-1)})$ while, when $0 < c < c_0$, there are $h$ more null isocline branches made by disjoint closed curves (between the left and the right branch), belonging to the region $(a,1)\times(0,c^{1/(p-1)})$. A qualitative representation is shown in Figure \ref{fig:ISOCLINESTYPEC}. From this analysis it is clear that if there exists a critical speed $c_{\ast}$, then it must be $c_{\ast} < c_0$.

\emph{Step 3: Monotonicity of $T_c(\cdot)$ w.r.t. $c > 0$}. In this crucial step we prove that
\[
\text{for all } \; 0 < c_1 < c_2 \; \text{ then } \; T_{c_2}(X) < T_{c_1}(X), \; \text{ for all } \; a < X < 1
\]
where, of course, $T_c$ is the trajectory ``coming into'' $S(1,0)$. Note that for $0 \leq X \leq a$, $T_c(\cdot)$ is not in general a function of $X$, so that we have to restrict our ``comparison interval'' to $(a,1)$. However, our statement holds true on the interval of definition of $T_c = T_c(X)$.

\noindent Now, fix $0 < c_1 < c_2$. First of all, we note that
\begin{equation}\label{eq:MONOTONICITYRESPECTCEQTRAJECTORIESTYPEC}
\frac{\partial H}{\partial c}(X,Z;c) = \frac{1}{(p-1)X|Z|^{p-2}} > 0, \quad \text{for all } 0 < X \leq 1, \; Z > 0, \; c > 0,
\end{equation}
which implies $H(X,Z;c_1) < H(X,Z;c_2)$.

\noindent Now, assume by contradiction $T_{c_1}$ and $T_{c_2}$ touch in a point $(X_0,T_{c_1}(X_0) = T_{c_2}(X_0))$, with $a < X_0 < 1$. Since $dT_{c_1}(X_0)/dX < dT_{c_2}(X_0)/dX$ by \eqref{eq:MONOTONICITYRESPECTCEQTRAJECTORIESTYPEC}, we have that $T_{c_2}$ stays above $T_{c_1}$ in a small right-neighbourhood $I_0$ of $X_0$ and so, by the continuity of the trajectories, there exists at least another ``contact point'' $X_0 < X_0^+ < 1$, with $T_{c_1}(X_0^+) = T_{c_2}(X_0^+)$. Consequently, for $h > 0$ small enough, we have
\[
\frac{T_{c_2}(X_0^+) - T_{c_2}(X_0^+ - h)}{h} \leq \frac{T_{c_1}(X_0^+) - T_{c_1}(X_0^+ - h)}{h}
\]
and taking the limit as $h \to 0$, we get the contradiction $dT_{c_2}(X_0^+)/dX \leq dT_{c_1}(X_0^+)/dX$. Our assertion follows from the arbitrariness of $a < X_0 < 1$.

\emph{Step 4: Existence and uniqueness of a critical speed $c = c_{\ast}$.} In \emph{Step 0}, we have shown that for $c = 0$ there are not \emph{admissible} TWs, and, in particular, the trajectory $T_0 = T_0(X)$ ``coming into'' $S(1,0)$ stays above the trajectories ``leaving'' the origin $O(0,0)$.
\\
Consequently, thanks to the continuity of the trajectories w.r.t. the parameter $c$ we can conclude the same, for \emph{small} values of $c > 0$, i.e., naming $T^+_c = T^+_c(X)$ and $T^-_c = T^-_c(X)$ the trajectories from $R_c(0,c^{1/(p-1)})$ and $O(0,0)$, respectively, we have that $T_c(X)$ is above $T^+_c(X)$ and $T^-_c(X)$ in $[0,1]$ (note that for $c = 0$, $R_0 = O$ and both $T^+_0$ and $T^-_0$ ``leave'' $O$).

\noindent In particular, the study of the \emph{null isoclines} carried out in \emph{Step 2} shows that $T^+_c(X=a) > c^{1/(p-1)}$ for all $c > 0$, and so, using the monotonicity of $T_c$ w.r.t. $c > 0$ proved in \emph{Step 3}, we conclude that for $c > 0$ large enough it must be $T_c(X = a) < T^+_c(X = a)$, which means that for large $c > 0$, $T_c(X)$ stays below $T^+_c(X)$, in $[0,1]$ by uniqueness of the trajectories. This means that there exists a critical speed $c_{\ast} = c_{\ast}(m,p,f)$ such that $T_{c_{\ast}}(X) = T^+_{c_{\ast}}(X)$ for all $X \in [0,1]$, and the uniqueness of $c_{\ast}$ follows by the strict inequality in \eqref{eq:MONOTONICITYRESPECTCEQTRAJECTORIESTYPEC}. In other words, the trajectories $T^+_c$ and $T_c$ approach as $c < c_{\ast}$ grows until they touch (i.e. they coincide) for $c = c_{\ast}$, while for $c > c_{\ast}$ they are ordered in the opposite way w.r.t. the range $c < c_{\ast}$, i.e. $T^+_c(X) > T_c(X)$ in $[0,1]$ for all $c > c_{\ast}$.

\noindent We conclude this step showing that the trajectory $T_{c_{\ast}}$ linking  $S(1,0)$ \normalcolor and $R_{c_{\ast}}(0,c_{\ast}^{1/(p-1)})$ corresponds to an \emph{admissible} TW \emph{profile}$X(\xi) = \varphi(\xi)$ of a \emph{finite} TW $u(x,t) = \varphi(x-ct)$, i.e., $\varphi(-\infty) = 1$ and $\varphi(\xi) = 0$, for all $\xi \geq \xi_0$, for some $-\infty < \xi_0 < +\infty$. The fact that $X(-\infty) = \varphi(-\infty) = 1$, follows by integrating the first ODEs in \eqref{eq:SYSTEMNONSINGULARTWTYPECPME1}:
\[
-m\frac{dX}{d\xi} = X^{1- \frac{\gamma}{p-1}}Z,
\]
by separation of variables and recalling the asymptotic behaviour of $T_c$ near $X =1$, given in formula \eqref{eq:ASYMPTOTICBEHAVIOURNEAR1TCTYPEC}. Indeed, fixing $0 < X_0 < X_1 < 1$ and taking $Z(X) = T_{c_{\ast}}(X)$, the local analysis around the saddle point $S(1,0)$ carried out in \emph{Step 1} allows us to estimate the time $\xi_1$ in which the profile reaches the level $u=1$
\[
\xi_0 - \xi_1 = m\int_{X_0}^{X_1}\frac{X^{\frac{\gamma}{p-1}-1}}{T_{c_{\ast}}(X)}dX \sim m\int_{X_0}^{X_1}\frac{dX}{T_{c_{\ast}}(X)} \qquad \text{for } X_1 \sim 1,
\]
and, since the second integral diverges for $X_1 \sim 1$, we can conclude that $\xi_1 = -\infty$, i.e., $X(-\infty) = \varphi(-\infty)$ = 1.
\\
On the other hand, since $T_{c_{\ast}}(X) \sim c_{\ast}^{1/(p-1)}$ for $X \sim 0$, proceeding as before we deduce that
\[
\xi_0 - \xi_1= m\int_{X_0}^{X_1}\frac{X^{\frac{\gamma}{p-1}-1}}{T_{c_{\ast}}(X)}dX \sim m
c_{\ast}^{-1/(p-1)}\int_{X_0}^{X_1} X^{\frac{\gamma}{p-1}-1} dX \qquad \text{for } X_0 \sim 0,
\]
which means that the time $\xi_0$, in which the profile gets to the level $u=0$, is \emph{finite} and, moreover, taking $X_0 \sim 0$ and relabeling $X_1 = X = \varphi$, it holds (up to a multiplicative constant)
\begin{equation}\label{eq:ADMISSIBLEPROFILECASTDARCYLAWHYP}
\xi_0 - \xi \sim X^{\frac{\gamma}{p-1}}(\xi) \quad \Leftrightarrow \quad X^{\frac{\gamma}{p-1}}(\xi) = \varphi^{\frac{\gamma}{p-1}}(\xi) \sim \xi_0 - \xi, \qquad \text{for } \xi \sim \xi_0,
\end{equation}
which gives the behaviour of the \emph{finite} TW near the \emph{free boundary} point $-\infty < \xi_0 < +\infty$, according to the Darcy law, see \cite{V2:book}.

\emph{Step 5: Non existence of TWs for $c > c_{\ast}$}. We are left to prove that there are not \emph{admissible} TW solutions when $c > c_{\ast}$. This follows from the fact that if the trajectory $T_c$ joins $O(0,0)$ and $S(1,0)$, then the resulting connection is not \emph{admissible} since the derivative of the corresponding \emph{profile} changes sign. Indeed, using the continuity of the trajectory w.r.t. the speed of propagation, we show that for all $c > 0$, there exists a unique trajectory $T^-_c = T^-_c(X)$ ``leaving'' $O(0,0)$ (see \emph{Step 0}) and a simple computation shows that
\[
T^-_c(X) \sim \frac{mf'(0)}{c}X^{\frac{\gamma}{p-1}}, \quad \text{for } X \sim 0^+.
\]
Hence, if $T_c$ links $O(0,0)$ and $S(1,0)$, it must coincide with $T^-_c$ and so, the derivative of its profile must change sign, i.e., it is not an \emph{admissible} profile.
\paragraph{Proof of Theorem \ref{THEOREMEXISTENCEOFTWSPMEREACTIONTYPEC}: Part (i), range $\boldsymbol{\gamma = 0}$.} Fix $m > 0$ and $p > 1$ such that $\gamma = 0$. Proceeding as before, we consider the variables \eqref{eq:NSTANDCVARTWSHYPPME}, which in the case $\gamma = 0$ read as
\[
X = \varphi \qquad \text{ and } \qquad Z = -m (\log X)' = -mX^{-1}X' \quad (\geq 0),
\]
and we get the system
\begin{equation}\label{eq:SYSTEMNONSINGULARTWTYPECPME1GAMMA0}
-m\frac{dX}{d\xi} = XZ, \quad\quad -|Z|^{p-2} \frac{dZ}{d\xi} = cZ - |Z|^p - F(X),
\end{equation}
where $F(X) = mX^{-1}f(X)$, and, after re-parametrization $d\xi = |Z|^{p-2}d\tau$,
\[
\frac{dX}{d\tau} = (p-1)X|Z|^{p-2}Z, \quad\quad \frac{dZ}{d\tau} = cZ - |Z|^p - F(X).
\]
Since $F(0) = mf'(0) < 0$, $F(a) = F(1) = 0$, with $F(X) < 0$ in $(0,a)$ and $F(X) > 0$ in $(a,1)$, the critical points are now $S(1,0)$, $A(a,0)$, $R_{\lambda_1}(0,\lambda_1)$, and $R_{\lambda_2}(0,\lambda_2)$, where $\lambda_1 = \lambda_1(c) < 0 < c^m < \lambda_2 = \lambda_2(c)$ are the solutions of the equation
\[
cZ - |Z|^p = F(0), \qquad c > 0.
\]
The ``new'' \emph{equation of the trajectories} is
\[
\frac{dZ}{dX} = \frac{cZ - |Z|^p - F(X)}{(p-1)X|Z|^{p-2}Z} := H(X,Z;c).
\]
In the next paragraph we repeat the scheme followed before looking for trajectories in the strip $[0,1]\times[0,\infty)$ connecting $S(1,0) \leftrightsquigarrow R_{\lambda_2}(0,\lambda_2)$ for a specific speed of propagation $c = c_{\ast}(m,p,f)$.

\emph{Step 0': Case c = 0}. If $c = 0$, the \emph{equation of the trajectories} reads
\[
\frac{dZ}{dX} = -\frac{|Z|^p + F(X)}{(p-1)X|Z|^{p-2}Z} := H(X,Z;0).
\]
The \emph{null isoclines} are composed by two branches, the upper one linking $R_{\lambda_2}(0,\lambda_2)$ and $A(a,0)$, and the lower one joining $R_{\lambda_1}(0,\lambda_1)$ and $A(a,0)$, where in this easier case
\[
\lambda_1 = \lambda_1(0) = -\sqrt[p]{-mf'(0)}, \qquad \lambda_2 = \lambda_2(0) = \sqrt[p]{-mf'(0)}.
\]
Employing the Lyapunov linearization method, it is not difficult to prove that $R_{\lambda_1}(0,\lambda_1)$ and $R_{\lambda_2}(0,\lambda_2)$ are two \emph{saddle} points. So, there are exactly two trajectories $T^- = T^-(X)$ and $T_+ = T_+(X)$ ``coming from'' $R_{\lambda_1}(0,\lambda_1)$ and $R_{\lambda_2}(0,\lambda_2)$, respectively, lying in the strip $[0,1]\times(-\infty,\infty)$ in the $(X,Z)$-plane. Moreover, since
\[
H(X,-Z;0) = -H(X,Z;0), \quad \text{ for all } 0 \leq X \leq 1, \; Z \in \RR,
\]
we deduce that $T^- \equiv T^+$. At the same time, exactly as in the case $\gamma > 0$ we have a trajectory $T_0 = T_0(X) > 0$ ``coming into'' $S(1,0)$ (see \emph{Step 1} of the case $\gamma > 0$). Assuming \eqref{eq:NOHOMOCLINICASSUMPTIONONF}, i.e.,
$\int_0^1u^{m-1}f(u)\,du > 0$, it follows that $T_+(X) < T_0(X)$ for all $0 \leq X \leq 1$, with $T_0(X) \sim +\infty$ for $X \sim 0$. This follows by using the same techniques of the case $\gamma > 0$. In particular, it is simple to see that the same construction works if we take $\gamma = 0$ and formula \eqref{eq:EQUATIONFORSXZPLANETYPEC} holds. Consequently, there are not \emph{admissible} TWs for $c = 0$.

\begin{figure}[!ht]
  \centering
  \includegraphics[scale =0.4]{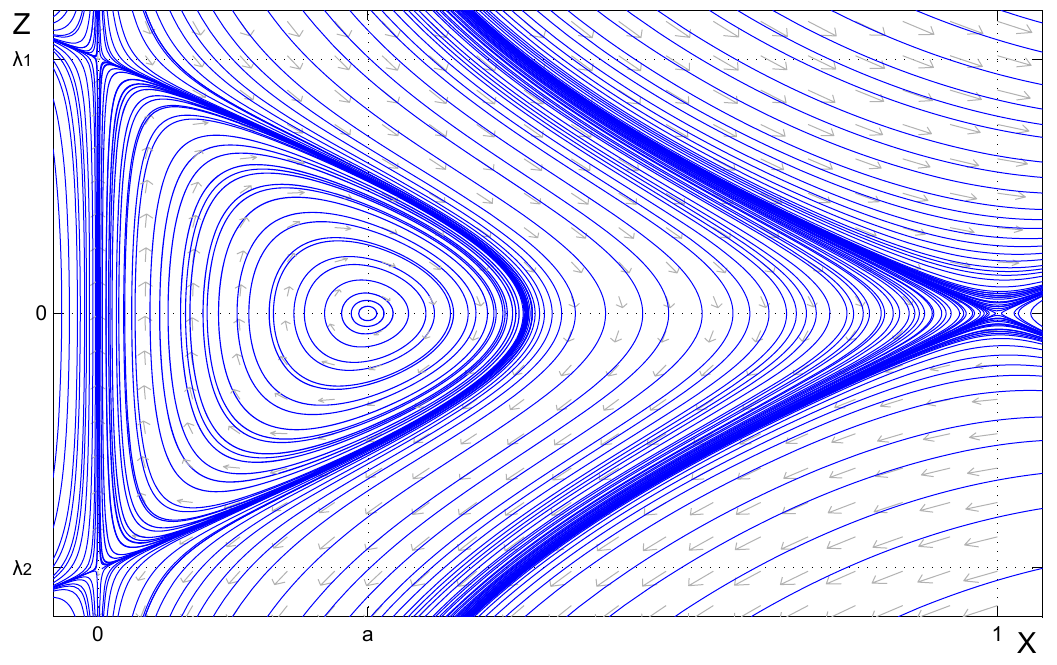}
  \caption{Reactions of type C, range $\gamma = 0$, case $c = 0$. Qualitative behaviour of the trajectories in the $(X,Z)$-plane for $f(u) = u(1-u)(u - a)$, $a = 0.3$.}\label{fig:QUALBEHAVTRAJPMETYPECGAMMA0C0}
\end{figure}

\emph{Step 1': Local analysis of $S(1,0)$.} This step coincides with \emph{Step 1} of the case $\gamma > 0$, since the nature of the critical point $S(1,0)$ does not change if we take $\gamma = 0$. This can be easily seen noting that $F(X) \sim f_{m,p}(X) \sim -mf'(1)(1-X)$ for $X \sim 1$.

\emph{Step 2': Study of the null isoclines.} We proceed as before by studying the solutions of the equation
\[
c\widetilde{Z} - |\widetilde{Z}|^p = mX^{-1}f(X), \qquad \text{in } [0,1]\times(-\infty,\infty).
\]
As before, we find that there exists $c_0 > 0$ such that for $c > c_0$, we have again two branches: the upper one, linking $R_{\lambda_1}(0,\lambda_1)$, $(a,c^m)$ and $(1,c^m)$, while the lower one joining $R_{\lambda_2}(0,\lambda_2)$, $(a,0)$ and $(1,0)$. Depending on the number $h$ of points in $(a,1)$ in which the global maximum of $F(\cdot)$ is attained, we have that the two branches approach as $c \to c_0$ until they touch at $h$ points in the region $(a,1)\times(0,c^{1/(p-1)})$ for $c = c_0$. Finally, for $0 < c < c_0$ the \emph{null isoclines} are composed by two branches: the left one, linking the points $R_{\lambda_1}(0,\lambda_1)$, $(a,c^m)$, $(a,0)$ and $R_{\lambda_2}(0,\lambda_2)$, whilst the second linking $(1,c^m)$ and $S(1,0)$, and $h$ disjoint closed curves (between these two branches), belonging to the region $(a,1)\times(0,c^{1/(p-1)})$ (cfr. with \emph{Step 2}, range $\gamma > 0$).

\emph{Step 3': Monotonicity of $T_c(\cdot)$ w.r.t. $c > 0$}. If we denote again with $T_c = T_c(X)$ the trajectory ``coming into'' $S(1,0)$, the proof of monotonicity property of $T_c$ w.r.t. to $c > 0$ coincides with the one done in \emph{Step 3} of the case $\gamma > 0$.

\emph{Step 4': Existence and uniqueness of a critical speed $c = c_{\ast}$.} The existence of a unique critical speed $c_{\ast} = c_{\ast}(m,p,f)$ with corresponding trajectory linking $S(1,0)$ and $R_{\lambda_2}(\lambda_2(c_{\ast}),0)$ follows exactly as in the case $\gamma > 0$. The unique (important) difference is the fact that the the TW is positive everywhere. Indeed, integrating the first equation in \eqref{eq:SYSTEMNONSINGULARTWTYPECPME1GAMMA0} along the trajectory $T_{c_{\ast}} = T_{c_{\ast}}(X) \sim \lambda_1(c_{\ast})$ for $X \sim 0$, we obtain
\[
\xi_0 - \xi_1= m\int_{X_0}^{X_1}\frac{1}{XT_{c_{\ast}}(X)}dX \sim m
\lambda_2(c_{\ast})\int_{X_0}^{X_1} X^{-1} dX \qquad \text{for } X_0 \sim 0,
\]
from which we deduce $\varphi(+\infty) = X(+\infty) = 0$, i.e., the TW profile $X(\xi) = \varphi(\xi)$ reaches the level $u = 0$ in infinite time.

\begin{figure}[!ht]
  \centering
  \includegraphics[scale =0.4]{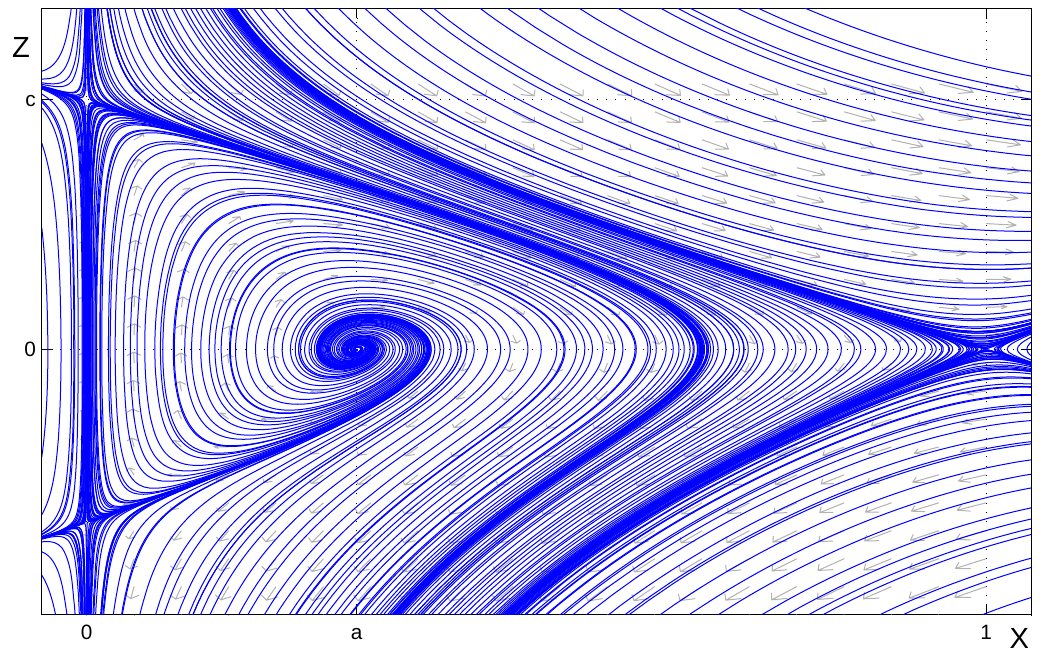}
  \includegraphics[scale =0.4]{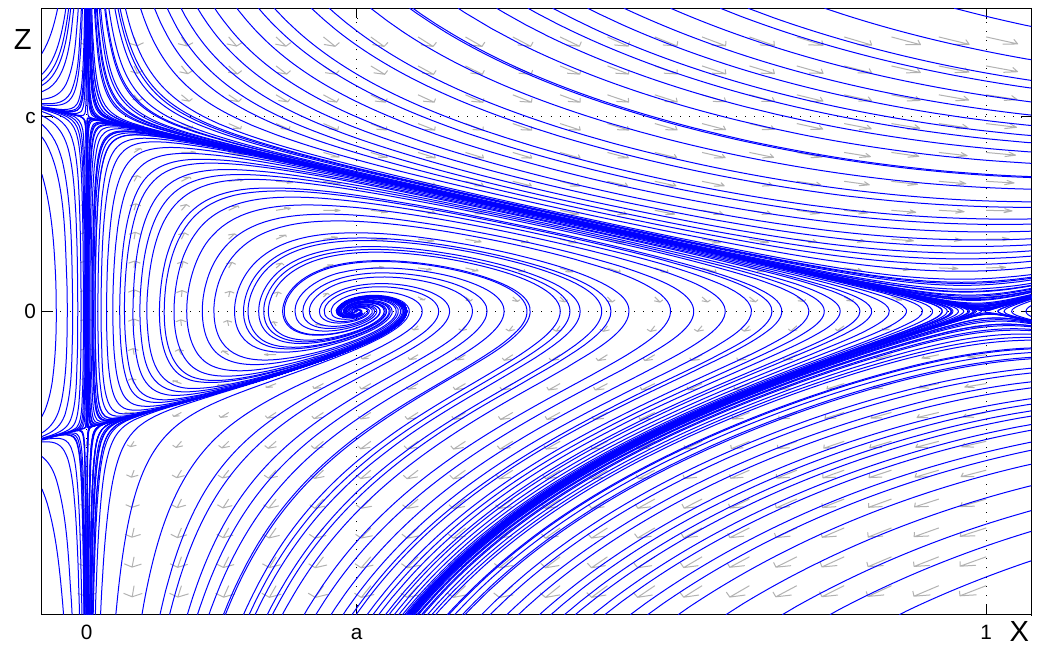}
  \includegraphics[scale =0.4]{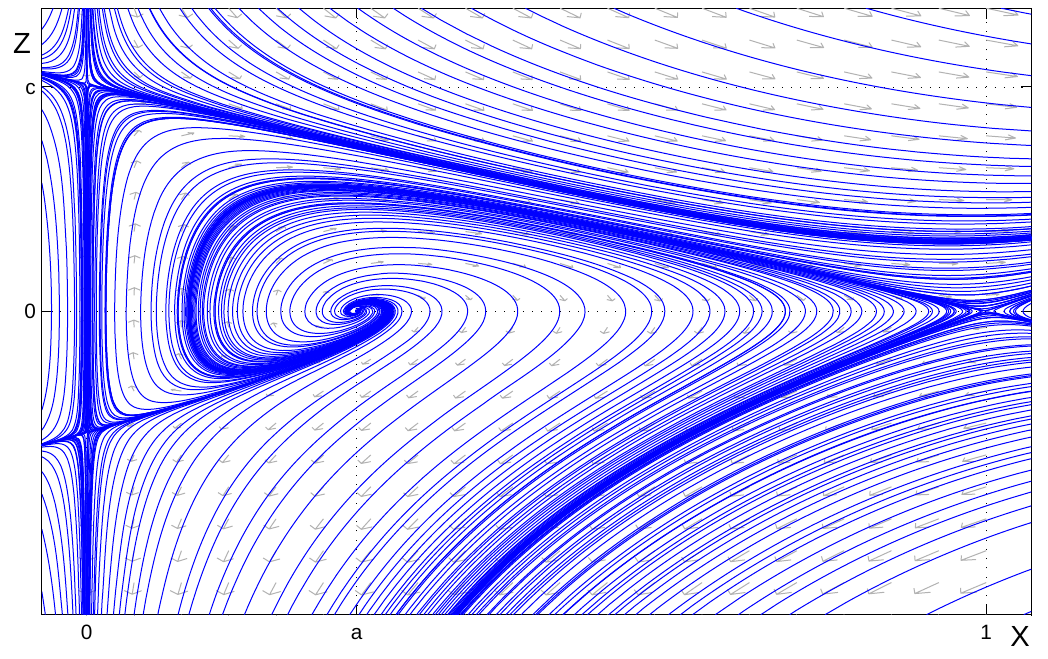}
  \caption{Reactions of type C, range $\gamma = 0$. Qualitative behaviour of the trajectories in the $(X,Z)$-plane for $f(u) = u(1-u)(u - a)$, $a = 0.3$. The first picture shows the case $0< c < c_{\ast}$, while the others the cases $c = c_{\ast}$ and $c > c_{\ast}$, respectively.}\label{fig:QUALBEHAVTRAJPMETYPECGAMMA0}
\end{figure}

\emph{Step 5: Non existence of TWs for $c > c_{\ast}$}. Proving the non existence of \emph{admissible} TW profiles is easier than the case $\gamma > 0$, since from the study of the critical points and the \emph{null isoclines} it follows that there cannot exist nonnegative trajectories linking $S(1,0)$ and $O(0,0)$. The qualitative behaviour of the trajectories in the $(X,Z)$-plane is reported in Figure \ref{fig:QUALBEHAVTRAJPMETYPECGAMMA0}.

\paragraph{Proof of Theorem \ref{THEOREMEXISTENCEOFTWSPMEREACTIONTYPEC}: Part (ii), range $\boldsymbol{\gamma >0}$.} Fix $m > 0$ and $p > 1$ such that $\gamma > 0$, and $0 < a < 1$. We study the existence of \emph{a-admissible} TW solutions for equation \eqref{eq:REACDIFFTYPCDIM1} with reaction satisfying \eqref{eq:ASSUMPTIONSONTHEREACTIONTERMTYPEDPME}. Respect to the previous case, our proof strongly relies on some results carried out in Theorem 2.1 of \cite{AA-JLV:art}. Indeed, as we have mentioned in the introduction, reaction terms satisfying \eqref{eq:ASSUMPTIONSONTHEREACTIONTERMTYPEDPME} are of the Fisher-KPP type if we restrict them to the interval $[0,a] \subset [0,1]$, in the sense that

\[
\begin{cases}
f(0) = f(a) = 0, \quad &f(u) > 0, \text{ in } (0,a) \\
f \in C^1([0,a]), \quad &f'(0) > 0, \; f'(a) < 0, \quad f(u) \leq f'(0)u, \; 0 \leq u \leq a. \\
\end{cases}
\]
\normalcolor
For this reason, it follows that the qualitative behaviour of the trajectories in the strip $[0,a]\times(-\infty,\infty)$ of the $(X,Z)$-plane is the same of the ones studied in Theorem 2.1 \cite{AA-JLV:art} in the larger strip $[0,1]\times(-\infty,\infty)$, where the Fisher-KPP case has been analyzed.  In this way, it is easily seen that the study of the trajectories corresponding to \emph{a-admissible} TW solutions of equation \eqref{eq:REACDIFFTYPCDIM1} (with reaction of type C') is reduced to the study of \emph{admissible} TWs for equation \eqref{eq:REACDIFFTYPCDIM1} with a reaction of Fisher-KPP type (or type A). In view of this explaination, some part of the following proof coincide one of Theorem 2.1 of \cite{AA-JLV:art}, so that, for the reader's convenience, we will try to report the most important ideas, quoting the specific paragraphs of \cite{AA-JLV:art} for each technical detail.

Now, following the beginning of the proof of part (i) case $\gamma > 0$, we obtain systems \eqref{eq:SYSTEMNONSINGULARTWTYPECPME1}:
\[-m\frac{dX}{d\xi} = X^{1- \frac{\gamma}{p-1}}Z, \quad\quad -m(p-1)X^{\frac{\gamma}{p-1}}|Z|^{p-2} \frac{dZ}{d\xi} = cZ - |Z|^p - mX^{\frac{\gamma}{p-1}-1}f(X),
\]
and \eqref{eq:SYSTEMNONSINGULARTWSTYPECPME2}:
\[
\frac{dX}{d\tau} = (p-1)X|Z|^{p-2}Z, \quad\quad \frac{dZ}{d\tau} = cZ - |Z|^p - mX^{\frac{\gamma}{p-1}-1}f(X),
\]
with critical points $O(0,0)$, $S(1,0)$, $A(a,0)$, $R_c(0,c^{1/(p-1)})$, and with the \emph{equation of the trajectories} \eqref{eq:EQUATIONOFTHETRAJECTORIESHYP}:
\[
\frac{dZ}{dX} = \frac{cZ - |Z|^p - mX^{\frac{\gamma}{p-1}-1}f(X)}{(p-1)X|Z|^{p-2}Z} := H(X,Z;c)
\]
Even though, they formally coincide, the reaction term is now of type C', i.e., it satisfies \eqref{eq:ASSUMPTIONSONTHEREACTIONTERMTYPEDPME}. The main structural difference between the reaction of type C case and the type C' case is that the study of the case $c = 0$ is not needed for our purposes. This is basically due to the fact that the critical point $A(a,0)$ is a \emph{saddle} type critical point, for all $0 < a < 1$, as we will see in a moment.

\begin{figure}[!ht]
  \centering
  \includegraphics[scale =0.4]{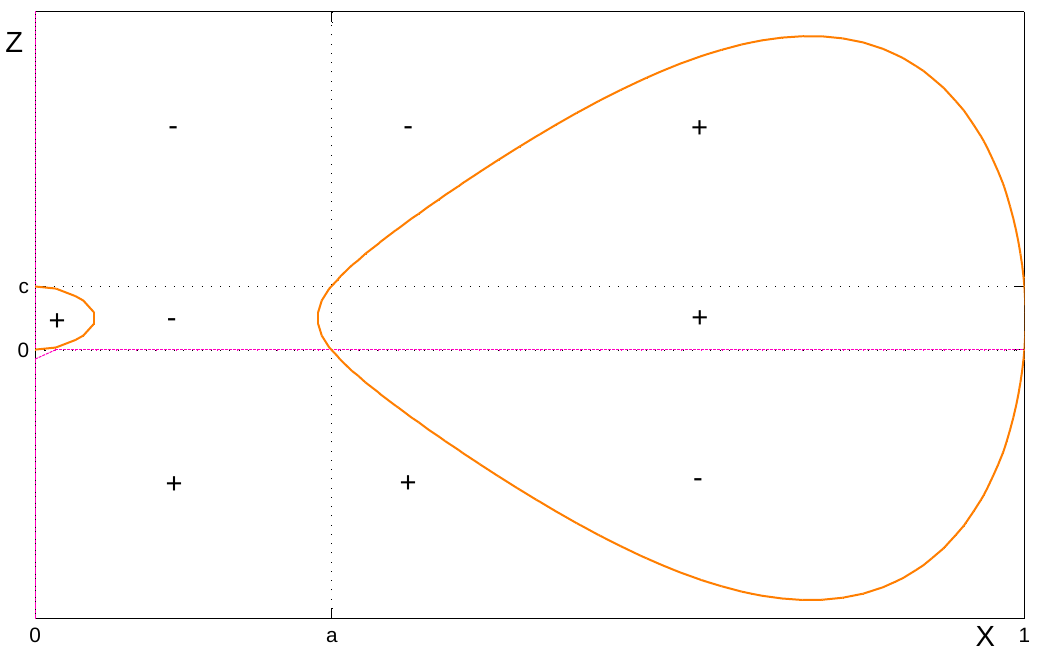}
  \includegraphics[scale =0.4]{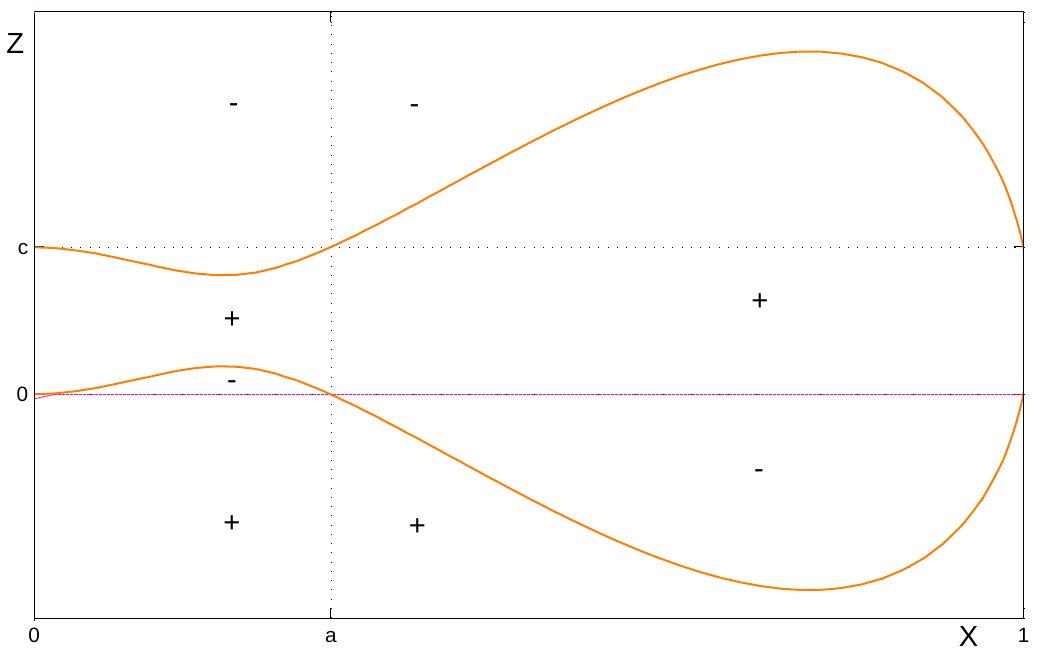}
  \caption{Reactions of type C', range $\gamma > 0$. Null isoclines in the $(X,Z)$-plane for $f(u) = u(1-u)(a - u)$, $a = 0.3$, in the cases $0< c < c_0$ and $c > c_0$, respectively.}\label{fig:ISOCLINESTYPECPRIME}
\end{figure}

\emph{Step 1: Local analysis of $A(a,0)$ and $S(1,0)$.} Let us take $c > 0$. Proceeding as in the proof of Theorem \ref{THEOREMEXISTENCEOFTWSPMEREACTIONTYPEC} (see \emph{Step 1}), and recalling that now $f'(1) > 0$, while $f'(a) < 0$, we deduce that $A(a,0)$ is a \emph{saddle} type critical point, and formulas \eqref{eq:ASYMPTOTICBEHAVIOURNEAR1TCTYPEC} hold replacing $f'(1)$ with $f'(a)$.
\\
For what concerns the point $S(1,0)$, we can conclude it has a \emph{focus}/\emph{node} nature from the study of the \emph{null isoclines} we perform in \emph{Step 2}.

Now, let $T_c = T_c(X)$ be the trajectory entering in $A(a,0)$ with $T_c(X) > 0$ for all $0 < X < a$. In the next paragraphs, following the proof of Theorem 2.1 of \cite{AA-JLV:art} and the ideas of part (i), we show that there exists a unique $c_{\ast} = c_{\ast}(m,p,f)$ such that $T_{c_{\ast}}$ links $R_{c_{\ast}}(0,c_{\ast}^{1/(p-1)})$ and $A(a,0)$ and we prove that this trajectory corresponds to a \emph{finite} TW profile. Secondly, we show that for all $c > c_{\ast}$, $T_c$ joins $O(0,0)$ and $A(a,0)$, and it corresponds to a \emph{positive} TW profile. Finally, we prove that there are not connections of the type $A(a,0) \leftrightsquigarrow R_c(0,c^{1/(p-1)})$ and/or $A(a,0) \leftrightsquigarrow O(0,0)$ for $c < c_{\ast}$, i.e. there are not any \emph{a-admissible} TW profiles for $c < c_{\ast}$.

\emph{Step 2: Study of the null isoclines}. We study the \emph{null isoclines} of system \eqref{eq:SYSTEMNONSINGULARTWSTYPECPME2}, i.e., the curve $\widetilde{Z} = \widetilde{Z}(X)$ satisfying
\[
c\widetilde{Z} - |\widetilde{Z}|^p = mX^{\frac{\gamma}{p-1}-1}f(X),
\]
exactly as in \emph{Step 2} of the proof of part (i). W.r.t. to the bistable framework, the situation here is ``inverted''. In the case in which there exists a unique global maximum point of $f_{m,p}(X) := mX^{\frac{\gamma}{p-1}-1}f(X)$ in $[0,a]$, we obtain that there exists $c_0 > 0$ such that for $0 < c < c_0$, the \emph{null isocline} is composed of two disjoint branches: the left one, linking the points $O(0,0)$ and $R_c(0,c^{1/(p-1)})$, and the right one, connecting $S(1,0)$, $A(a,0)$, $(a,c^{1/(p-1)})$ and $(1,c^{1/(p-1)})$. For $c > c_0$, we have again two branches: the upper one linking $R_c(0,c^{1/(p-1)})$, $(a,c^{1/(p-1)})$ and $(1,c^{1/(p-1)})$, whilst the lower one joining $O(0,0)$, $A(a,0)$ and $S(1,0)$.  As before, the two branches approach as $c \to c_0$, and they touch at a point when $c = c_0$.
\\
If $f_{m,p}(\cdot)$ has more than one global maximum point in $[0,a]$ the analysis is very similar to the bistable case and we refer the reader to \emph{Step 2} of the proof of Theorem \ref{THEOREMEXISTENCEOFTWSPMEREACTIONTYPEC}. Again we see that if our $c_{\ast}$ exists, then it has to be $c_{\ast} < c_0$. The qualitative shape of the \emph{null isoclines} for reactions of type C' in the cases $0< c < c_0$ and $c > c_0$ is reported in Figure \ref{fig:ISOCLINESTYPECPRIME}. We stress that the shape of the \emph{null isoclines} in the rectangle $[0,a]\times[0,c^{1/(p-1)}]$ is (of course) the same of the one found for Fisher-KPP reactions in the rectangle $[0,1]\times[0,c^{1/(p-1)}]$ (cfr. with \emph{Step 1} of Theorem 2.1 of \cite{AA-JLV:art}). \normalcolor

\begin{figure}[!ht]
  \centering
  \includegraphics[scale =0.4]{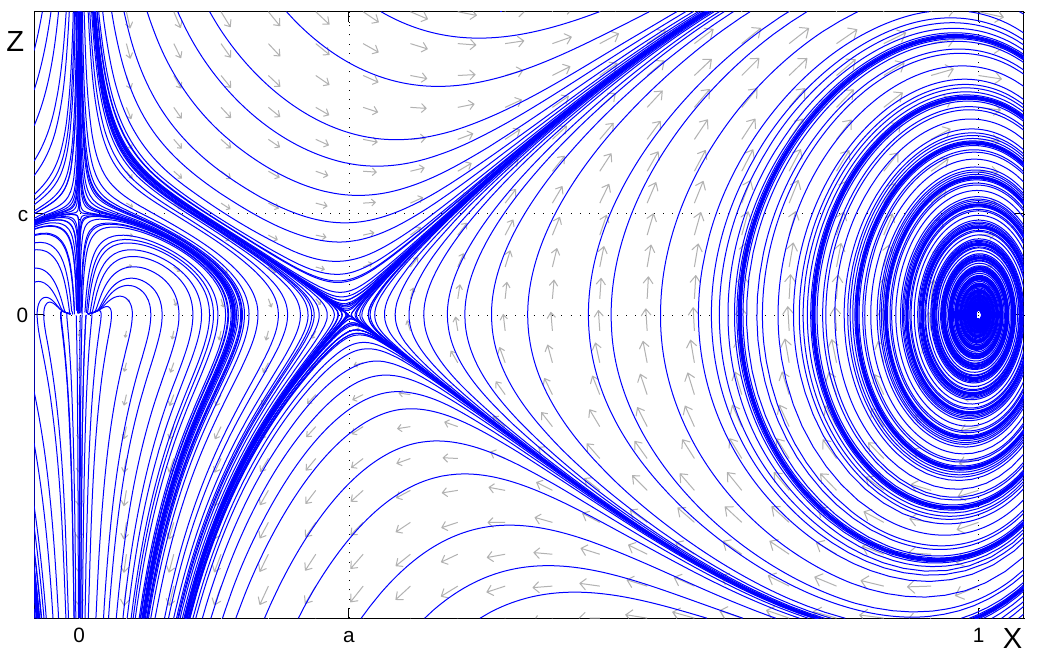}
  \includegraphics[scale =0.4]{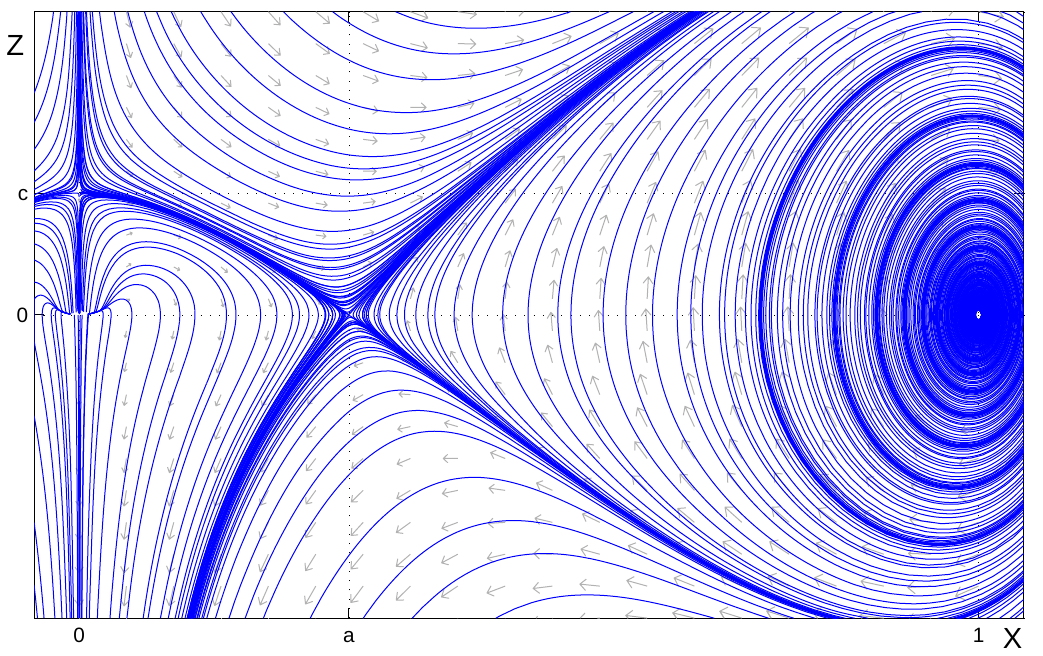}
  \includegraphics[scale =0.4]{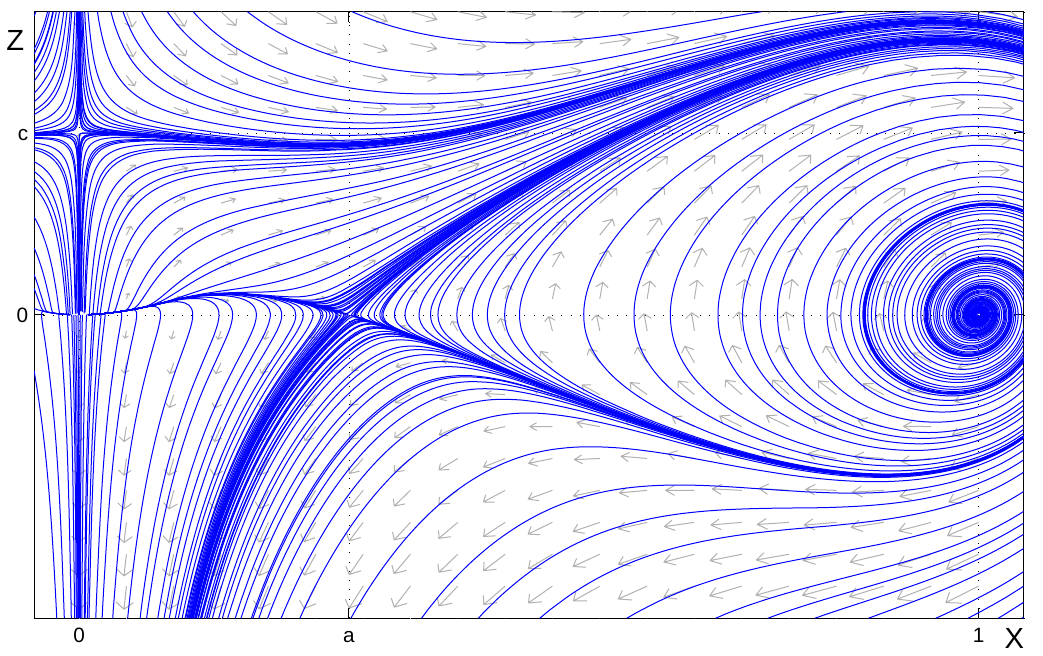}
  \caption{Reactions of type C', range $\gamma > 0$. Qualitative behaviour of the trajectories in the $(X,Z)$-plane for $f(u) = u(1-u)(a - u)$, $a = 0.3$, in the ranges $0< c < c_{\ast}$, $c = c_{\ast}$ and $c > c_{\ast}$, respectively.}\label{fig:QUALBEHAVTRAJPMETYPECPRIME}
\end{figure}

\emph{Step 3: Existence and uniqueness of a critical speed $c = c_{\ast}$}. As we have explained in \emph{Step 1}, we have to prove the existence and the uniqueness of a speed $c_{\ast} = c_{\ast}(m,p,f)$ such that $T_{c_{\ast}}$ links $R_{c_{\ast}}(0,c_{\ast}^{1/(p-1)})$ and $A(a,0)$ with corresponding TW profile  vanishing in a half-line. Consequently, the proof of this fact coincides with what was proved in \emph{Step 3} of Theorem 2.1 of \cite{AA-JLV:art}, substituting the point $S(1,0)$ with $A(a,0)$.

\emph{Step 4: The cases $0 < c < c_{\ast}$ and $c > c_{\ast}$}. We have to show that for $0 < c < c_{\ast}$, the are not \emph{a-admissible} TW, while to each $c > c_{\ast}$, it corresponds exactly one \emph{a-admissible} TW and it is \emph{positive}. Again it is sufficient to adapt \emph{Step 4} and \emph{Step 5} of Theorem 2.1 of \cite{AA-JLV:art} and we conclude the proof. A qualitative representation of the trajectories for $c < c_{\ast}$, $c = c_{\ast}$ and $c > c_{\ast}$ is shown in Figure \ref{fig:QUALBEHAVTRAJPMETYPECPRIME}. $\Box$

\paragraph{Proof of Theorem \ref{THEOREMEXISTENCEOFTWSPMEREACTIONTYPEC}: Part (ii), range $\boldsymbol{\gamma = 0}$.} Fix $m > 0$ and $p > 1$ such that $\gamma = 0$, and $0 < a < 1$. As in the previous part, we base our proof on what proved in \cite{AA-JLV:art}, see Theorem 2.2.

In this setting (cfr. with Part (i), range $\gamma = 0$), \eqref{eq:SYSTEMNONSINGULARTWTYPECPME1} and \eqref{eq:SYSTEMNONSINGULARTWSTYPECPME2} can be written as
\begin{equation}\label{eq:SYSTEMNONSINGULARTWTYPECPME1GAMMA0PRIME}
-m\frac{dX}{d\xi} = XZ, \quad\quad -|Z|^{p-2} \frac{dZ}{d\xi} = cZ - |Z|^p - F(X),
\end{equation}
where $F(X) = mX^{-1}f(X)$,
\[
\frac{dX}{d\tau} = (p-1)X|Z|^{p-2}Z, \quad\quad \frac{dZ}{d\tau} = cZ - |Z|^p - F(X),
\]
and the \emph{equation of the trajectories} is
\[
\frac{dZ}{dX} = \frac{cZ - |Z|^p - F(X)}{(p-1)X|Z|^{p-2}Z} := H(X,Z;c).
\]
This time it holds $F(0) = mf'(0) > 0$, $F(a) = F(1) = 0$, with $F(X) > 0$ in $(0,a)$ while $F(X) < 0$ in $(a,1)$, and for all $c > 0$, $S(1,0)$, $A(a,0)$ are critical points. Now, studying the equation
\[
cZ - |Z|^p = F(0), \qquad c > 0,
\]
for $X = 0$, and defining
\[
c_{\ast}(m,p,f) := p(m^2f'(0))^{\frac{1}{mp}},
\]
it follows that if $c < c_{\ast}$ then there are not other critical points, if $c = c_{\ast}$ there is one more critical point
\[
R_{\lambda_{\ast}}(0,\lambda_{\ast}), \qquad \lambda_{\ast}:=(c_{\ast}/p)^m = (m^2f'(0))^{\frac{1}{p}},
\]
while if $c > c_{\ast}$ there are two more critical points $R_i(0,\lambda_i)$, $i=1,2$ where $\lambda_i = \lambda_i(c)$ and $0 < \lambda_1 < \lambda_{\ast} < \lambda_2 < c^m$.

\emph{Step 1': Local analysis of $A(a,0)$ and $S(1,0)$.} This step coincides with \emph{Step 1} of part (ii), case $\gamma > 0$.

\emph{Step 2': Study of the null isoclines}. As always, we study the solutions of the equation
\[
c\widetilde{Z} - |\widetilde{Z}|^p = F(X), \qquad c > 0,
\]
finding that for $c > c_{\ast}$, there are two branches: the upper one, linking $R_{\lambda_2}(0,\lambda_2)$, $(a,c^m)$ and $(1,c^m)$, whilst the lower one joining $R_{\lambda_1}(0,\lambda_1)$, $A(a,0)$ and $S(1,0)$, while when $c = c_{\ast}$ the branch touches the $Z$-axis at the point $R_{\lambda_{\ast}}(0,\lambda_{\ast})$. Finally, for $0 < c < c_{\ast}$, the \emph{null isoclines} are composed by a branch linking the points $(1,c^m)$, $(a,c^m)$, $A(a,0)$, and $S(1,0)$, and a certain number (depending on the number of global maximum point of $F(\cdot)$ on $[0,a]$) of disjoint closed curves between the $Z$-axis and the first branch, belonging to the region $(0,a)\times(0,c^m)$ (cfr. with \emph{Step 2} of the proof of Theorem \ref{THEOREMEXISTENCEOFTWSPMEREACTIONTYPEC}).

\emph{Step 3': Existence and uniqueness of a critical speed $c = c_{\ast}$}. In this step, we have to prove the existence of a trajectory $T_{c_{\ast}}$ linking $A(a,0)$ with $R_{\lambda_{\ast}}(0,\lambda_{\ast})$, corresponding to an \emph{a-admissible positive} TW profile. This easily follows remembering the scaling property we explained before and substituting $S(1,0)$ with $A(a,0)$ in the proof of Theorem 2.2 of \cite{AA-JLV:art}.

\noindent We anticipate that in the PDEs part we will need more information about the asymptotic behaviour of the TW \emph{profile} $X(\xi) = \varphi(\xi)$ corresponding to the critical speed $c_{\ast}$. This was studied in \cite{AA-JLV:art} (cfr. with Theorem 2.2 and Section 11) where it was proved that the ``critical'' TW satisfies
\begin{equation}\label{eq:ASYMPTOTICSOFCRITICALPROFILEGAMMA0}
\varphi(\xi) \sim a_0 |\xi|^{\frac{2}{p}}e^{-\frac{\lambda_{\ast}}{m}\xi} = a_0 |\xi|^{\frac{2}{p}}\exp\Big(-m^{\frac{2-p}{p}}f'(0)^{\frac{1}{p}}\xi\Big) \quad \text{for }\; \xi \sim +\infty,
\end{equation}
where as before $\lambda_{\ast}:=(c_{\ast}/p)^m$ and $a_0 > 0$ is a suitable constant.

\begin{figure}[!ht]
  \centering
  \includegraphics[scale =0.4]{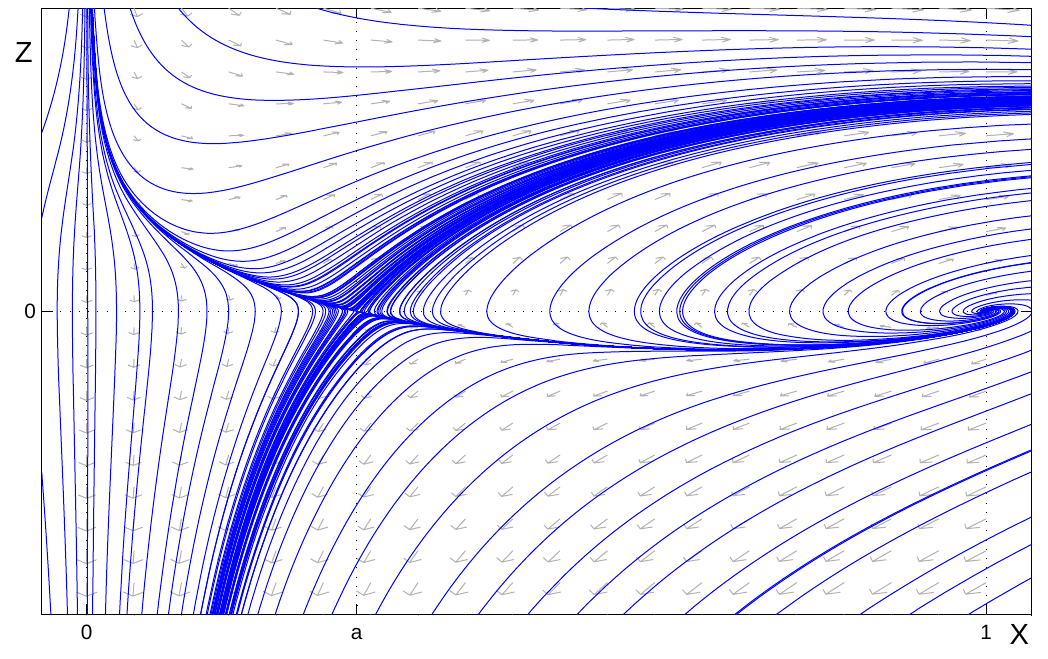}
  \includegraphics[scale =0.4]{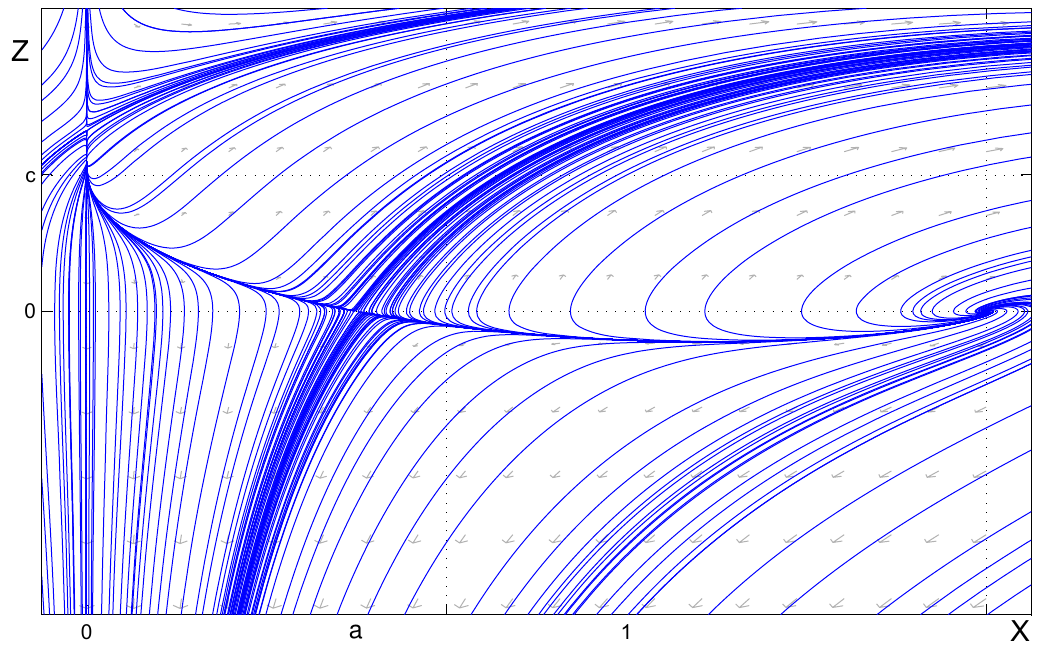}
  \includegraphics[scale =0.4]{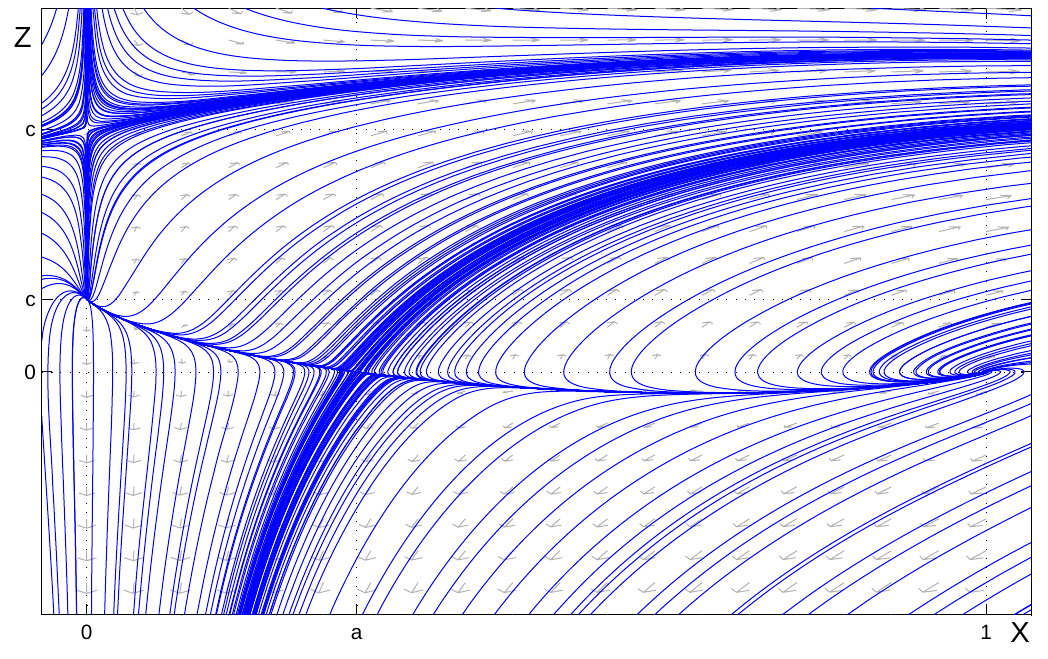}
  \caption{Reactions of type C', range $\gamma = 0$. Qualitative behaviour of the trajectories in the $(X,Z)$-plane for $f(u) = u(1-u)(a - u)$, $a = 0.3$, in the ranges $0< c < c_{\ast}$, $c = c_{\ast}$ and $c > c_{\ast}$, respectively.}\label{fig:QUALBEHAVTRAJPMETYPECPRIMEGAMMA0}
\end{figure}

\emph{Step 4': The cases $0 < c < c_{\ast}$ and $c > c_{\ast}$}. If $0 < c < c_{\ast}$, there are not \emph{a-admissible} TW. The proof of this fact easily follows from the study of the \emph{null isoclines} and from the non existence of critical points on the $Z$-axis.
\\
To the other hand, at each $c > c_{\ast}$, it corresponds exactly one \emph{a-admissible} TW and it is \emph{positive}. This is proved by showing the existence of a trajectory $T_c$ linking $A(a,0)$ and $R_{\lambda_1}(0,\lambda_1)$ corresponding to an \emph{a-admissible positive} TW profile. Again we refer to the proof of Theorem 2.2 of \cite{AA-JLV:art} for all the technical details. See Figure \ref{fig:QUALBEHAVTRAJPMETYPECPRIMEGAMMA0} for a qualitative representation of the trajectories in the $(X,Z)$-plane. $\Box$

\subsection{Reactions of type C, range \texorpdfstring{\boldmath}{}\texorpdfstring{$\gamma \geq 0$}{gamma}. Analysis of some special trajectories}\label{FINALREMARKODEPART}
In the PDEs part (see Theorem \ref{ASYMPTOTICBEHAVIOURTHEOREMTYPEC}, Part (ii)) we will need to compare general solutions to problem \eqref{eq:ALLENCAHNPME} with specific barriers which will be essentially constructed using TWs studied in the proof of the previous theorem. In the case of reactions of type C (satisfying \eqref{eq:ASSUMPTIONSONTHEREACTIONTERMTYPECPME}), we will employ TW profiles $\varphi(\xi) = \varphi(x-ct)$ satisfying
\[
\varphi(0) = a + \delta, \qquad \varphi(\xi_0) = 0 = \varphi(\xi_1), \qquad \varphi'(\xi) \;
\begin{cases}
> 0, \quad \text{ if } \xi_0 \leq \xi < 0 \\
< 0, \quad \text{ if } 0 < \xi \leq \xi_1,
\end{cases}
\]
for all $0 \leq c < c_{\ast}$ and suitable $\delta > 0$ depending on $c$ (cfr. with Figure \ref{fig:QUALBEHAVTRAJPMETYPEC0} and \ref{fig:QUALBEHAVTRAJPMETYPEC}). These TWs have been called ``change sign'' TWs of type 2 (CS-TWs) and their existence was proved in \cite{AA-JLV:art}, Subsection 3.1). In particular, the fact that $\varphi(\xi_0) = 0 = \varphi(\xi_1)$ comes from the fact that
\begin{equation}\label{eq:CSTWNEARCSPOINTS}
Z(X) \sim \pm X^{-\frac{1}{p-1}},
\end{equation}
for any trajectory with $Z(X) \sim \pm\infty$ for $X \sim 0$. This can be seen from the \emph{equation of the trajectories} \eqref{eq:EQUATIONOFTHETRAJECTORIESHYP}:
\[
\frac{dZ}{dX} = \frac{cZ - |Z|^p - f_{m,p}(X)}{(p-1)X|Z|^{p-2}Z},
\]
since for $Z(X) \sim \pm\infty$ for $X \sim 0$ it holds
\[
\frac{Z}{|Z|^2}\frac{dZ}{dX} \sim - \frac{1}{p-1} \frac{1}{X}, \qquad \text{for } \; Z(X) \sim \pm\infty, \; X \sim 0,
\]
which gives \eqref{eq:CSTWNEARCSPOINTS} (the accurate analysis is done in Section 3.1 of \cite{AA-JLV:art}).
We stress that these profiles exist only for speeds $0 \leq c < c_{\ast}$ and for all $\underline{\delta}_c \leq \delta < 1-a$, where $\underline{\delta}_c > 0$ is suitably chosen depending on $c$. Note that, from the monotonicity of the trajectory $T_c = T_c(X)$ studied before and the analysis of the \emph{nullisoclines}, we have that $\underline{\delta}_c \to \underline{\delta}_0$, as $c \to 0$, for some $0 < \underline{\delta}_0 < 1-a$. The fact that $\underline{\delta}_0 > 0$ is very important in the PDEs analysis.

In study of the so called ``threshold'' results for problem \eqref{eq:ALLENCAHNPME}-\eqref{eq:ASSUMPTIONSONTHEREACTIONTERMTYPECPME})(see Theorem \ref{ASYMPTOTICBEHAVIOURTHEOREMTYPEC}, Part (i)) will be employed other two important families of TWs, found in the ODEs analysis. The first one is composed by TW profiles $\varphi(\xi) = \varphi(x-ct)$ with the following properties:
\[
\varphi(0) = 1-\varepsilon, \qquad \varphi(\xi_0) = 0, \qquad \varphi(\xi_1) = a, \qquad \varphi(\xi) > a, \quad \text{for all } \; 0 < \xi < \xi_1,
\]
for some $\xi_0 < 0$, $\xi_1 > 0$, $\varepsilon > 0$ small and $c \geq c_{\ast}$. The property $\varphi(\xi_0) = 0$ is obtained exploiting again the fact that $Z(X) \sim \pm X^{-\frac{1}{p-1}}$, for $X \sim 0$, $Z \sim \pm \infty$, in the $(X,Z)$-plane, while the others come from the analysis of the null isoclines
(cfr. with Figure \ref{fig:QUALBEHAVTRAJPMETYPEC}). We will call them ``0-to-a'' TWs. As always, for any $c \geq c_{\ast}$, we can consider their reflections $\psi(\xi) = \psi(x + ct)$ satisfying
\[
\psi(0) = 1-\varepsilon, \qquad \psi(\xi_0) = a, \qquad \psi(\xi_1) = 0, \qquad \psi(\xi) > a, \quad \text{for all } \;\xi_0 < \xi < 0,
\]
for some $\xi_0 < 0$, $\xi_1 > 0$, $\varepsilon > 0$ small, to which will refer as ``a-to-0'' TW.

\subsection{Reactions of type C', range \texorpdfstring{\boldmath}{} \texorpdfstring{$\gamma \geq 0$}{gamma}. Analysis of some special trajectories}\label{FINALREMARKODEPARTCPRIME}
For what concerns the reactions of type C' (satisfying \eqref{eq:ASSUMPTIONSONTHEREACTIONTERMTYPEDPME}), we will consider TW profiles $\varphi(\xi) = \varphi(x-ct)$ with the following properties:
\begin{equation}\label{eq:INCREASINGATO1TWS}
\varphi(-\infty) = a, \qquad \varphi(\xi_0) = 1, \qquad \varphi'(\xi) > 0 \quad \text{ for all } \xi \leq \xi_0,
\end{equation}
where $0 < a < 1$ with $f(a) = 0$, $\xi_0 \in \RR$ is suitably chosen, and $c > 0$. The existence of these TW profiles follows from the analysis in the $(X,Z)$-plane (see part (i) and (ii) of Theorem \ref{THEOREMEXISTENCEOFTWSPMEREACTIONTYPEC}). Indeed, the study of the null isoclines and local behaviour of the critical point $A(a,0)$ show the existence of two trajectories ``coming from'' $A(a,0)$ and crossing the line $X=1$ in the $(X,Z)$-plane. The first one, lying in the strip $[a,1]\times(-\infty,0]$ satisfies \eqref{eq:INCREASINGATO1TWS}. The second one, lying in $[a,1]\times[0,+\infty)$, has symmetric properties but less significative for our purposes. We will call ``increasing a-to-1'' TWs the profiles satisfying \eqref{eq:INCREASINGATO1TWS}. These special solutions and their reflections will be used to prove that solutions to problem \eqref{eq:ALLENCAHNPME} converge to the steady state $u = a$ as $t \to +\infty$.

Finally, we point out that there are CS-TWs of type 2 in this setting too, but now they satisfy
\[
\varphi(0) = \delta, \qquad \varphi(\xi_0) = 0 = \varphi(\xi_1), \qquad \varphi'(\xi) \;
\begin{cases}
> 0, \quad \text{ if } \xi_0 \leq \xi < 0 \\
< 0, \quad \text{ if } 0 < \xi \leq \xi_1,
\end{cases}
\]
where $0 < \overline{\delta}_0 \leq \delta < a$, $0 < c < c_{\ast}$ and suitable $\xi_0 < \xi_1$, and $\underline{\delta}_0 > 0$. Their existence follows by analysis in the $(X,Z)$-plane or, as always, recalling the scaling property that links problem \eqref{eq:ALLENCAHNPME} with reaction of Fisher-KPP type to the one with reaction of type C' (cfr. with \cite{AA-JLV:art} subsection 3.1).

\paragraph{Continuity of $\boldsymbol{c_{\ast} = c_{\ast}(m,p,f)}$.} We end this section by stressing that the critical speed of propagation $c_{\ast} = c_{\ast}(m,p,f)$ is continuous in the region w.r.t. $(m,p)$ in the region $\{(m,p): m(p-1)-1 = \gamma \geq 0\}$, both for reactions of type C and of type C'. It can be proved using the methods of the Fisher-KPP setting adapting the proof given for the Fisher-KPP case (see Theorem 2.3, Corollary 3.2, and Lemma 4.1 of \cite{AA-JLV:art}).
%
%
%
%
%
%
%
%
%
%
%
%
\section{Reactions of Type C: threshold results and asymptotic behaviour}\label{THRESHOLDREACTIONSOFTYPEC}
This section is devoted to the proof of Theorem \ref{ASYMPTOTICBEHAVIOURTHEOREMTYPEC}, which is concerned on the asymptotic behaviour of solutions to problem \eqref{eq:ALLENCAHNPME} with initial data satisfying \eqref{eq:ASSUMPTIONSONTHEINITIALDATUMBISTABLE} and reaction terms of type C (satisfying \eqref{eq:ASSUMPTIONSONTHEREACTIONTERMTYPECPME}) and, as anticipated, on the stability/instability of the steady states $u = 0$, $u = 1$, depending on the initial data. Thus, before starting with the proof, we introduce two classes of initial data which generate solutions to problem \eqref{eq:ALLENCAHNPME} evolving to $u=0$ or $u=1$, respectively.
\begin{definition}\label{NONREACTIONGINITIALDATA}
We divide this definition depending on the dimension $N=1$ or $ N \geq 2$.

\noindent $\bullet$ Let $N = 1$. An initial data $u_0 = u_0(x)$ satisfying \eqref{eq:ASSUMPTIONSONTHEINITIALDATUMBISTABLE} is called ``not-reacting'' if there are $c_1, c_2 \geq c_{\ast}$ such that
\[
u_0(x) \leq \min \{\overline{\varphi},\overline{\psi}\}(x), \quad \text{ for all } x \in \RR,
\]
where
\[
\overline{\varphi}(\xi) :=
\begin{cases}
0 \quad & \text{if } \xi \leq \xi_0^{c_1} \\
\varphi_{c_1}(\xi) \quad & \text{if } \xi_0^{c_1} < \xi < \xi_1^{c_1}  \\
a \quad & \text{if } \xi \geq \xi_1^{c_1}
\end{cases}
\qquad
\overline{\psi}(\xi) :=
\begin{cases}
a \quad & \text{if } \xi \leq \xi_0^{c_2} \\
\psi_{c_2}(\xi) \quad & \text{if } \xi_0^{c_2} < \xi < \xi_1^{c_2}  \\
0 \quad & \text{if } \xi \geq \xi_1^{c_2},
\end{cases}
\]
and $\varphi_{c_1} = \varphi_{c_1}(x - ct)$ is a ``0-to-a'' TW corresponding to $c_1$ and $\psi_{c_2} = \psi_{c_2}(x + ct)$ is a ``a-to-0'' TW corresponding to $c_2$ (see Subsection \ref{FINALREMARKODEPART}).

\noindent $\bullet$ Let $N \geq 2$. An initial data $u_0 = u_0(x)$ satisfying \eqref{eq:ASSUMPTIONSONTHEINITIALDATUMBISTABLE} is called ``not-reacting'' if
\[
u_0(x) \leq \widetilde{u}_0(|x|), \quad \text{ for all } x \in \RR^N,
\]
where $\widetilde{u}_0 = \widetilde{u}_0(y)$ ($y \in \RR$) is a radial ``not-reacting'' initial datum in $N=1$.
\end{definition}

\begin{definition}\label{REACTIONGINITIALDATA}
Again we separate the cases $N=1$ or $ N \geq 2$.

\noindent $\bullet$ Let $N = 1$. An initial data $u_0 = u_0(x)$ satisfying \eqref{eq:ASSUMPTIONSONTHEINITIALDATUMBISTABLE} is called ``reacting'' if there is $0 < \widetilde{c} < c_{\ast}$ such that for all $0 \leq c \leq \widetilde{c}$, it holds
\[
u_0(x) \geq \max\{\underline{\varphi},\underline{\psi}\}(x),  \quad \text{ for all } x \in \RR,
\]
where
\[
\underline{\varphi}(\xi) :=
\begin{cases}
\varphi_c(\xi) \quad & \text{if } \xi_0^c \leq \xi \leq \xi_1^c \\
0              \quad &  \text{otherwise}
\end{cases}
\qquad
\underline{\psi}(\xi) :=
\begin{cases}
\psi_c(\xi) \quad & \text{if } \xi_0^{c'} \leq \xi \leq \xi_1^{c'} \\
0              \quad &  \text{otherwise}
\end{cases}
\]
and $\varphi_c = \varphi_c(x-ct)$ is a ``change-sign'' TW (of type 2) corresponding to $c$ and $\psi_c = \psi_c(x+ct)$ is its reflection  (see Subsection \ref{FINALREMARKODEPART}).

\noindent $\bullet$ Let $N \geq 2$. An initial data $u_0 = u_0(x)$ satisfying \eqref{eq:ASSUMPTIONSONTHEINITIALDATUMBISTABLE} is called ``reacting'' if there is $0 < c < c_{\ast}$ such that it holds
\[
u_0(x) > \underline{\varphi}(|x|-c\overline{t}) \quad \text{ for all } x \in \RR^N,
\]
where
\[
\underline{\varphi}(\xi) :=
\begin{cases}
\varphi(0)              \quad & \text{if } \xi \leq 0 \\
\varphi(\xi) \quad & \text{if } 0 \leq \xi \leq \xi_1 \\
0              \quad &  \text{otherwise}
\end{cases}
\qquad \text{ with } \varphi(0) = \max \varphi(\xi),
\]
$\overline{t} > 0$ is large enough and $\varphi = \varphi(\xi)$ is a ``change-sign'' TW (of type 2) corresponding to $c$ (the minimum value of the time $\overline{t} > 0$ will be specified in the proof of Theorem \ref{ASYMPTOTICBEHAVIOURTHEOREMTYPEC}).
\end{definition}
We are now ready to prove Theorem \ref{ASYMPTOTICBEHAVIOURTHEOREMTYPEC}.

\paragraph{Proof of Theorem \ref{ASYMPTOTICBEHAVIOURTHEOREMTYPEC}: Part (i).} We take a ``non-reacting'' initial datum $u_0 = u_0(x)$ (see Definition \ref{NONREACTIONGINITIALDATA}) and we prove that the solution $u = u(x,t)$ to problem \eqref{eq:ALLENCAHNPME} satisfies
\[
u(x,t) \to 0 \text{ uniformly in } \RR^N, \quad \text{ as } t \to +\infty.
\]
Let us firstly consider the case $N=1$. Since $u_0= u_0(x)$ is ``not-reacting'' there are $c_1,c_2 \geq c_{\ast}$ such that $u_0(x) \leq \min \{\overline{\varphi},\overline{\psi}\}(x)$ in $\RR$, as in Definition \ref{NONREACTIONGINITIALDATA}.

\noindent Note that both $\overline{\varphi}(\xi) = \overline{\varphi}(x-c_1t)$ and $\overline{\psi}(\xi) = \overline{\psi}(x+c_2t)$ are solutions to the equation in \eqref{eq:ALLENCAHNPME}, and at time $t = 0$, we have $u_0(x) \leq \overline{\varphi}(x)$ and $u_0(x) \leq \overline{\psi}(x)$ for all $x \in \RR$. Consequently, from the Comparison Principle we deduce $u(x,t) \leq \overline{\varphi}(x-c_1t)$ and $u(x,t) \leq \overline{\psi}(x+c_2t)$ for all $x \in \RR$ and $t > 0$, and, since $\overline{\varphi}(x-c_1t) = 0$ for all $x \leq \xi_0^{c_1}+c_1t$ and $\overline{\psi}(x+c_2t) = 0$ for all $x \geq \xi_1^{c_2}-c_2t$, we deduce that there is a time $t_{c_1,c_2} > 0$, such that $u(x,t) = 0$ for all $t \geq t_{c_1,c_2}$. This conclude the proof for the case $N=1$.

\noindent Before moving forward, we show that if $N=1$ and $u = u(x,t)$ is a solution to
\[
\begin{cases}
\partial_t u = \partial_x\left(|\partial_xu^m|^{p-2}\partial_xu^m\right) + f(u) \quad &\text{in } \RR\times(0,\infty) \\
u(x,0) = u_0(x) \quad &\text{in } \RR,
\end{cases}
\]
with initial data $u_0 = u_0(x)$ satisfying \eqref{eq:ASSUMPTIONSONTHEINITIALDATUMBISTABLE}, $u_0(x) = u_0(-x)$ and $u_0(\cdot)$ non-increasing for all $x \geq 0$, then $u(\cdot,t)$ is non-increasing w.r.t. $x \geq 0$, for all $t > 0$. So, fix $h > 0$ and let $v = v(x,t)$ be the solution to the problem
\[
\begin{cases}
\partial_t v = \partial_x\left(|\partial_xv^m|^{p-2}\partial_xv^m\right) + f(v) \quad &\text{in } \RR\times(0,\infty) \\
v(x,0) = v_0(x) := u_0(x+h) \quad &\text{in } \RR.
\end{cases}
\]
Hence, since $v_0(x) \leq u_0(x)$ we deduce $v(x,t) \leq u(x,t)$ and, by uniqueness of the solutions, it follows $v(x,t) = u(x+h,t)$. Hence, we obtain that $u(\cdot,t)$ is non-increasing for all $t \geq 0$ thanks to the arbitrariness of $x \geq 0$ and $h \geq 0$.

Now, assume $N \geq 2$ and consider radial solutions to problem  \eqref{eq:ALLENCAHNPME}, i.e., solutions $u = u(r,t)$ to the problem
\begin{equation}\label{eq:RADIALPROBLEMTYPECPRIME}
\begin{cases}
\partial_t u = \partial_r\left(|\partial_ru^m|^{p-2}\partial_ru^m\right) + \frac{N-1}{r}|\partial_ru^m|^{p-2}\partial_ru^m + f(u) \quad &\text{in } \RR_+\times(0,\infty) \\
u(r,0) = u_0(r) \quad &\text{in } \RR_+\times\{0\},
\end{cases}
\end{equation}
where $r = |x|$, $x \in \RR^N$, and $u_0(\cdot)$ is a radially decreasing ``not-reacting'' initial datum. Moreover, let $\overline{u} = \overline{u}(r,t)$ be a solution to the problem
\[
\begin{cases}
\partial_t \overline{u} = \partial_r\left(|\partial_r\overline{u}^m|^{p-2}\partial_r\overline{u}^m\right) + f(\overline{u}) \quad &\text{in } \RR_+\times(0,\infty) \\
\overline{u} = u \quad &\text{in } \{0\}\times(0,\infty) \\
\overline{u}(r,0) = u_0(r) \quad &\text{in } \RR_+\times\{0\}.
\end{cases}
\]
For what explained before, we have $\partial_r u(r,t) \leq 0$ in $\RR_+\times(0,\infty)$, and so $\overline{u} = \overline{u}(r,t)$ is a super-solution to \eqref{eq:RADIALPROBLEMTYPECPRIME}. But $\overline{u} = \overline{u}(r,t)$ is a solution of the one-dimensional equation with ``not-reacting'' initial data, and so, from the case $N=1$, it follows
\[
\overline{u}(r,t) = 0 \; \text{ in } \RR_+\cup\{0\}, \quad \text{ for all } t \geq t_{c_1,c_2},
\]
and by the comparison, we deduce the same for $u = u(r,t)$, concluding the proof of Part (i). $\Box$

\paragraph{Proof of Theorem \ref{ASYMPTOTICBEHAVIOURTHEOREMTYPEC}: Part (ii), case $\boldsymbol{N=1}$.} We fix $N =1$ and we proceed in two steps.

\emph{Step 1: Propagation of minimal super-level sets.} Consider a ``reacting'' initial datum $u_0 = u_0(x) \geq \max \{\underline{\varphi},\underline{\psi}\}(x)$, $x \in \RR$, where $\underline{\varphi}(\xi) = \underline{\varphi}(x-ct)$ and $\underline{\psi}(\xi) = \underline{\psi}(x+ct)$ are defined in Definition \ref{REACTIONGINITIALDATA}. From the ODEs analysis of section \ref{EXISTENCEOFTWSTYPECDPME} and from what explained in Subsection \ref{FINALREMARKODEPART}, we can assume that for all $0 \leq c \leq \widetilde{c}$, it holds
\[
\underline{\varphi}(0) = \underline{\psi}(0) = a + \underline{\delta}_c < 1, \qquad \underline{\delta}_c \geq \underline{\delta}_0,
\]
for some $0 < \underline{\delta}_0 < 1-a$. Now, since $u_0(x) \geq \underline{\varphi}(x)$ and $u_0(x) \geq \underline{\psi}(x)$, we deduce by comparison $u(x,t) \geq \underline{\varphi}(x-ct)$ and $u(x,t) \geq \underline{\psi}(x+ct)$ for all $x \in \RR^N$ and $t >0$. Hence, by the arbitrariness of $0 \leq c \leq \widetilde{c}$, we obtain that
\[
u(x,t) \geq a + \underline{\delta}_0 \quad \text{in } \{|x|\leq \widetilde{c} t\}, \; \text{ for all } t > 0.
\]
\emph{Step 2: Convergence to 1 on compact sets.} Now, fix $\varepsilon > 0$ small and $\widetilde{\varrho} > 0$ arbitrarily large. Then, we have
\[
u(x,t) \geq a + \underline{\delta}_0 \quad \text{in } \{|x| \leq \widetilde{\varrho}\}, \; \text{ for all } t \geq t_{\widetilde{\varrho},\widetilde{c}} := \widetilde{\varrho}/\widetilde{c},
\]
which implies
\begin{equation}\label{eq:LINEARBOUNDFROMBELOWFORFDELTA0}
f(u) \geq q(1-u) \quad \text{in } \{|x| \leq \widetilde{\varrho}\}\times[t_{\widetilde{\varrho},\widetilde{c}},\infty),
\end{equation}
for some suitable $q = q_{\underline{\delta}_0} > 0$ (cfr. with the remark at the end of the proof). Thus, the solution $\underline{u} = \underline{u}(x,t)$ to the problem
\begin{equation}\label{eq:SUBSOLUTIONBOUNDEDDOMAINTYPECN1}
\begin{cases}
\partial_t\underline{u} = \partial_x\left(|\partial_x\underline{u}^m|^{p-2}\partial_x\underline{u}^m\right) + q_{\underline{\delta}_0}(1-\underline{u})  \quad & \text{in } \{|x| \leq \widetilde{\varrho}\}\times[t_{\widetilde{\varrho},\widetilde{c}},\infty) \\
\underline{u}(x,t) = a + \underline{\delta}_0 \quad & \text{in } \partial\{|x| \leq \widetilde{\varrho}\}\times[t_{\widetilde{\varrho},\widetilde{c}},\infty) \\
\underline{u}(x,t_{\widetilde{\varrho},\widetilde{c}}) = a + \underline{\delta}_0 \quad & \text{in } \{|x| \leq \widetilde{\varrho}\}
\end{cases}
\end{equation}
is a sub-solution to problem \eqref{eq:ALLENCAHNPME} in $\{|x| \leq \widetilde{\varrho}\}\times[t_{\widetilde{\varrho},\widetilde{c}},\infty)$ and so, by the Comparison Principle, we obtain $\underline{u}(x,t) \leq u(x,t)$ in $\{|x| \leq \widetilde{\varrho}\}\times[t_{\widetilde{\varrho},\widetilde{c}},\infty)$. Now, following the proof of Lemma 7.1 of \cite{AA-JLV:art}, it is not difficult to see that $a + \underline{\delta}_0 \leq \underline{u}(x,t) \leq 1$ and the function $t \to \underline{u}(x,t)$ is non-decreasing for any $x \in\RR$, so that there exists its uniform limit $\underline{u}_{\infty}(x) := \lim_{t \to +\infty} \underline{u}(x,t)$ which solves
\begin{equation}\label{eq:ELLIPTICCOMPARISONSUBSOLUTIONASYMPTOTICBEHAVIOURPGREATER2}
\begin{cases}
- \Delta_p \underline{u}_{\infty}^m = q(1-\underline{u}_{\infty}) \quad &\text{in } \{|x| \leq \widetilde{\varrho}\} \\
\underline{u}_{\infty} = a + \underline{\delta}_0  \quad &\text{in } \partial\{|x| \leq \widetilde{\varrho}\}.
\end{cases}
\end{equation}
Now, let us define
\[
w^m(r) := A \left[e^{g(r)} - 1\right], \qquad g(r) := 1 - \left(\frac{r}{\widetilde{\varrho}}\right)^{\lambda}, \qquad  \lambda:= \frac{p}{p-1},
\]
where $r = |x|$ and $(1-\varepsilon/2)^m (e-1)^{-1} < A < (e-1)^{-1}$. It is immediate to check that $w = w(|x|)$ is well-defined in $\{|x| \leq \widetilde{\varrho}\}$, $w = 0$ in $\partial\{|x| \leq \widetilde{\varrho}\}$. In Lemma 7.1 of \cite{AA-JLV:art} it is proved that $w = w(|x|)$ is a sub-solution to the equation in \eqref{eq:ELLIPTICCOMPARISONSUBSOLUTIONASYMPTOTICBEHAVIOURPGREATER2}, i.e., $-\Delta_{p,r}w^m \leq q (1-w)$ in $\{|x| \leq \widetilde{\varrho}\}$ whenever $\widetilde{\varrho} > 0$ is large enough, namely
\[
\widetilde{\varrho}^{\,p} \geq \widetilde{\varrho}_{\varepsilon}^{\,p} := \frac{(Ae\lambda)^{p-1}}{q\{1 - [A(e-1)]^{1/m} \}}.
\]
Furthermore, it follows that $w(r) \geq 1 - \varepsilon/2$ in $\{r \leq \widetilde{A}_{\varepsilon}\widetilde{\varrho}\}$, where
\[
\widetilde{A}_{\varepsilon}^{\,\lambda} := 1 - \log\left[\frac{A + (1 - \varepsilon/2)^m}{A} \right],
\]
and, by the Elliptic Comparison Principle, we obtain $\underline{u}_{\infty}(x) \geq 1 - \varepsilon/2$ in $\{|x| \leq \widetilde{A}_{\varepsilon}\widetilde{\varrho}\}$. Consequently, by uniform convergence $\underline{u}(\cdot,t) \to \underline{u}_{\infty}(\cdot)$, we deduce the existence of a time $t_1 > 0$ large enough such that $\underline{u}(x,t) \geq 1-\varepsilon$ in $\{|x| \leq \widetilde{A}_{\varepsilon} \widetilde{\varrho}\}$ for all $t \geq t_1$. Since $\widetilde{\varrho} > 0$ can be taken larger and recalling that $u(x,t) \geq \underline{u}(x,t)$ for all $t \geq t_{\widetilde{p},\widetilde{c}}$, the thesis follows. $\Box$
\paragraph{Remark.} Note that in \cite{AA-JLV:art}, the authors worked with a \emph{concave} Fisher-KPP reaction term (cfr. with formula (1.2) of \cite{AA-JLV:art}). However, the concavity assumption on $f(\cdot)$ is not used to prove Lemma 7.1 of \cite{AA-JLV:art} which we basically apply in the above proof. As we said, for our purposes it is enough to prove the claim:
\[
\text{For all } 0 < \delta < 1-a, \text{ there exists } q_{\delta} > 0 \text{ such that } f(u) \geq q_{\delta}(1-u), \text{ for all } a + \delta \leq u \leq 1,
\]
from which \eqref{eq:LINEARBOUNDFROMBELOWFORFDELTA0} it is easily deduced. To see this, let us fix $0 < \delta < 1-a$  and take
\[
q_{\delta} := \frac{f(a+\delta/n)}{1-a-\delta} > 0,
\]
for some integer $n \geq 1$ large enough. Assume by contradiction that for any integer $n \geq 1$, there exist $u_n \in (a+\delta,1)$, such that
\[
f(u_n) < \frac{f(a + \delta/n)}{1-a-\delta}(1-u_n).
\]
Since $\{u_n\}_n$ is a bounded sequence we can assume it converges to some limit $l \in [a+\delta,1]$ up to passing to a suitable subsequence that we rename $u_n$. If $l \in [a+\delta,1)$, taking the the limit as $n \to +\infty$ in the above inequality to obtain that the l.h.s. converges to $f(l) > 0$, while the r.h.s. converges to zero since $f(a) = 0$. If $u_n \to 1$, the same argument hold using the fact that $f(u_n) \sim |f'(1)|(1-u_n)$ as $n \to +\infty$ and $f'(1) < 0$ by assumption. This gives the desired contradiction and proves the claim.

\paragraph{Proof of Theorem \ref{ASYMPTOTICBEHAVIOURTHEOREMTYPEC}: Part (ii), case $\boldsymbol{N \geq 2}$.} We fix $N \geq 2$ and, proceeding as in part (i), we consider the radial problem
\[
\begin{cases}
\partial_t u = \partial_r\left(|\partial_ru^m|^{p-2}\partial_ru^m\right) + \frac{N-1}{r}|\partial_ru^m|^{p-2}\partial_ru^m + f(u) \quad &\text{in } \RR_+\times(0,\infty) \\
u(r,0) = u_0(r) \quad &\text{in } \RR_+\times\{0\},
\end{cases}
\]
where $r = |x|$, $x \in \RR^N$, and $u_0(\cdot)$ is a radially decreasing ``reacting'' initial datum. By definition, we can assume that for any fixed $\varepsilon > 0$ (small), there is $0 < c < c_{\ast}$, such that
\[
u_0(r) > \underline{\varphi}(r -(c+\varepsilon)\overline{t}) \quad \text{ for all } r > 0,
\]
where $\underline{\varphi} = \underline{\varphi}(\xi)$ is as in Definition \ref{REACTIONGINITIALDATA} part (ii) and $\overline{t} > 0$ is large enough and will be chosen later. Now, setting $\delta := 1-\varepsilon$, we define
\[
\begin{cases}
\underline{u}(r,t) = \underline{\varphi}(\delta^{1/p}r - c\delta t) \quad &\text{ if } m > 1, \; p > 1 \; (\gamma > 0) \\
\underline{u}(r,t) = \underline{\varphi}(r - c\delta t)             \quad &\text{ if } 0 < m < 1, \; p > 2 \; (\gamma > 0),
\end{cases}
\]
as in the proof of Theorem 2.6 of \cite{AA-JLV:art} (see Section 9). Repeating the same procedure of that proof, it is easily seen that $\underline{u} = \underline{u}(r,t)$ is a sub-solution to the equation in problem \eqref{eq:RADIALPROBLEMTYPECPRIME} in $\RR_+\times[\widetilde{t},\infty)$ where $\widetilde{t} > 0$ is suitably chosen (large enough). Consequently, taking $\overline{t} := \widetilde{t} \gg 0$ and noting that we can assume $u_0(r) \geq \underline{u}(r,\overline{t})$ for any $r > 0$ (since $\delta = 1-\varepsilon$ and $\varepsilon > 0$ is arbitrarily small), we conclude by comparison
\[
u(r,t) \geq \underline{u}(r,t + \overline{t}) \geq \varphi(0) = a + \underline{\delta}_0 \quad \text{in } \{r = |x|\leq \widetilde{c} t\}, \; \text{ for all } t > 0,
\]
for some $\underline{\delta}_0 > 0$ and $\widetilde{c} := c\delta$ (see also Lemma 5.1 of \cite{Aro-Wein2:art} for the linear setting). Moreover, exactly as in \emph{Step 2} of the case $N=1$, we have
\[
f(u) \geq q(1-u) \quad \text{in } \{|x| \leq \widetilde{\varrho}\}\times[t_{\widetilde{\varrho},\widetilde{c}},\infty),
\]
for some suitable $q = q_{\underline{\delta}_0} > 0$, all $\widetilde{\varrho} > 0$ and $t_{\widetilde{\varrho},\widetilde{c}} >0$ large enough, and so we can repeat that construction to show that for all $\varepsilon > 0$ (small) and $\widetilde{\varrho} > 0$ (large), there exists $t_1 > 0$ such that
\[
u(r,t) \geq 1 - \varepsilon \quad \text{ in } \{r = |x| \leq \widetilde{a}_{\varepsilon}\widetilde{\varrho}\} \; \text{ for all } t \geq t_1,
\]
where $0 < \widetilde{a}_{\varepsilon} < 1$ is as in the case $N=1$, concluding the proof of the case $N \geq 2$. $\Box$
\paragraph{Remark.} The above proof strongly relies on the fact that the function
\[
\begin{cases}
\underline{u}(r,t) = \underline{\varphi}(\delta^{1/p}r - c\delta t) \quad &\text{ if } m > 1, \; p > 1 \; (\gamma > 0) \\
\underline{u}(r,t) = \underline{\varphi}(r - c\delta t)             \quad &\text{ if } 0 < m < 1, \; p > 2 \; (\gamma > 0),
\end{cases}
\]
defined depending on the value $m > 0$ and $p > 1$ such that $\gamma > 0$, is a sub-solution to problem \eqref{eq:RADIALPROBLEMTYPECPRIME} for large times $t \gg 0$. As we have mentioned before, this fact can be easily showed by repeating the proof of Theorem 2.6 of \cite{AA-JLV:art} (cfr. with Section 9) for the Fisher-KPP framework. This parallelism is due to the fact that the main difficulty is to study the sign of the quantity
\[
\partial_t \underline{u} - \partial_r\left(|\partial_r\underline{u}^m|^{p-2}\partial_r\underline{u}^m\right) - \frac{N-1}{r}|\partial_r\underline{u}^m|^{p-2}\partial_r\underline{u}^m - f(\underline{u})
\]
near the points in which $\underline{u} = 0$, i.e. $\varphi = 0$ (here $\varphi$ denotes the profile of a ``change-sign'' TW, cfr. with Subsection \ref{FINALREMARKODEPART}). The behaviour of $\varphi$ near the ``change-sign'' points is completely understood and is the same for both reactions of type C and Fisher-KPP reactions (cfr. with \ref{FINALREMARKODEPART} and Subsection 3.1 of \cite{AA-JLV:art}).

\paragraph{Proof of Theorem \ref{ASYMPTOTICBEHAVIOURTHEOREMTYPEC}: Part (iii).} Let us prove that for all radially decreasing initial data $u_0 = u_0(x)$ satisfying \eqref{eq:ASSUMPTIONSONTHEINITIALDATUMBISTABLE} and for all $c > c_{\ast}(m,p,f)$ it holds
\[
u(x,t) \to 0 \text{ uniformly in } \{|x| \geq ct\}, \quad \text{ as } t \to +\infty.
\]
This part is the easiest and it actually coincides with Proposition 8.1 (for the case $N=1$) and Theorem 2.6 (for $N \geq 2$) of \cite{AA-JLV:art}. Here we just explain the main ideas and we refer the reader to the just mentioned references for the details. If $N=1$, we fix $c > c_{\ast}$, $\varepsilon > 0$, and we consider the functions
\[
\overline{v}(x,t) := \varphi(x - c_{\ast}t), \qquad\qquad \overline{w}(x,t) := \psi(x + c_{\ast}t),
\]
where $\varphi = \varphi(\xi)$ is the \emph{finite admissible} TW studied in Theorem \ref{THEOREMEXISTENCEOFTWSPMEREACTIONTYPEC}, part (i), with its reflection $\psi(\xi) = \psi(x + c_{\ast}t)$. Since $u_0 = u_0(x)$ satisfies \eqref{eq:ASSUMPTIONSONTHEINITIALDATUMBISTABLE}, we can assume $u_0(x) \leq \varphi(x)$ and $u_0(x) \leq \psi(x)$ for all $x \in \RR$, and so, thanks to the Comparison Principle, we  obtain both $u(x,t) \leq \overline{v}(x,t)$ and $u(x,t) \leq \overline{w}(x,t)$. Thus, since $\overline{v}(x,t) \leq \varepsilon$ for $x \geq c_{\ast}t + \xi_0$ and $\overline{w}(x,t) \leq \varepsilon$ for $x \leq -c_{\ast}t + \xi_0$ and $c > c_{\ast}$, we deduce that $u(x,t) \leq \varepsilon$ in $\{|x| \geq ct\}$ for $t > 0$ large enough.
\\
We point out that if $\gamma > 0$ then $\overline{v}(x,t) = 0$ for $x \geq c_{\ast}t + \xi_0$ and $\overline{w}(x,t) = 0$ for $x \leq -c_{\ast}t + \xi_0$ which implies that $u = u(x,t)$ has a \emph{free boundary}, too, whilst this does not happen when $\gamma = 0$, since the TW solutions are positive everywhere.

When $N \geq 2$, we follow the proof of Part (i), using that solutions to problem \eqref{eq:ALLENCAHNPME} with $N=1$ are super-solution for radial solutions of the same problem and so, by comparison, the thesis follows.

Now, we show that for all ``reacting'' initial data $u_0 = u_0(x)$ and for all $0 < c < c_{\ast}(m,p,f)$, it holds
\[
u(x,t) \to 1 \text{ uniformly in } \{|x| \leq ct\}, \quad \text{ as } t \to +\infty.
\]
Let us consider the case $N=1$. From part (ii) we obtain that for all $\varepsilon > 0$, $\widetilde{\varrho} > 0$, and all ``reacting'' initial data $u_0 = u_0(x)$, there exist $t_1 > 0$, such that
\[
u(x,t) \geq 1 - \varepsilon \quad \text{ in } \{|x| \leq \widetilde{\varrho}\} \; \text{ for all } t \geq t_1.
\]
Hence, for all $0 \leq c < c_{\ast}$, taking eventually $\widetilde{\varrho} > 0$ larger, there is a ``change-sign'' TW $\underline{\varphi}(\xi) = \underline{\varphi}(x-ct)$, i.e. a wave solution with profile satisfying (cfr. with Subsection \ref{FINALREMARKODEPART})
\[
\underline{\varphi}(0) = 1-\varepsilon, \qquad \varphi(\xi_0) = 0 = \varphi(\xi_1), \qquad \varphi'(\xi) \;
\begin{cases}
> 0, \quad \text{ if } \xi_0 \leq \xi < 0 \\
< 0, \quad \text{ if } 0 < \xi \leq \xi_1,
\end{cases}
\]
for suitable $\xi_0 < 0 < \xi_1$, such that $u(x,t_1) \geq \underline{\varphi}(x)$ for all $x \in \RR$. Of course, we can assume that also its reflection $\underline{\psi}(\xi) = \underline{\psi}(x+ct)$ (with $\underline{\psi}(0) = 1 - \varepsilon$) satisfies $u(x,t_1) \geq \underline{\psi}(x)$ for all $x \in \RR$.  Consequently, by comparison we have $u(x,t_1+t) \geq \underline{\varphi}(x-ct)$ and $u(x,t_1+t) \geq \underline{\psi}(x+ct)$ for all $x \in \RR$ and $t >0$, and the level $1-\varepsilon$ propagates with speed $c$. Hence, using again the arbitrariness of $0 \leq c < c_{\ast}$, we deduce
\[
u(x,t) \geq 1 - \varepsilon \quad \text{ in } \{|x| \leq ct\} \; \text{ for all } t \geq t_2,
\]
for some $t_2 = t_2(\varepsilon,c)$ large enough. This shows our statement, since $\varepsilon > 0$ has been chosen arbitrarily small.

Finally, when $N \geq 2$, following the proof of part (ii), case $N \geq 2$ and using again the sub-solution constructed in the proof of Theorem 2.6 of \cite{AA-JLV:art} with speed $c < c_{\ast}$ and $\underline{\varphi}(0) = 1-\varepsilon$, we conclude as in the case $N=1$. $\Box$

%
%
%
%
%
%
%
%
%
%
%
%
\section{Reactions of Type C': asymptotic behaviour}\label{REACTIONSOFTYPECPRIME}
This section is devoted to the proof of Theorem \ref{ASYMPTOTICBEHAVIOURTHEOREMTYPECCPRIME} part (ii). We then consider reactions of type C', i.e. satisfying \eqref{eq:ASSUMPTIONSONTHEREACTIONTERMTYPEDPME}. As in the ODEs part, some of our proofs rely on the results obtained in \cite{AA-JLV:art} that can be recovered by scaling (see the beginning of the proof of Theorem \ref{THEOREMEXISTENCEOFTWSPMEREACTIONTYPEC}, Part (ii), range $\gamma >0$). We recall that, as always, $0 < a < 1$ satisfies $f(a) = 0$.

We proceed by proving Theorem \ref{ASYMPTOTICBEHAVIOURTHEOREMTYPECCPRIME} part (ii), taking spatial dimension $N=1$. The reduction to dimension $N=1$ is necessary to compare solutions $u = u(x,t)$ to problem \eqref{eq:ALLENCAHNPME}-\eqref{eq:ASSUMPTIONSONTHEINITIALDATUMBISTABLE}-\eqref{eq:ASSUMPTIONSONTHEREACTIONTERMTYPEDPME} with TW solutions studied in Section \ref{EXISTENCEOFTWSTYPECDPME}. As we will see in a moment, we construct two super-solutions to prove that $u = u(x,t)$ reaches the level $0 < a < 1$ in finite time and a third super-solution combined to a scaling technique, to show that $u = u(x,t)$ converges uniformly to zero in the ``outer sets'' $\{|x|\geq ct\}$ as $t \to +\infty$.
\paragraph{Proof of Theorem \ref{ASYMPTOTICBEHAVIOURTHEOREMTYPECCPRIME}: Case $\boldsymbol{N=1}$, range $\boldsymbol{\gamma > 0}$.} Fix $m > 0$ and $p > 1$ such that $\gamma > 0$. We begin with two preliminary steps, crucial in the rest of the proof.

\emph{Step 0.} We first prove that for all $\varepsilon > 0$, there exists a waiting time $t_{\varepsilon} > 0$ such that
\[
u(x,t) \leq a + \varepsilon, \quad \text{for all } x \in \RR, \quad t \geq t_{\varepsilon}.
\]
To do this we employ the ``increasing a-to-1'' TWs and their reflections, found in Theorem \ref{THEOREMEXISTENCEOFTWSPMEREACTIONTYPEC}, cfr. with Subsection \ref{FINALREMARKODEPARTCPRIME}. To be more specific, we fix $c = 1$ and we consider a TW profile $\varphi(\xi) = \varphi(x-t)$ moving toward the right direction, satisfying
\[
\varphi(-\infty) = a, \qquad \varphi(\xi_0) = 1, \qquad \varphi'(\xi) > 0 \quad \text{ for all } \xi \leq \xi_0,
\]
and its ``reflection'' $\psi(\xi) = \psi(x + t)$, moving toward the left direction, satisfying
\[
\psi(+\infty) = a, \qquad \psi(\xi_1) = 1, \qquad \psi'(\xi) < 0 \quad \text{ for all } \xi \geq \xi_1,
\]
for some $\xi_0,\xi_1 \in \RR$ (cfr with formula \eqref{eq:INCREASINGATO1TWS}). Defining
\[
\overline{\varphi}(\xi) :=
\begin{cases}
\varphi(\xi) \quad &\text{if } \xi \leq \xi_0 \\
1            \quad &\text{if } \xi \geq \xi_0,
\end{cases}
\qquad \qquad \qquad
\overline{\psi}(\xi) :=
\begin{cases}
1            \quad &\text{if } \xi \leq \xi_1 \\
\psi(\xi) \quad &\text{if } \xi \geq \xi_1,
\end{cases}
\]
and recalling that $u_0 \in \mathcal{C}_c(\RR)$ with $0 \leq u_0 \leq 1$ we can assume both $u_0(x) \leq \overline{\varphi}(x)$ and $u_0(x) \leq \overline{\psi}(x)$ for all $x \in \RR$.

\noindent Now, we fix $\varepsilon > 0$ small, such that $1-\varepsilon > 0$. Defining the function $\overline{v}(x,t) := \overline{\varphi}(x-(1-\varepsilon)t)$ and using the definition of $\overline{\varphi} = \overline{\varphi}(\xi)$, we get that
\[
\begin{aligned}
&\partial_t\overline{v} - \partial_x\left(|\partial_x\overline{v}^m|^{p-2}\partial_x\overline{v}^m\right) - f(\overline{v}) \\
& = -(1-\varepsilon)\overline{\varphi}' - \left[|\left(\overline{\varphi}^m\right)'|^{p-2}\left(\overline{\varphi}^m\right)'\right]' - f(\overline{\varphi}) = \varepsilon\overline{\varphi}' - \overline{\varphi}' - \left[|\left(\overline{\varphi}^m\right)'|^{p-2}\left(\overline{\varphi}^m\right)'\right]' - f(\overline{\varphi}) \\
& = \varepsilon\overline{\varphi}' \geq 0, \quad \text{for all } \xi \leq \xi_0,
\end{aligned}
\]
where $\overline{\varphi}'$ stands for the derivative of $\overline{\varphi}(\cdot)$ w.r.t. $\xi$. Note that when $\xi \geq \xi_0$, $\overline{v}(x,t) = 1$, i.e., it is just a stationary state of the equation in \eqref{eq:ALLENCAHNPME} and the equality holds in the last inequality for  $\xi \geq \xi_0$. In particular, it follows that the function $\overline{v} = \overline{v}(x,t)$ is a super-solution for the equation in \eqref{eq:ALLENCAHNPME}.

\noindent Similarly, one can define $\overline{w}(x,t) = \overline{\psi}(x+(1+\varepsilon)t)$ and prove it is a super-solution too. In this case the function $\overline{w} = \overline{w}(x,t)$ is wave moving toward the left direction.

\noindent Hence, thanks to the comparison principle and remembering that $u_0(x) \leq \overline{\varphi}(x)$ and $u_0(x) \leq \overline{\psi}(x)$, we deduce
\[
u(x,t) \leq \overline{v}(x,t), \quad \text{ and } \quad u(x,t) \leq \overline{w}(x,t), \quad \text{ in } \RR^N\times[0,\infty).
\]
Moreover, thanks to the properties of $\varphi = \varphi(\xi)$ and $\psi = \psi(\xi)$, we deduce the existence of $\xi_{\varepsilon} > 0$, such that
\[
\begin{aligned}
\overline{v}(x,t) \leq a + \varepsilon, \quad \text{ for all } x \leq -\xi_{\varepsilon} + (1-\varepsilon)t \\
\overline{w}(x,t) \leq a + \varepsilon, \quad \text{ for all } x \leq \xi_{\varepsilon} - (1+\varepsilon)t.
\end{aligned}
\]
Thus, we get $u(x,t) \leq a +\varepsilon$ in $\RR^N$ if
\[
-\xi_{\varepsilon} + (1-\varepsilon)t \geq \xi_{\varepsilon} - (1+\varepsilon)t,
\]
i.e. $t \geq t_{\varepsilon}:= \xi_{\varepsilon}$.

\emph{Step 1.} In this step, we construct a global super-solution to problem \eqref{eq:ALLENCAHNPME} to show that our solution $u = u(x,t)$ propagates with \emph{finite} speed of propagation, i.e., $u=0$ outside an interval of $\RR$ with radius expanding in time. Consider the solution to the problem
\[
\begin{cases}
\partial_t\overline{u} = \partial_x\left(|\partial_x\overline{u}^m|^{p-2}\partial_x\overline{u}^m\right) + f'(0)\overline{u} \quad &\text{in } \RR\times(0,\infty) \\
\overline{u}(x,0) = u_0(x) \quad &\text{in } \RR.
\end{cases}
\]
Observe that $u(x,t) \leq \overline{u}(x,t)$ in $\RR\times(0,\infty)$ since $f(u) \leq f'(0)u$, thanks to the first assumption in \eqref{eq:ASSUMPTIONSONTHEREACTIONTERMTYPEDPME} and the comparison principle. Furthermore, the function defined by
\[
\widetilde{u}(x,\tau) = e^{-f'(0)t}\overline{u}(x,t), \qquad \text{with } \; \tau(t) = \frac{1}{f'(0)\gamma}\left(e^{f'(0)\gamma t} - 1\right), \quad t \geq 0,
\]
satisfies the purely diffusive equation
\[
\begin{cases}
\partial_{\tau}\widetilde{u} = \partial_x\left(|\partial_x\widetilde{u}^m|^{p-2}\partial_x\widetilde{u}^m\right) \quad &\text{in } \RR\times(0,\infty) \\
\widetilde{u}(x,0) = u_0(x) \quad &\text{in } \RR.
\end{cases}
\]
Consequently, since $\widetilde{u} = \widetilde{u}(x,\tau)$ has \emph{finite} speed of propagation (see for instance \cite{V1:book,V2:book}), we deduce the same for $\overline{u} = \overline{u}(x,t)$, and so for $u = u(x,t)$. We conclude this step pointing out that the same procedure can be adapted (with obvious changes) to the case $N \geq 1$.

\emph{Step 2.} In this part of the proof, we show that for all $c > c_{\ast}(m,p,f)$, there exists $t_1 = t_1(c) > 0$ such that
\[
u(x,t) = 0, \quad \text{ in } \{|x| \geq ct\}, \quad \text{ for all } t \geq t_1.
\]
So, fix $\varepsilon > 0$ small and $c > c_{\ast}(m,p,f)$.  We assume for a moment that $0 \leq u(x,t_{\varepsilon}) < a$ for all $x \in \RR$, where $t_{\varepsilon} > 0$ is the one found in \emph{Step 0}. Moreover, we know that $u(x,t_{\varepsilon}) = 0$ outside an interval of $\RR$ of radius large enough (see \emph{Step 1}). Hence, we define
\[
\overline{v}(x,t) := \varphi(x - c_{\ast}t), \qquad\qquad \overline{w}(x,t) := \psi(x + c_{\ast}t),
\]
where $\varphi = \varphi(\xi)$ is the \emph{finite $a$-admissible} TW studied in Theorem \ref{THEOREMEXISTENCEOFTWSPMEREACTIONTYPEC} (part (ii), range $\gamma > 0$), satisfying $\varphi(-\infty) = a$, $\varphi(\xi) = 0$ for all $\xi \geq \xi_0$, and $\psi = \psi(\xi)$ is its ``reflection''. Consequently, up to a left-right shift, we can assume $u(x,t_{\varepsilon}) \leq \varphi(x)$ and $u(x,t_{\varepsilon}) \leq \psi(x)$, and so, by the comparison principle we deduce
\[
u(x,t + t_{\varepsilon}) \leq \overline{v}(x,t), \qquad u(x,t + t_{\varepsilon}) \leq \overline{w}(x,t), \quad \text{ in } \RR^N\times(0,\infty).
\]
Thus, since $\overline{v}(x,t) = 0$ for $x \geq c_{\ast}t + \xi_0$ and $\overline{w}(x,t) = 0$ for $x \leq -c_{\ast}t + \xi_0$ and $c > c_{\ast}$, we deduce that $u(x,t) = 0$ in $\{|x| \geq ct\}$ for large times exactly as in the proof of Theorem \ref{ASYMPTOTICBEHAVIOURTHEOREMTYPEC}, Part (iii).

\noindent Now, if $u(x,t_{\varepsilon}) \geq a$ for some $x \in \RR$, it must be $u(x,t_{\varepsilon}) \leq a + \varepsilon$ in $\RR$, from what proved in \emph{Step 0}. We consider the re-scaling of $u = u(x,t)$ defined by
\[
u_{\varepsilon}(y,\tau) = a^{-1}(a + 2\varepsilon)u(x,t), \quad \text{ where } \; y = \left[a^{-1}(a + 2\varepsilon)\right]^{\frac{m(p-1)}{p}}x, \quad \tau = a^{-1}(a + 2\varepsilon)t,
\]
which satisfies the equation
\begin{equation}\label{eq:RESCALEDEQUATIONTYPECPRIMECASEN1}
\partial_{\tau}u_{\varepsilon} = \partial_y\left(|\partial_yu_{\varepsilon}^m|^{p-2}\partial_yu_{\varepsilon}^m\right) + f\left(a(a+2\varepsilon)^{-1}u_{\varepsilon}\right), \quad \text{in } \RR^N\times(0,\infty).
\end{equation}
Note that now $f_{\varepsilon}(u_{\varepsilon}) := f\left(a(a+2\varepsilon)^{-1}u_{\varepsilon}\right)$ satisfies $f_{\varepsilon}(a + 2\varepsilon) = f(a) = 0$. Hence, from Theorem \ref{THEOREMEXISTENCEOFTWSPMEREACTIONTYPEC} (part (ii), range $\gamma > 0$), there exists a critical speed $c_{\ast}^{\varepsilon} = c_{\ast}(m,p,\varepsilon) > 0$ and a corresponding $(a+2\varepsilon)$\emph{-admissible} TW with \emph{finite} profile $\varphi_{\varepsilon} = \varphi_{\varepsilon}(\xi)$, and $\xi = x - c_{\ast}^{\varepsilon}t$:
\[
\varphi_{\varepsilon}(-\infty) = a + 2\varepsilon, \quad \varphi_{\varepsilon}(\xi) = 0 \; \text{ for all } \xi \geq \xi_0^{\varepsilon},
\]
for some $\xi_0^{\varepsilon} \in \RR$, satisfying equation \eqref{eq:RESCALEDEQUATIONTYPECPRIMECASEN1}. Thus, if $\widetilde{u}_{\varepsilon} = \widetilde{u}_{\varepsilon}(y,\tau)$ denotes the solution to equation \eqref{eq:RESCALEDEQUATIONTYPECPRIMECASEN1} with $\widetilde{u}_{\varepsilon}(y,0) = u(y,t_{\varepsilon}) \leq a + \varepsilon$ and $\overline{u}_{\varepsilon}(y,\tau) = \varphi_{\varepsilon}(y - c_{\ast}^{\varepsilon}\tau)$, we can repeat the comparison procedure with the assumption $u(x,t_{\varepsilon}) \leq a$, since we can now assume $u(y,t_{\varepsilon}) \leq \varphi_{\varepsilon}(y)$ and so, $\widetilde{u}_{\varepsilon}(y,\tau) \leq \overline{u}_{\varepsilon}(y,\tau)$. We finally obtain $u(x,t) = 0$ in $\{|x| \geq ct\}$ for large times for the arbitrariness of $\varepsilon > 0$.

\emph{Step 3}. In this final step, we prove that for all for all $0 < \varepsilon < a$ and for all $0 < c < c_{\ast}(m,p,f)$, there exists $t_1' = t_1'(\varepsilon,c) > 0$ such that the solution $u = u(x,t)$ satisfies
\[
u(x,t) \geq a -\varepsilon, \quad \text{ in } \{|x| \leq ct\}, \quad \text{ for } t \geq t_1'.
\]
This follows by considering the solution $\underline{u} = \underline{u}(x,t)$ to the problem
\[
\begin{cases}
\partial_t\underline{u} = \partial_x\left(|\partial_x\underline{u}^m|^{p-2}\partial_x\underline{u}^m\right) + f(\underline{u}) \quad &\text{in } \RR\times(0,\infty) \\
\underline{u}(x,0) = \underline{u}_0(x) \quad &\text{in } \RR,
\end{cases}
\]
where $\underline{u}_0 \in \mathcal{C}_c(\RR)$ is defined by $\underline{u}_0(x) := \min\{a,u_0(x)\}$. Consequently, we deduce $\underline{u}(x,t) \leq u(x,t)$ and $0 \leq \underline{u}(x,t) \leq a$
in $\RR\times(0,\infty)$ thanks to the comparison principle, and, furthermore:
\[
\underline{u}(x,t) \geq a-\varepsilon, \quad \text{ in } \{|x| \leq ct\}, \quad \text{ for } t \text{ large enough}.
\]
This last property easily follows by applying Proposition 2.4, Proposition 8.1, Theorem 2.6 of \cite{AA-JLV:art} to $\underline{u} = \underline{u}(x,t)$ and remembering the scaling property quoted at the beginning of this section (we could even repeat the construction done in \cite{AA-JLV:art} using the ``change sign'' TWs introduced in Subsection \ref{FINALREMARKODEPARTCPRIME}. Note that this procedure applies to higher dimensions $N \geq 1$ too, as explained in Theorem 2.6 of \cite{AA-JLV:art}. $\Box$
\paragraph{Remark.} Note that in \emph{Step 3} of the above proof we have applied some results of \cite{AA-JLV:art} even if the reaction $f(\cdot)$ (satisfying \eqref{eq:ASSUMPTIONSONTHEREACTIONTERMTYPEDPME}) is not assumed to be \emph{concave} in $[0,a]$, as it was in the Fisher-KPP setting (cfr. with the assumptions in (1.2) of \cite{AA-JLV:art}). The validity of our procedure follows from the fact that the concavity assumption was used only in Proposition 2.4 of \cite{AA-JLV:art} to prove that solutions (to the Fisher-KPP problem) possess ``minimal non-contracting level sets'' for large times. More precisely: for any $\widetilde{\varrho}_1 > 0$, there exist $\widetilde{\varepsilon}_1 > 0$ and $t_0 > 0$ (depending on the initial data) such that
\begin{equation}\label{eq:REMARKSMALLSUPERLEVELSETSTUYPECPRIME}
u(x,t) \geq \widetilde{\varepsilon}_1 > 0 \quad \text{ in } \{|x| \leq \widetilde{\varrho}_1\} \; \text{ for all } t \geq t_0.
\end{equation}
So, if $f(\cdot)$ satisfies just the assumptions in \eqref{eq:ASSUMPTIONSONTHEREACTIONTERMTYPEDPME}, we can infer as in Remark 3.5 of \cite{C-R2:art}, taking a new reaction $\widetilde{f} = \widetilde{f}(u)$ defined as the primitive of
\[
h(u) := \min_{v \in [0,u]} f'(v),
\]
satisfying $\widetilde{f}(0) = 0$. It easily seen that $\widetilde{f}(\cdot)$ satisfies
\[
\begin{cases}
\widetilde{f}(0) = \widetilde{f}(\theta) = 0, \quad &f(u) \geq \widetilde{f}(u) > 0 \text{ in } (0,\theta) \\
\widetilde{f} \in C^1([0,\theta]), \qquad\qquad\quad\;  &(\widetilde{f})'(0) = f'(0) \\
\widetilde{f}(\cdot) \text{ is concave in } (0,\theta),
\end{cases}
\]
for some $0 < \theta < a$, cfr. with formula (3.20) of \cite{C-R2:art}. Now, since the proof of Proposition 2.4 of \cite{AA-JLV:art} just concerns the ``small'' level sets of $u = u(x,t)$ (i.e. $ \widetilde{\varepsilon}_1 > 0$ can be taken smaller), we can substitute $f(\cdot)$ with $\widetilde{f}(\cdot))$ (which is now concave) and the argue by comparison, since $f \geq \widetilde{f}$ in $(0,\theta)$. Once it is proved the claim in \eqref{eq:REMARKSMALLSUPERLEVELSETSTUYPECPRIME}, the only estimate used in Proposition 8.1 of \cite{AA-JLV:art} is $f(u) \geq q_{\varepsilon}(1-u)$ for all $\varepsilon \leq u \leq 1$, for some $q_{\varepsilon} > 0$ and all $\varepsilon > 0$ fixed. This estimate does not need the concavity assumption as showed in the remark at the end of the proof of Theorem \ref{ASYMPTOTICBEHAVIOURTHEOREMTYPEC}: Part (ii), case $N=1$. Finally, in the proof of Theorem 2.6 of \cite{AA-JLV:art}, the concavity of the reaction term plays no role and so, can be easily adapted to the present setting (cfr. also with the remark at the end of proof of Theorem \ref{ASYMPTOTICBEHAVIOURTHEOREMTYPEC}: Part (ii), case $N \geq 2$). $\Box$
\paragraph{Proof of Theorem \ref{ASYMPTOTICBEHAVIOURTHEOREMTYPECCPRIME}: Case $\boldsymbol{N=1}$, range $\boldsymbol{\gamma = 0}$.} Fix $m > 0$ and $p > 1$ such that $\gamma = 0$. The proof in this range is similar to the previous one, with some modifications.

\emph{Step 0'.} This step coincides with \emph{Step 0} of the range $\gamma > 0$.

\emph{Step 1'.} In this step we proceed as in \emph{Step 1} of the range $\gamma > 0$, considering the super-solution given by the problem
\[
\begin{cases}
\partial_t\overline{u} = \partial_x\left(|\partial_x\overline{u}^m|^{p-2}\partial_x\overline{u}^m\right) + f'(0)\overline{u} \quad &\text{in } \RR\times(0,\infty) \\
\overline{u}(x,0) = u_0(x) \quad &\text{in } \RR,
\end{cases}
\]
and the function $\widetilde{u}(x,t) = e^{-f'(0)t}\overline{u}(x,t)$ satisfying
\[
\begin{cases}
\partial_t\widetilde{u} = \partial_x\left(|\partial_x\widetilde{u}^m|^{p-2}\partial_x\widetilde{u}^m\right) \quad &\text{in } \RR\times(0,\infty) \\
\widetilde{u}(x,0) = u_0(x) \quad &\text{in } \RR.
\end{cases}
\]
This time $\widetilde{u} = \widetilde{u}(x,t)$ does not generally propagate with finite speed of propagation, but it is everywhere positive for all $t > 0$. In the next paragraphs, we provide a bound from above for $\widetilde{u} = \widetilde{u}(x,t)$ which will be enough for our purposes. The main tool are the Barenblatt solutions presented in the introduction (see Subsection \ref{SUBSECTIONPRELIMINARIESTYPECCPRIME}). Indeed, since $u_0$ has compact support, there are a mass $M > 0$ large enough and delay $\theta > 0$ such that $u_0(x) \leq B_M(x,\theta)$ for all $x \in \RR$. Thus, from the Comparison Principle, we obtain $\widetilde{u}(x,t) \leq B_M(x,t+\theta)$ for all $x \in \RR$ and $t > 0$. Coming back to the solution $u = u(x,t)$, this gives
\[
\begin{aligned}
u(x,t) &\leq \overline{u}(x,t) = e^{f'(0)t}\widetilde{u}(x,t) \leq e^{f'(0)t}B_M(x,t + \theta) = e^{f'(0)t} (t+\theta)^{-1/p}F_M(x(t+\theta)^{-1/p}) \\
& = C_M (t+\theta)^{-1/p} \exp\left[f'(0)t - k(t+\theta)^{-1/(p-1)}|x|^{p/(p-1)}\right],
\end{aligned}
\]
where $k = (p-1)p^{-p/(p-1)}$ (cfr. with Subsection \ref{SUBSECTIONPRELIMINARIESTYPECCPRIME}, range $\gamma = 0$). In particular, we obtain
\begin{equation}\label{eq:GLOBALBOUNDFROMABOVETYPECPRIME}
u(x,t) \leq C_{M,\theta}(t) \exp\left[-k_{\theta}(t)|x|^{p/(p-1)}\right], \quad \text{ in } \RR\times(0,\infty),
\end{equation}
where $C_{M,\theta}(t) = C_M (t+\theta)^{-1/p} e^{f'(0)t}$ and $k_{\theta}(t) = k(t+\theta)^{-1/(p-1)}$, $t >0$. Again, this bound can be easily extended to the case $N \geq 1$, with minor changes in the functions $C_{M,\theta}(\cdot)$ and $k_{\theta}(\cdot)$.

\emph{Step 2'.} As before, in this step we prove that for all $\varepsilon > 0$ and for all $c > c_{\ast}(m,p,f)$, there exists $t_1 = t_1(\varepsilon,c) > 0$ such that
\[
u(x,t) \leq \varepsilon, \quad \text{ in } \{|x| \geq ct\}, \quad \text{ for all } t \geq t_1.
\]
So, we fix $\varepsilon > 0$, $c > c_{\ast}(m,p,f)$, and we consider $t_{\varepsilon} > 0$ given by \emph{Step 1'-Step 1}. As before, we can assume $u(x,t_{\varepsilon}) < a$, since the scaling technique exploited in \emph{Step 2} of the range $\gamma > 0$ works also in the present setting. Again, we consider
\[
\overline{v}(x,t) := \varphi(x - c_{\ast}t), \qquad\qquad \overline{w}(x,t) := \psi(x + c_{\ast}t),
\]
where $\varphi = \varphi(\xi)$ is the \emph{positive $a$-admissible} TW studied in Theorem \ref{THEOREMEXISTENCEOFTWSPMEREACTIONTYPEC} (part (ii), range $\gamma = 0$), satisfying $\varphi(-\infty) = a$, $\varphi(+\infty) = 0$ , and $\psi = \psi(\xi)$ is its ``reflection''. The main difference w.r.t. the range $\gamma > 0$ is neither $u = u(x,t)$ nor $\overline{v} = \overline{v}(x,t)$ (resp. $\overline{w} = \overline{w}(x,t)$) have compact support in $\RR$, and so we cannot immediately conclude $u(x,t_{\varepsilon}) \leq \overline{v}(x,0)$ (resp. $u(x,t_{\varepsilon}) \leq \overline{w}(x,0)$) in $\RR$ (up to a right/left shift) of the profile $\varphi = \varphi(x)$ (resp. $\psi = \psi(x)$).

\noindent However, we known that the asymptotic behaviour of the tails of $\varphi = \varphi(x)$ and $\psi = \psi(x)$ (cfr. with formula \eqref{eq:ASYMPTOTICSOFCRITICALPROFILEGAMMA0}):
\[
\varphi(x) \sim a_0 |x|^{\frac{2}{p}}e^{-\frac{\lambda_{\ast}}{m}x}, \quad \text{for } x \sim +\infty,
\qquad
\psi(x) \sim a_0 |x|^{\frac{2}{p}}e^{-\frac{\lambda_{\ast}}{m}|x|}, \quad \text{for } x \sim -\infty,
\]
where $\lambda_{\ast}:=(c_{\ast}/p)^m$, $a_0 > 0$, and, at the same time,
\[
u(x,t_{\varepsilon}) \leq C_{M,\theta}(t_{\varepsilon}) \exp\left[-k_{\theta}(t_{\varepsilon})|x|^{p/(p-1)}\right], \quad x \in \RR\times(0,\infty),
\]
from the global bound \eqref{eq:GLOBALBOUNDFROMABOVETYPECPRIME} of \emph{Step 1'}. Consequently, since $p > 1$, $u(x,t_{\varepsilon})$ decays faster than $\varphi(x)$ and $\psi(x)$ when $|x| \sim \infty$, and so we can now assume $u(x,t_{\varepsilon}) \leq \overline{v}(x,0)$ and $u(x,t_{\varepsilon}) \leq \overline{w}(x,0)$ for all $x \in \RR$ and applying the Comparison Principle to have $u(x,t + t_{\varepsilon}) \leq \overline{v}(x,t)$ and $u(x,t + t_{\varepsilon}) \leq \overline{w}(x,t)$ for all $x \in \RR$ and $t > 0$. Thus, using that $\overline{v}(x,t) \leq \varepsilon$ for $x \geq c_{\ast}t + \xi_0$ and $\overline{w}(x,t) \leq \varepsilon$ for $x \leq -c_{\ast}t + \xi_0$ and $c > c_{\ast}$, we deduce that $u(x,t) \leq \varepsilon$ in $\{|x| \geq ct\}$ for large times (see also the proof of Theorem \ref{ASYMPTOTICBEHAVIOURTHEOREMTYPEC}, Part (iii)). $\Box$
\paragraph{Proof of Theorem \ref{ASYMPTOTICBEHAVIOURTHEOREMTYPECCPRIME}: Case $\boldsymbol{N \geq 2}$.} Fix $m > 0$ and $p > 1$ such that $\gamma > 0$ (the range $\gamma = 0$ is almost identical and we skip it). Again, we focus on radial solutions to problem  \eqref{eq:ALLENCAHNPME}, i.e., solutions $u = u(r,t)$ to  problem \eqref{eq:RADIALPROBLEMTYPECPRIME}:
\[
\begin{cases}
\partial_t u = \partial_r\left(|\partial_ru^m|^{p-2}\partial_ru^m\right) + \frac{N-1}{r}|\partial_ru^m|^{p-2}\partial_ru^m + f(u) \quad &\text{in } \RR_+\times(0,\infty) \\
u(r,0) = u_0(r) \quad &\text{in } \RR_+\times\{0\},
\end{cases}
\]
where $r = |x|$, $x \in \RR^N$, and $u_0(\cdot)$ is a radial decreasing initial datum.

\emph{Step 1: Convergence to zero in ``outer'' sets.} Proceeding as in the proof of Theorem \ref{ASYMPTOTICBEHAVIOURTHEOREMTYPEC} (Part (i)), we can assume $\partial_ru^m \leq 0$ in $\RR_+\times(0,\infty)$. Consequently, the solution $\overline{u} = \overline{u}(r,t)$ to the problem
\[
\begin{cases}
\partial_t \overline{u} = \partial_r\left(|\partial_r\overline{u}^m|^{p-2}\partial_r\overline{u}^m\right) + f(\overline{u}) \quad &\text{in } \RR_+\times(0,\infty) \\
\overline{u} = u \quad &\text{in } \{0\}\times(0,\infty) \\
\overline{u}(r,0) = u_0(r) \quad &\text{in } \RR_+\times\{0\},
\end{cases}
\]
is a super-solution to \eqref{eq:RADIALPROBLEMTYPECPRIME} and, at the same time, it is a solution of the one-dimensional equation with compactly supported initial data. Thus, for all $c > c_{\ast}(m,p,f)$, it follows
\[
\overline{u}(r,t) = 0 \text{ uniformly in } \{r \geq ct\}, \quad \text{ as } t \to +\infty,
\]
and by the comparison, we deduce the same for $u = u(r,t)$. Of course, if $\gamma = 0$, the solutions are always positive and it holds $u(r,t) \leq \varepsilon$ uniformly in $\{r \geq ct\}$ for large times $t > 0$.

We ask the reader to note that with the same comparison technique we can prove that
for all $\varepsilon > 0$, it holds
\begin{equation}\label{eq:NDIMPROPOPERTYSMLLTHENAPLUSVEP}
u(x,t) \leq a + \varepsilon, \quad \text{for all } x \in \RR^N, \quad t \geq t_{\varepsilon},
\end{equation}
for some suitable waiting time $t_{\varepsilon} > 0$ (as we have seen before, this property holds for the case $N=1$).

\emph{Step 2: Convergence to $a$ in ``inner'' sets.} In this second step, we have to prove that for all $\varepsilon > 0$ and $0 < c < c_{\ast}(m,p,f)$, the solution to problem \eqref{eq:RADIALPROBLEMTYPECPRIME} satisfies
\begin{equation}\label{eq:CONVERGENCETOAINNERSETS}
u(r,t) \geq a -\varepsilon, \text{ uniformly in } \{r \leq ct\}, \quad t \to +\infty.
\end{equation}
Following \cite{AA-JLV:art}, we have to proceed in three main steps. In the first one we have to show that the solution $u = u(r,t)$ does not extinguish and actually lifts-up to a small level $\widetilde{\varepsilon} > 0$ in compact sets of $\RR^N$ for large times. This follows from Proposition 2.4 and 2.5 of \cite{AA-JLV:art} of the Fisher-KPP setting and recalling the scaling property linking reactions of type C' to Fisher-KPP reactions.
\\
Then following the proof of Theorem \ref{ASYMPTOTICBEHAVIOURTHEOREMTYPEC}, Part (ii) case $N=1$, we have
\[
f(u) \geq q(a-u) \quad \text{in } \{|x| \leq \widetilde{\varrho}\}\times[t_{\widetilde{\varrho}},\infty),
\]
for a suitable choice of $q = q_{\widetilde{\varepsilon}} > 0$, all $\widetilde{\varrho} > 0$ and $t_{\widetilde{\varrho}} > 0$ large enough. We point out that the previous inequality holds true only when $\widetilde{\varepsilon} \leq u \leq a$, which is an assumption we can make thanks to \eqref{eq:NDIMPROPOPERTYSMLLTHENAPLUSVEP} and the scaling technique employed in \emph{Step 2} of proof of the case $N=1$ (see range $\gamma > 0$). Thus, exactly as before, we get that for all $\varepsilon > 0$ (small) and $\widetilde{\varrho} > 0$ (large)
\[
u(r,t) \geq a - \varepsilon \quad \text{ in } \{r = |x| \leq \widetilde{\varrho}\} \; \text{ for all } t \geq t_1,
\]
for some (large) $t_1 > 0$. Finally, we get \eqref{eq:CONVERGENCETOAINNERSETS} by constructing a sub-solution to problem \eqref{eq:RADIALPROBLEMTYPECPRIME} through ``change sign'' TWs (cfr. with the proof of Theorem 2.6 of \cite{AA-JLV:art} and Subsection \ref{FINALREMARKODEPARTCPRIME}). Recalling the scaling property linking reactions of type C' to Fisher-KPP reactions, we consider a barrier (from below) built with the function
\[
\underline{\varphi}(\xi) :=
\begin{cases}
a-\varepsilon              \quad & \text{if } \xi \leq 0 \\
\varphi_c(\xi) \quad & \text{if } 0 \leq \xi \leq \xi_1^c  \\
0              \quad &  \text{otherwise}
\end{cases}
\qquad \text{ with } \quad a-\varepsilon = \max \varphi_c(\xi),
\]
where $\varphi_c = \varphi_c(x-(c+\varepsilon)t)$ is a ``change-sign'' TW (of type 2) corresponding to the speed $0 < c < c_{\ast}$ (see Subsection \ref{FINALREMARKODEPARTCPRIME}). Thus, the barrier propagate level $a-\varepsilon$ with speed $c$, and so, using the arbitrariness of $0 < c < c_{\ast}$ obtain \eqref{eq:CONVERGENCETOAINNERSETS} (cfr. with the proof of Theorem 2.6 of \cite{AA-JLV:art} for all the details). $\Box$
%
%
%
%
%
%
%
%
%
%
%
\section{Comments and open problems}\label{SECTIONCOMMENTANDOPENPROBLEMSTYPECCPRIME}
We end the paper with some comments and open problems.
\paragraph{``Fast'' diffusion range.} A first possible extension of our work consists in studying problem \eqref{eq:ALLENCAHNPME} with a different assumption on the parameters $m > 0$ and $p > 1$, i.e.
\[
-p/N < \gamma < 0,
\]
see the final sections of \cite{V1:book} for more information on this range, also called ``fast'' diffusion range. In the Fisher-KPP setting it has been done in \cite{AA-JLV:art1} where it has been proved the solutions spread exponentially fast for large times, in sharp contrast w.r.t. the ``slow'' diffusion range.
\paragraph{An interesting limit case.} As already pointed out in \cite{AA-JLV:art}, keeping $m(p-1) = \theta$ with $\theta > 0$ fixed, we can formally compute the limit of the doubly nonlinear operator:
\[
\Delta_pu^m = m^{\frac{\theta}{m}}\nabla\cdot\left(u^{\theta+1-p}|\nabla u|^{p-2}\nabla u\right) \to \nabla\cdot\left(u^{\theta}\frac{\nabla u}{|\nabla u|}\right)  \quad \text{as } m \to \infty \text{ and } p \to 1.
\]
Consequently, a very interesting open problem is the study of the existence of admissible TWs for the equation
\[
\partial_t u = \partial_x\left(u^{\theta}|\partial_x u|^{-1}\partial_x u\right) + f(u) \quad \text{in } \RR\times(0,\infty),
\]
for different values of the parameter $\theta \geq 0$, and $f(\cdot)$ of type C or C'. In the Fisher-KPP setting there are interesting results in the series of papers \cite{A-Cas-M:art,C-C-Cas-S-S:art, C-G-S-Sol:art, Cam-Sol:art}, where the authors showed the existence of discontinuous TWs which are a very interesting novelty w.r.t. to the doubly nonlinear diffusion. However, the flux-limited operators they consider do not cover this new class introduced above.
\paragraph{Reactions of type B.} In the literature, other kind of reactions have been intensively investigated. Possibly, the most famous are the so called reactions of type B
\begin{equation}\label{eq:ASSUMPTIONSONTHEREACTIONTERMTYPECOMBUSTION}
\begin{cases}
f(0) = f(a) = f(1) = 0, \quad f(u) \leq 0 \text{ in } (0,a), \;\; f(u) > 0 \text{ in } (a,1) \\
f \in C^1([0,1]), \qquad\qquad\quad\;  f'(1) < 0,
\end{cases}
\end{equation}
which emerge from combustion models (see for instance the famous works \cite{BerestNiren1992:art,Roquejoffre1997:art} and the interesting survey \cite{Sire2014:art}). We have not considered this framework in this paper but we want to point out, thanks to a simple comparison with reaction of type C, that part (ii) of Theorem \ref{ASYMPTOTICBEHAVIOURTHEOREMTYPEC} hold even for reaction of type B. Also part (i) holds if we take initial data $0 \leq u_0 \leq a$ (this is true even for reactions of type C thanks to a straightforward  comparison technique, but we have not insisted on it since it goes out of our purposes).
\paragraph{Sharp threshold results.} As we have pointed out in the presentation of the results of this paper, Theorem \ref{ASYMPTOTICBEHAVIOURTHEOREMTYPEC} has not a sharp threshold statement. As already explained, the problem has been studied and solved in dimension $N=1$ and very general reaction terms by Du and Matano \cite{DuMatano2010:art}. In this work, it is essential the existence of nontrivial solutions to
\[
-\partial_{xx}u = f(u) \quad \text{ in } \RR,
\]
which correspond to stationary solutions to the corresponding parabolic problem and eventual limit configurations (see for instance Theorem 1.1 of \cite{DuMatano2010:art}). The study of these stationary solutions is clearly more complicated in the doubly nonlinear framework and seems to be a very challenging open problem (see \cite{MuratovZhong2017:art,Polacik2011:art} for the case $N \geq 1$).

%
%
%
%
%
%
%
%
%
%
%

\bigskip

%
%

\noindent {\textbf{\large \sc Acknowledgments.}} The author has been partially funded by Projects MTM2011-24696 and MTM2014-52240-P (Spain), by local projects ``Equazioni differenziali lineari e non lineari'', ``Equazioni differenziali non lineari e applicazioni'' (Italy), by the GNAMPA project ``Equazioni diffusive non-lineari in contesti non-Euclidei e disuguaglianze funzionali associate'' (Italy) and by the ERC Advanced Grant 2013 n.~339958 ``Complex Patterns for Strongly Interacting Dynamical Systems - COMPAT''.

\noindent I want to thank Juan Luis V\'azquez for having introduced me to reaction diffusion equations with nonlinear diffusion, and for his constant support to my research activities. Thanks to Nicola Vassena too for sharing his ideas during a really interesting discussion.

%
%
%
\vskip 1cm

\

2000 \textit{Mathematics Subject Classification.}
35K57,  
35K65, 
35C07,  	
35K55. 

\medskip

\textit{Keywords and phrases.}  Bistable equations, Doubly nonlinear diffusion, Propagation of level sets.

\bigskip

\bigskip

Contact information:

\smallskip

\begin{tabular}{ll}

Address:   & Dipartimento di Matematica ``Giuseppe Lugi Lagrange'', Politecnico di Torino. \\
           & Corso Duca degli Abruzzi, 24, 10129, Torino, Italy. \\
E-mail(s): & alessandro.audrito@polito.it \\

\end{tabular}

\end{document}